\documentclass[a4paper,11pt]{article}
\usepackage[margin=3cm]{geometry}
\usepackage{bm}
\usepackage{amsfonts,amssymb,amsmath,amsthm}
\usepackage{stmaryrd} 
\usepackage{latexsym}
\usepackage{graphics}
\usepackage{epsfig}
\usepackage{psfrag}
\usepackage{bbm}
\usepackage{soul} 
\usepackage{enumerate}
\usepackage{bbm}
\usepackage{nicefrac}
\usepackage{hyperref}
\DeclareFontEncoding{FMS}{}{}
\DeclareFontSubstitution{FMS}{futm}{m}{n}
\DeclareFontEncoding{FMX}{}{}
\DeclareFontSubstitution{FMX}{futm}{m}{n}
\DeclareSymbolFont{fouriersymbols}{FMS}{futm}{m}{n}
\DeclareSymbolFont{fourierlargesymbols}{FMX}{futm}{m}{n}
\DeclareMathDelimiter{\VERT}{\mathord}{fouriersymbols}{152}{fourierlargesymbols}{147}

\def\E{\mathbb{E}}

\def\1{\mathbbm{1}}

\usepackage{todonotes}
\presetkeys{todonotes}{size=\scriptsize}{}
\newcommand\denis[1]{\todo[color=green!40]{1}}
\newcommand\denisIL[1]{{\todo[inline,color=green!40]{1}}}
\newcommand\nicolas[1]{\todo[color=blue!40]{1}}
\newcommand\nicolasIL[1]{{\todo[inline,color=blue!40]{1}}}
\newcommand\edouard[1]{\todo[color=red!40]{1}}
\newcommand\edouardIL[1]{{\todo[inline,color=red!40]{1}}}

\newtheorem{theoreme}{Theorem}[section]
\newtheorem{lemme}[theoreme]{Lemma}
\newtheorem{proposition}[theoreme]{Proposition}

\newtheorem{definition}[theoreme]{Definition\rm}

\newtheorem{remarque}{\bf Remark}

\usepackage{epsf}
\usepackage{enumitem}
\usepackage[english]{babel}
\usepackage{graphicx}
\usepackage{mathptmx}
\usepackage{makeidx}
\usepackage{amsmath}
\usepackage{amsfonts}
\usepackage[bottom]{footmisc}
\usepackage{multicol}
\usepackage{type1cm}
\usepackage{graphicx}
\usepackage{url}
\usepackage{color}
\usepackage[pagewise]{lineno}
\usepackage[francais]{minitoc}
\usepackage{setspace}
\usepackage{float}
\usepackage{bm}
\usepackage{authblk}
\usepackage[T1]{fontenc}
\usepackage[applemac]{inputenc}
\title{Untangling the role of temporal and spatial variations in persistence of populations}

\author[1]{Michel Bena{\"i}m }

\author[2]{Claude Lobry}

\author[3]{Tewfik Sari}

\author[4]{\'Edouard Strickler}

\affil[1]{ Institut de Math\'ematiques, Universit\'e de Neuch{\^a}tel, Switzerland }
\affil[2]{C.R.H.I, Universit\'e Nice Sophia Antipolis, France}
\affil[3]{ITAP, University of Montpellier, INRAE, Institut Agro, Montpellier, France}
\affil[4]{ Universit\'e de Lorraine, CNRS, Inria, IECL, Nancy, France }
\date{\today}

\newcommand{\og}{\guillemotleft}
\newcommand{\fg}{\guillemotright}
\newcommand{\po}{\mathrm{o}}
\newcommand{\disp}{\displaystyle}
\newcommand{\bitbul}{\begin{itemize}[label = \textbullet]}
\newcommand{\bittiret}{\begin{itemize}[label = -]}
\newcommand{\bito}{\begin{itemize}[label =$\circ$]}
\newcommand{\bit}{\begin{itemize}}
\newcommand{\fit}{\end{itemize}}
\newcommand{\ben}{\begin{enumerate}}
\newcommand{\fen}{\end{enumerate}}

\newcommand{\fin}{\end{document}}
\newcommand{\beq}{\begin{equation}}
\newcommand{\feq}{\end{equation}}
\newcommand{\dcom}{\begin{quote}\begin{small}}
\newcommand{\fcom}{\end{small}\end{quote}}

\newcommand{\bc}{\begin{center}}
\newcommand{\fc}{\end{center}}
\newcommand{\emat}{\mathrm{e}}
\newcommand{\ch}{\mathrm{cosh}}
\newcommand{\sh}{\mathrm{sinh}}
\newcommand{\eps}{\varepsilon}
\def\1{{\rm 1\mskip-4.4mu l}}
\newcommand\mb{\mathbf}

\begin{document}
\maketitle

\abstract{We consider a population distributed between two habitats, in each of which it experiences a growth rate that switches periodically between two values, $1- \varepsilon > 0$  or $ - (1 + \varepsilon) < 0$. We study the specific case where the growth rate is positive in one habitat and negative in the other one for the first half of the period, and conversely for the second half of the period, that we refer as the $(±1)$ model. In the absence of migration, the population goes to $0$ exponentially fast in each environment. In this paper, we show that, when the period is sufficiently large, a small dispersal between the two patches is able to produce a very high positive exponential growth rate for the whole population, a phenomena called inflation. We prove in particular that the threshold of the dispersal rate at which the inflation appears is exponentially small with the period. We show that inflation is robust to random perturbation, by considering a model where the values of the growth rate in each patch are switched at random times: we prove that inflation occurs for low switching rate and small dispersal. We also consider another stochastic model, where after each period of time $T$, the values of the growth rates in each patch is chosen randomly, independently from the other patch and from the past. Finally, we provide some extensions to more complicated models, especially epidemiological and density dependent models.}


\tableofcontents
\section{Introduction}
It is a ubiquitous fact that  a population has the ability to migrate between several patches  which  have different environmental conditions. A patch is called a {\em source }  when, in the absence of migration, the environmental conditions lead to the persistence of the population, and a {\em sink} when, on the contrary, they lead to the extinction of the population. A question of primary importance is obviously to understand how environmental conditions and migration interact so that a set of patches is or is not globally favorable to persistence. Mathematical modeling by dynamical systems is one of the tools used to address this question and the papers that use it are innumerable, so we give up reporting on them here. We present only those that we consider important for our purpose, which is to investigate the conditions under which migration between two patches can increase or decrease the abundance of the metapopulation.

The simplest case of a continuous time model,  two patches with logistic dynamics and linear migration, has been extensively studied in the case of a fixed environment (i.e. the parameters of the model do not depend on time) \cite{ ARD15,ARD18, DeA79, DeA14, FRI77}. This very elementary (and therefore unrealistic) model is now well understood mathematically and it appears that for certain values of the growth parameters the total population at equilibrium is not a monotonic function of the migration intensity \cite{ ARD15,ARD18, DeA14}, a phenomenon that we will find again in the case of variable environments that we study here.

In the case of migration between a {\em source}  and a {\em sink}  it is intuitively clear that migration from the  {\em source}  to the {\em sink} can prevent the extinction of the population on the latter. On the other hand, it seems paradoxical that: 
\dcom
Populations can persist in an environment consisting of sink habitats only.
\fcom
as announced in the title of the article \cite{JAN98} by Jansen et al.
Our article is a contribution to the clarification of this paradox. 

Jansen et al. consider the implicitly spatial  discrete-time model:
\beq \label{discrete1} 
N(t+1) =[m f S_1(t) + m (1-f) S_2(t)]N(t).
\feq
which represents the fact that each individual give $m$ offspring, that are then dispersed according to the fractions $f$  and $1-f$ on sites 1 and 2, where they survive at rate $S_1(t)$ and $S_2(t)$, respectively. Jansen et al. assume that:
\bito
\item $S_1$ is a sequence of independent random variables taking two values $S_a < 1/m < S_b$ with probabilities $p$ (bad years) and $1-p$ (good years), the parameters being such that on the long term the patch 1  is a {\em sink}.
\item $S_2$ is {\em constant} and strictly smaller than $1/m$, so that the second patch is also a {\em sink} 
\fit
Intuitively persistence is possible with a little dispersion on patch 2 which means building up reserves for bad years.  This is indeed what Jansen et al. show by calculation: for values of $m$ that are neither too large nor too small, the meta-population is persistent. 

{\em Migration can therefore have an "inflationary" effect}, an expression coined by Holt in \cite{HOLT97}.
That  ''inflationary'' effect  noticed by Holt was of another nature.  In \cite{HOLT97} the author also considers a model of the form \eqref{discrete1} and assumes that $S_1$ is deterministic but density dependent $S_1 = S_1(N(t))$ (the notations have been changed) as, for example, in the logistic or Ricker models. Based on the classical results of May and Oster (see \cite{MAY76}) on the appearance of periodic and then chaotic solutions in discrete density-dependent dynamics, Holt remarks that when the population has a stable equilibrium in the absence of migration, the presence of migration to a  {\em sink} only decreases the value of this equilibrium, but, on the other hand, if the population has periodic solutions, the migration to a {\em sink}  can significantly increase the mean of the metapopulation in the long run. Since variations in the population $N(t)$ in the density-dependent model $N(t+1) = S_1(N(t),t)N(t)$ can be interpreted (if denoting $R_1(t) = S_1(N(t))$) as fluctuations in the fitness of the linear model $N(t+1) = R(t)N(t)$, Holt concludes that the presence of autocorrelation in the variations in the sequence of replacement rates $R(t)$ can be a cause of inflation.

Gonzalez and Holt in \cite{HOLT02} have also highlighted an ''inflationary'' effect of the environment fluctuations in the case of the continuous time model:
$\frac{dN}{dt} = f(t)N(t) +I$. Here it is assumed that in the absence of immigration $I(t)$ the population $N(t)$ tends towards extinction and that its persistence is ensured only in the presence of immigration $I$. When the growth rate is fixed and negative, $f(t) \equiv -\mu$ the population tends towards a stationary population $N^* = \frac{I}{\mu}$. The objective of the authors is to compare this equilibrium with the mean $\displaystyle \overline{N} = \lim_{t \rightarrow +\infty} \frac{1}{T}\int_0^TN(t)dt$ of the population when $f(t)$ is no longer constant and to show that, in a certain sense, $\overline{N}$ exceeds $N^*$ all the more as the fluctuations of $f$ are important. To do this their strategy is to consider piecewise constant periodic functions for which they can make explicit calculations that describe the inflation phenomenon; then, as they notice that  \og square-wave deterministic variation is of course a rather artificial pattern of temporal variation \fg$\,$ they turn to more realistic models on which they show by simulation the existence of inflation.

These three studies and others (see \cite{SCH10} for a more detailled discussion of this topic) where the spatialization is ''implicit'' have in common, whether they are discrete or continuous in time, deterministic or random, to use {\em one dimensional} models  where the mathematical properties are easier to determine. The next step is to consider an "explicit" spatialization with two or more sites. This is what Roy, Holt and Barfield do in \cite{ROY05} where they consider the probabilistic  model on $n $ sites:
\beq \label{Roy}
N_i(t+1) = R_i(t)N_i(t) + I_i(t)-E_i(t)
\feq
where $I_i$ and $E_i$ represent respectively the immigration and emigration on the site.  They demonstrate the inflation effect through unformal reasonings and numerous simulations. They can conclude their discussion with: 
\dcom
Given temporal variability and positive temporal autocorrelation in local growth rates, moderate rates of dispersal can enhance the ability of a sink metapopulation to persist; moreover, given persistence, temporal autocorrelation can inflate metapopulation abundance.
\fcom
 This pioneering work of Holt and his colleagues has been further refined and mathematically clarified  by  Schreiber.
In \cite{SCH10}, Schreiber considers the model \eqref{Roy} of Roy et al. in the form:
\beq \label{SCH}
 N^i_{t+1} = \left( 1 - \sum_{k≠i}d_{ki}\right) f_t^iN_t^i+ \sum_{j≠i} d_{ij}f^j_t N^j_t = \sum_{j=1}^nd_{ij}f^j_t N^j_t
 \feq 
From precise mathematical developments based on a probabilistic version of the Perron-Frobenius theorem which allows to show that for this type of model the concept of growth rate of the meta population is well defined he can conclude:
\dcom
When environmental fluctuations have positive temporal autocorrelations and the population is partially
mixing, the metapopulation growth rate can be positive despite the arithmetic mean of fitness being less than 1
in every patch. (...)Furthermore, in the presence of these positive autocorrelations, the analysis reveals
that the maximal metapopulation growth rate occurs at intermediate dispersal rates,(...)
\fcom 
Unlike the case of discrete-time stochastic models that we have just briefly examined, the case of continuous-time deterministic models has been little studied.
 With the exception of \cite{HOLT03} which treats one site only, \cite{KLA08} which discusses the case of two sites on an example and \cite{EVA13} which treats the case where the growths on each site obey diffusion processes, we only know of a very recent\footnote{Which was published while our paper was under review.} article \cite{KAT22} by Katriel: {\em Dispersal-induced growth in a time-periodic environment} which treats the question of inflation on continuous time models.  Katriel considers the following model (the notations are modified to remain consistent with the previous notations): 
\beq \label{Katriel}
\displaystyle \frac{dN_i}{dt} = r_i(t)N_i + m\sum_{j≠i } L_{ij}(N_j-N_i)\quad
\feq
where the functions $r_i(t)$ are continuous periodic functions of period $T$ that define the growth rates on $n$ isolated sites, the $L_{ji} = L_{ij} \geq 0$ describe the geometry between the sites (assumed to be connected) and $m$ measures the strength of the migration. As in the discrete case, extended versions of the Perron-Frobenius theorem allows to define the metapopulation growth  and relying on a theorem of  \cite{LIU22} 
he shows that for inflation to occur it is necessary that the $r_i(t)$ are desynchronized,
that the period $T$ is large enough and that the migration intensity $m$ is large enough, but not too large.
This result is important because it concerns the case of any number of sites.

In turn, we consider Katriel's \eqref{SCH} model in the case of {\em only two sites}, which is obviously less general, but for more general, deterministic discontinuous and then random functions $r_i(t)$ and not necessarily  symmetric  migrations.

We start from a remarkably simple particular case which we call the periodic $(±1)$-model (where $r_i(t)$ are piecewise constant  equals to $+1-\eps$ or $-1-\eps$). This model depends on three parameters $\eps$, which represents the decay rate of each {\em sink}, $m$, which represents the migration rate and the period $T$.  We explicitly compute the growth rate of the metapopulation, and give explicit bounds on  $m$ and $T$ for inflation to occur.
In particular we show that that the threshold $m^*$ for the appearance of inflation can be very small, precisely exponentially small with respect to the period. This last point may have a practical importance. Indeed, in a recent paper  \cite{HOLT20}, Kortessisa et al. (Holt is the last co-author), about the COVID pandemic, based on simulations, draw the attention to the possibility of inflation in case of migration between two patches when the sanitary policies are not synchronous; our result confirms analytically these simulations and indicates that the migration threshold from which the pandemic develops, can be very small.

After having analyzed in detail the mathematical properties of the periodic $(±1)$-model we extend them to more realistic models, and, most importantly, we show that periodicity is not essential in the following sense. Our $(±1$) deterministic model can be understood as a situation where two different environmental regimes $+1$ and $-1$ follow each other for a fixed duration $T$. We  study what happens when the regimes succeed each other over random durations or with random  growth rates $r_i$, and show on various probabilistic models how the results obtained for the periodic $(±1)$ model can be extended.
These last results cannot be deduced from \cite{KAT22} which proposes a purely deterministic framework.

Mathematically we make essential use of the following change of variable. If $N_1(t)$ and $N_2(t)$ are the abundances at each patch we pose 
$$ U =\frac{ \log(N_1)+\log(N_2)}{2}= \log(\sqrt{N_1N_2})$$
$$ V= \frac{ \log(N_2)-\log(N_1)}{2}= \log(\sqrt{N_1/N_2})$$
Thus $U$ is the logarithm of the geometric mean of the abundances and $|V|$ is the logarithm of their geometric standard deviation\footnote{see e.g. \url{https://en.wikipedia.org/wiki/Geometric_standard_deviation}}.
This change of variable has the merit of translating the analytical properties of the model($±1$) into  "visible" and "robust" geometrical properties observed in the $(U,V)$ plane.

The paper is organized as follows. Section \ref{plusoumoins} contrains a detailed mathematical treatment of the $(±1)$-model for both deterministic and stochastic environment. Section \ref{morecomplex} contains extensions to more realistic models; the mathematical treatment is less precise and sometimes just outlined. The section \ref{discussion} contains an attempt to give an informal explanation of our view of the inflation phenomenon in the case of continuous time models and its relation to previous work on discrete time models.  An appendix contains technical details regarding mathematical proofs. A symbol index is provided at the end of the appendix.

\section{The ($±1$)-model}\label{plusoumoins}
\subsection{Some results of G. Katriel}
In \cite{KAT22}, Katriel considers the model :
\beq \label{KAT1}
\frac{dx_i}{dt} = r_i(\omega t)x_i +m \sum_{i≠j}L_{ij}(x_j-x_i ),\quad \quad i = 1,\cdots,N\quad L_{ij} \geq 0,\quad m\geq 0
\feq
where the function $r_i$ is $2 \pi$ - periodic {\em continuous} and  represents the varying growth rate within patch $i$, while, for $i \neq j$, $L_{ij} = L_{ji}$ represents the ability of dispersal from patch $i$ to patch $j$. We denote $L_{ii} = - \sum_{j \neq i} L_{ij}$ and assume that the matrix $L = (L_{ij})$ is \emph{irreducible}, meaning that the population can spread in all patches. With vector notations we can write \eqref{KAT1} in the form :
\beq \label{KAT1}
\frac{dX}{dt} = \big(R(\omega t) +mL\big)X
\feq
where $R$ is the diagonal matrix whose diagonal elements are $r_i(\omega t)$. We summarize here a part of the results of \cite{KAT22}.
We let
\bito
\item $\displaystyle \overline{r_i}=\frac{1}{2\pi}\int_0^{2\pi}r_i(\theta)d\theta $,
\item $\displaystyle r_{\mathrm{max}}(\theta) = \max_{i=1..N} r_i(\theta)$, and
\item $\displaystyle \chi = \frac{1}{2\pi}\int_0^{2\pi}r_{\mathrm{max}}(\theta)d\theta$.
\item The growth rate of a positive function $t \mapsto x(t)$ is the limit, if it exists,
$$ \Lambda[x] = \lim_{t \mapsto \infty} \frac{1}{t} \ln(x(s)).$$
\fit
When $m=0$, all equations of \eqref{KAT1} are decoupled,
$$\Lambda[x_i] = \overline{r_i},$$
and depending on whether $ \overline{r_i}$ is positive or negative we will say that the patch $i$ is a ''source'' or a ''sink''.
As soon as $m$ is strictly positive, since  $L$ is irreducible, from Perron-Frobenius theory, it follows that all the growth rates $ \Lambda[x_i]$ are equal. This common rate is the growth rate of the metapopulation and is noted $\Lambda(m,\omega)$.
\begin{definition} \label{defKat} \cite{KAT1} One says  that there is ''Dispersal Induced Growth'' (DIG) if, while all $\overline{r_i} $ are strictly negative, we have $\Lambda(m,\omega) > 0$.
\end{definition}
\textbf{Comment} Note that DIG is the phenomenon which is called {\em inflation} in \cite{EVA13, HOLT02, HOLT03, HOLTPNAS20, KLA08, ROY05, SCH10} with more or less formalized definitions associated to each context. In the present paper when we say ''inflation'' we refer to the above formalized   definition.

The problem is to characterize the set of values of $m$ and $\omega$ for which there is  DIG ({\em inflation}).

\begin{figure}
  \begin{center}
 \includegraphics[width=0.5\textwidth]{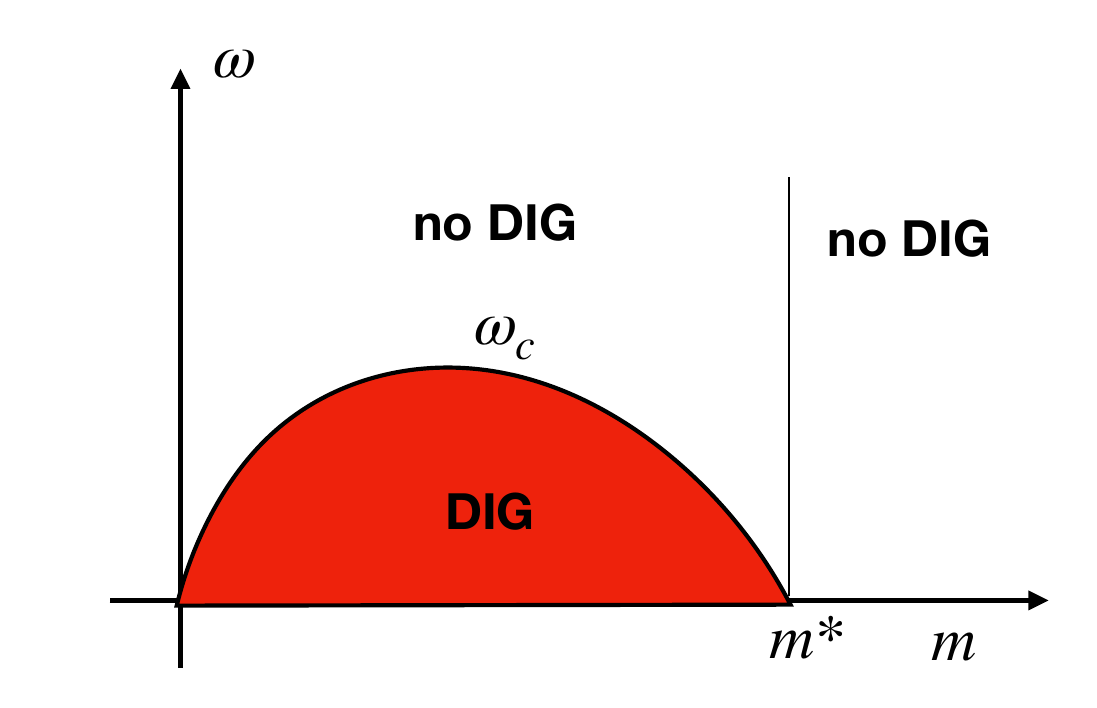}
 \caption{Inflation when $\chi > 0$}\label{ThKat}
 \end{center}
 \end{figure}

We note:
\beq \label{lambdazero}
 \Lambda_0(m) = \frac{1}{2\pi} \int_0^{2\pi} \lambda(R(\theta)+mL) d\theta 
\feq
where $\lambda(R(\theta)+mL) $ is the dominant eigenvalue of the matrix $(R(\theta)+mL)$.
\begin{theoreme}\label{theoremKatriel} \cite{KAT1} (see Figure  \ref{ThKat}).
\bittiret
	\item If $\chi < 0$, for all $m$ and all $\omega$ we have $\Lambda(m,\omega) < 0$. There is no DIG.
	\item If $\chi >0$ the equation $ \Lambda_0(m) = 0,\; m>0$ has a unique solution $m^*$ and there exists a  function $\omega_c$ continuous on $[0, 	m^*]$, strictly positive on $]0, m^*[$ , such that $\omega_c(0) = \omega_c(m^*) = 0$ for which :
	\ben
		\item If $m < m^*$,
		$$\omega < \omega_c \Longrightarrow \Lambda(m,\omega) > 0 \;\; \mathrm{there\;is\;DIG}\quad $$
		$$\omega > \omega_c \Longrightarrow \Lambda(m,\omega) < 0 \;\;\mathrm{ there\; is \;no\; DIG} $$
		\item If $m> m^*$,
		$$ \Lambda(m,\omega) < 0 \;\; \mathrm{there \;is\;no\;DIG }$$ 
	\fen
\fit
\end{theoreme}

In general it is not possible to compute effectively $\Lambda_0(m)$ except in the case of two sites. In \cite{KAT22} the following formula for $ \Lambda_0(m) $ is given :
$$ \Lambda_0(m) = \frac{1}{2} \left[ \bar{r}_1 + \bar{r}_2 + \frac{1}{2\pi} \int_0^{2\pi}\sqrt{(r_1(\theta) - r_2(\theta)^2+ 4 m^2}d\theta \right] - m$$
 
 In the present paper we give closed expression of $\omega_c$ for specific $r_i$ that we precise now.

\subsection{The two-patches model in the $(U,V)$ variables.}

We consider the model :
\beq \label{BLSS1}
\Sigma(r_1,r_2,m,T)\quad \quad\left\{
\begin{array}{lcl} 
\displaystyle \frac{dx_1}{dt} &=& r_1(t) x_1 +m( x_2-x_1) \\[6pt]
\displaystyle \frac{dx_2}{dt}& =&r_2(t) x_2 +m(x_1-x_2)
\end{array} 
\right.
\feq
where $r_1(t)$ and $r_2(t)$ are the growth rates at time $t \geq 0$ in patch 1 and 2, respectively. The functions $t \mapsto r_1(t) $ and $t \mapsto r_2(t) $ can be deterministic or random. In this paper, we will be interested in the case where these functions are piecewise constant, and change of values at periodic (see Section \ref{basic}) or random  times (see Section \ref{sec:sto}). Nevertheless, for the rest of this subsection, the precise form of $r_1$ and $r_2$ do not matter.

Thanks to the fact that we have only two patches we can define
\beq \label{varUV}
\begin{array}{l}
\displaystyle U = \ln(\sqrt{x_1x_2}) = \frac{1}{2}(\ln(x_1) +\ln(x_2) )\\[8pt]
\displaystyle V = \ln(\sqrt{x_1/x_2}) = \frac{1}{2}(\ln(x_1) -\ln(x_2) )
 \end{array}
\feq
This is legitimate since  the solutions of  \eqref{BLSS1} remain strictly positive as soon as the initial conditions are.
In these new variables the system becomes:
\footnote{We recall the notations $\sh(x) = \frac{\emat^x-\emat^{-x}}{2} \quad \quad \ch(x) = \frac{\emat^x+\emat^{-x}}{2} $}  :
 \beq \label{systvarUV}
  S(r_1,r_2,m,T)\quad \quad\quad \quad
\left\{
\begin{array}{lcl}
\displaystyle  \frac{dU}{dt} &=&\displaystyle \frac{r_1(t)+r_2(t)}{2} +m \,( \ch(2V) -1)  \\[8pt]
\displaystyle  \frac{dV}{dt} &=&\displaystyle  \frac{r_1(t)- r_2(t)}{2}-m\,\sh(2V)
 \end{array} 
 \right. 
\feq
One observes that the variable $V$ is decoupled from the variable $U$.
 By the way, once the solution $V(t)$ is known, the solution $U(t)$ is obtained by the simple quadrature
  \beq 
  \label{eq:Ugeneral}
\begin{array}{lcl}
\displaystyle  U(t)  &=&\displaystyle  U_0 + \int_0^{t}  \frac{r_1(\theta )+r_2(\theta )}{2}d\theta  +m \,\int_0^{t}( \ch(2V(\theta) )-1) d\theta  
 \end{array} 
\feq
 We can already make some remarks. When both sites are sinks ($\bar r_1 <0$ and $\bar r_2<0$), if there is no migration, as expected, the metapopulation is decreasing. As the quantity $( \ch(2V)-1)$ is strictly positive as soon as $V$ is different from $0$ we see that the more $V(\theta)$ will be (on average) different from $0$, the bigger the second integral term will be and so the more $U(t)$ will have the possibility to become positive. A quick look at the equation of $V$ shows that the larger $|r_1(t)-r_2(t)|$ is, the larger it will be  ; and the bigger is $m$ the smaller is $V$. This immediately tells us two ingredients favorable to inflation: $r_1(t)$ and $r_2(t)$ must be different, $m$ must be strictly positive for the second integral to be taken into account, but not too large so that the solutions of the second equation of \eqref{BLSS1} are not too small. 
 
Note that there are many ways to transform linear systems in a cascade of non linear systems like, for instance, polar coordinates, but this one seems the most appropriate to the study of migration between two patches (see Remark \ref{rem:avantageU} below). Moreover $U$ is the logarithm of the geometric mean of the abundances on the two patches and the absolute value of $V$ is the logarithm of their geometric standard deviation which have biological meaning. 

\begin{remarque}
\label{rem:commongrowth}
{\rm Assume that $r_1$ and $r_2$ are bounded by $R > 0$. Then, from the second equation  in \eqref{systvarUV}, it is easily seen that, as soon as $m > 0$,  $V$ will eventually enter and remain in the compact interval $[ -\sh^{-1}(R/m), \sh^{-1}(R/m)]$. As a consequence, $\lim_{t \to \infty}  V(t) /t = 0$, and thus (provided the limits exist)
\[
\Lambda[x_1] := \lim_{t \to \infty} \frac{\ln(x_1(t))}{t} =  \lim_{t \to \infty} \frac{\ln(x_2(t))}{t} =: \Lambda[x_2]
\] 
In other words, the long term growth rate is common in the two patches. Moreover, 
$$\lim_{t \to \infty}\frac{U(t)}{t} = \frac{1}{2} \left(\lim_{t \to \infty}\frac{\ln(x_1(t))}{t}+\lim_{t \to \infty} \frac{\ln(x_2(t))}{t}\right) = \Lambda[x_i]$$
and so we are interested in the growth or decay of $U$.
}
\end{remarque}
$\,$\\
\begin{remarque}
\label{rem:avantageU}
{\rm Assume that the following limits exist
$$\displaystyle \overline{r}_i:= \lim_{t \to \infty} \frac{1}{t}\int_0^t r_i(s) ds\quad \quad i = 1, 2 $$
Note that these limits indeed exist if $r_i$ are periodic or semi-Markov (see Remark \ref{rem:semiMarkov} in Section \ref{sec:random_switch}). Then, $\overline{r}_i$ is the long term average growth rate on each patch: if $m=0$, $x_i(t)$ tend to $0$ or infinity depending on whether $\overline{r}_i$ is negative or positive. From Equation \eqref{eq:Ugeneral} and the fact that $\ch(x) \geq 1$ for all $x \in \mathbb{R}$, we deduce that, for all $m \geq 0$,
\[
\liminf \frac{U(t)}{t} \geq \frac{\bar r_1 + \bar r_2}{2}.
\]
This means that the common growth rate of the two patches is always higher than the mean of the growth rates within each patch. This is straightforward in the variables $U - V$, while it seems difficult to conclude that from the classical "polar" decomposition, $S(t) = x_1(t) + x_2(t)$ et $y(t) = x_1(t) / S(t)$, which leads to 
\[
\frac{d S(t)}{dt} = S(t) \left( r_1(t) y(t) + r_2(t) ( 1 - y(t)) \right),
\]
and thus
\[
\frac{1}{t} \ln( x_1(t) + x_2(t)) = \frac{S(0)}{t} + \frac{1}{t} \int_0^t r_1(s) y(s) + r_2(s) ( 1 - y(s))ds.
\]
In Section \ref{nonsym}, we prove by adapting conveniently the variables $U - V$, that a similar result holds true in the case of a non symmetric dispersal, with a weighted mean of $\overline
r_1$ and $\overline{r}_2$ taking into account the asymmetry in the dispersal.
}
\end{remarque}

\subsection{The $(±1)$model in periodic environment}
\label{basic}
\subsubsection{The model}
\label{modelpm}
Our idea is to understand the mathematics of the simplest possible model of the form $\Sigma(r_1,r_2,m,T)$ and complicate it thereafter. Thus, we consider the system 
 \beq \label{systeme1}
\Sigma(\varepsilon,m,T) \quad
\left\{
\begin{array}{lcr}
\displaystyle  \frac{dx_1}{dt} &=&(+u(t)-\varepsilon)x_1+m(x_2-x_1)\\[8pt]
\displaystyle  \frac{dx_2}{dt} &=& (-u(t)-\varepsilon)x_2+m(x_1-x_2) 
 \end{array} 
 \right.
\feq
where $0 \leq \varepsilon \leq 1$, $0 \leq m$, $0 \leq T$ and the function   $t \mapsto u(t)$ is periodic of period $2T$, with :
$$ t \in [0,\,T[ \Rightarrow u(t) = 1,\quad \quad  t \in [T,\,2T[ \Rightarrow u(t) = -1$$
We call the system $\Sigma(\eps,m,T)$ the \textbf{ periodic $(±1)$-model}. From that definition one sees that for $u(t) = +1$  we are integrating the system 
 \beq 
 \label{Sigma+}
 \begin{array}{lcl}
	\Sigma^+(\varepsilon,m)&\quad\quad \quad&
	\left\{
	\begin{array}{lcl}
		\displaystyle  \frac{dx_1}{dt}& = &(+1-\varepsilon)x_1+m(x_2-x_1) \\[8pt]
		\displaystyle  \frac{dx_2}{dt}&=&  (-1-\varepsilon)x_2+m(x_1-x_2) 
 	\end{array} 
	 \right.
 \end{array}
 \feq
while for $u(t) = -1$ we are integrating the system 
 \beq 
 \label{Sigma-}
 \begin{array}{lcl}
\Sigma^-(\varepsilon,m)&\quad\quad \quad&
\left\{
	\begin{array}{lcl}
	\displaystyle  \frac{dx_1}{dt}& = &(-1-\varepsilon)x_1+m(x_2-x_1) \\[8pt]
	\displaystyle  \frac{dx_2}{dt}&=&  (+1-\varepsilon)x_2+m(x_1-x_2) 
	 \end{array} 
 \right.
 \end{array}
 \feq
Thus we are switching, each $T$ units of time, from system $\Sigma^+(\varepsilon,m)$ to system $\Sigma^-(\varepsilon,m)$ and vice versa; such systems are called  {\em switched} systems. Switched systems where intensively investigated in the context of control theory during the seventies and later (see for instance \cite{JUR97}) and more recently, in a probabilistic context, under the name of PDMP (Piecewise Deterministic Markov Processes) \cite{BH12, BMZIHP,BS19,DAV84,HS19}.
\begin{remarque}{ \rm Let us remark that the $(\pm 1)$-model defined by \eqref{systeme1} is a bit more general than it looks since it includes the case of two identical patches that are simply in ''phase opposition'', given by the following switched system :
\begin{equation}\label{+r-d}
\begin{array}{ll}
\mbox{For }t\in[0,T),&
\left\{
	\begin{array}{lcr}
 \frac{dx_1}{dt}& = &rx_1+m(x_2-x_1) \\[4pt]
 \frac{dx_2}{dt}&=&  -dx_2+m(x_1-x_2) 
 	\end{array} 
	 \right.
\\[4pt]
\mbox{For }t\in[T,2T),&
\left\{
	\begin{array}{lcr}
 \frac{dx_1}{dt}& = &-dx_1+m(x_2-x_1) \\[4pt]
 \frac{dx_2}{dt}&=&  rx_2+m(x_1-x_2) 
	 \end{array} 
 \right. 
\end{array}
\end{equation}
This system is of the form 
$\Sigma(r_1,r_2,m,T)$
given in \eqref{BLSS1}, where $r_1(t)$ and $r_2(t)$ are the $2T$-periodic functions defined by
$$
r_1(t)=\left\{
\begin{array}{rll}
r&\mbox{if}&t\in[0,T)\\
-d&\mbox{if}&t\in[T,2T)
\end{array}
\right.
\quad
r_2(t)=r_1(t+T)=\left\{
\begin{array}{rll}
-d&\mbox{if}&t\in[0,T)\\
r&\mbox{if}&t\in[T,2T)
\end{array}
\right.
$$
We assume that $d>r>0$. We have
$$
\bar{r}_1=\bar{r}_2=\frac{1}{2T}\int_0^{2T}r_i(t)dt=\frac{r-d}{2}<0,\quad
\chi=\frac{1}{2T}\int_0^{2T}\max(r_1(t),r_2(t))dt=r>0,
$$
which means that each patch is a sink, while $\chi>0$. According to the theorem of Katriel (see Theorem \ref{theoremKatriel})
inflation can occur. 
Let $\varepsilon=\frac{d-r}{d+r}$. We have $1-\varepsilon=\theta r$ and 
$1+\varepsilon=\theta d$, where $\theta=\frac{2}{d+r}$. Therefore, the change of time $t=\theta s$ transforms the switched system \eqref{+r-d} into the system
\begin{equation*}
\begin{array}{ll}
\mbox{For }s\in\left[0,\frac{T}{\theta}\right),&
\left\{
	\begin{array}{lcr}
 \frac{dx_1}{ds}& = &(1-\varepsilon)x_1+\theta m(x_2-x_1) \\[4pt]
 \frac{dx_2}{ds}&=&  -(1+\varepsilon)x_2+\theta m(x_1-x_2) 
 	\end{array} 
	 \right.
\\[4pt]
\mbox{For }s\in\left[\frac{T}{\theta},2\frac{T}{\theta}\right),&
\left\{
	\begin{array}{lcr}
 \frac{dx_1}{ds}& = &-(1+\varepsilon)x_1+\theta m(x_2-x_1) \\[4pt]
 \frac{dx_2}{ds}&=&  (1-\varepsilon)x_2+\theta m(x_1-x_2) 
	 \end{array} 
 \right. 
\end{array}
\end{equation*}
This is the periodic $(\pm 1)$-model given by \eqref{systeme1}
, with $T$ replaced by $T/\theta$ and $m$ replaced by $\theta m$, that is
$\Sigma(\varepsilon,\theta m,T/\theta)$
}
\end{remarque}

\subsubsection{The $(±1)$model in the variables  $(U,V)$}
In the new variables $(U,V)$ the system becomes
 \beq \label{plusoumoins1enUV}
  S(\varepsilon,m,T)\quad \quad\quad \quad
\left\{
\begin{array}{lcl}
\displaystyle  \frac{dU}{dt} &=&\displaystyle  m\,\ch(2V)-m-\varepsilon  \\[8pt]
\displaystyle  \frac{dV}{dt} &=&u(t)-m\,\sh(2V)
 \end{array} 
 \right. 
\feq
The non autonomous system
 \beq \label{SV}
  F(m,T)\quad \quad\quad \quad
\left\{
\begin{array}{lcl}
\displaystyle  \frac{dV}{dt} &=&u(t)-m\,\sh(2V)
 \end{array} 
 \right. 
\feq
  is a one dimensional switched system between the two autonomous equations
\beq \label{S1}
 F_m^+ \quad \quad\quad\quad \quad\quad \quad \left\{\frac{dV}{dt} = +1 -m\,\sh(2V) \right. \quad 
\feq
and
\beq \label{S2}
F_m^- \quad \quad\quad \quad \quad\quad \quad  \left\{ \frac{dV}{dt} =-1-m\,\sh(2V)\right.\quad 
\feq 
\begin{figure}
 \begin{center}
 \includegraphics[width=1\textwidth]{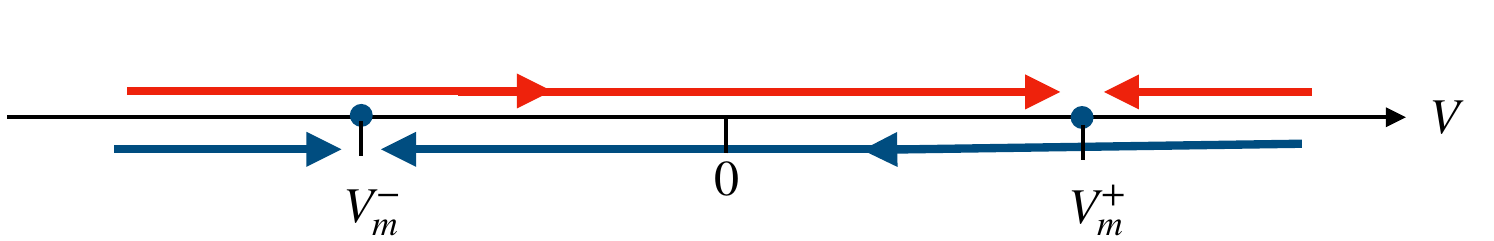}
 \caption{The switched system  $F(m,T)$: in red $F^+_m$,  in blue $F^-_m$, defined by
\eqref{S1} and \eqref{S2} respectively}\label{F(m,T)}
 \end{center}
 \end{figure}
 We let $\varphi^+_t(v)$ and $\varphi^-_t(v)$ the solutions to \eqref{S1} and \eqref{S2} at time $t \geq 0$, starting from $v$ at time $0$. The two differential equations $ F_m^+ $ and $ F_m^- $  have respectively the points 
\begin{equation}
\label{Vm+Vm-}
V_m^+ =  \frac{1}{2}\sh^{-1}(+1/m)  \quad \quad V_m^- =  \frac{1}{2}\sh^{-1}(-1/m)
\end{equation}
 as globally asymptotically   stable equilibria (that is, $\varphi^+_t(v)$ converges to $V_m^+$ and $\varphi^-_t(v)$ converges to $V_m^-$). From Figure \ref{F(m,T)} it is evident that the solutions  of $F(m,T)$ are trapped in the interval  
$[V_m^-\,,\,V_m^+]$.
The following proposition is easy to prove with elementary calculus means.  Since we don't know reference for it a proof is done in Appendix \ref{apendiceperiodic} :
{\proposition \label{prop1} The switched system $F(m,T)$ has a unique periodic solution, denoted by $P_{m,T}(t)$, globally asymptotically stable, which oscillates between two values $P^-_{m,T}$, and $P^+_{m,T}$ contained in the interval $[V_m^-,\; V_m^+]$ ;  $P^-_{m,T}= - P^+_{m,T}$ and 
the function $T \mapsto P^+_{m,T}$ is an increasing function of $T$ which tends to $V_m^+$ when $T$ tends to infinity.}

Let  us denote :
\beq 
\begin{array}{lcl} \label{deltaT2}
\displaystyle  \Delta(\eps,m,T)   = \displaystyle \frac{1}{2T}  \int_0^{2T} m\,\ch(2P_{m,T}(s))-m-\varepsilon ds . 
 \end{array} 
\feq
Then,
{\proposition \label{prop2} 
\beq
\label{eq:UtoDelta}
\lim_{t \to \infty} \frac{U(t)}{t} = \Delta( \varepsilon, m, T).
\feq
Hence,  $U(t)$ tends to $± \infty$ according to the sign of $\Delta(\eps,m,T)$.}\\
\textbf{Proof} From the first equation of \eqref{plusoumoins1enUV} one has 
\beq
U(t) = U(0) + \int_0^t \varphi(V(s))ds 
\feq
with 
\beq
\varphi(V) = m \cosh(2V)-m - \varepsilon
\feq
We claim that the following limits exist and are equal
 $$\lim_{t \to \infty} \frac{U(t)}{t} = \lim_{n \to \infty} \frac{U(n2T)}{n 2 T}.$$
Therefore, for all $n_1 \geq 0$,
\[
\lim_{t \to \infty} \frac{U(t)}{t} =\lim_{n \to \infty} \frac{1}{n 2T}\left(U(0)+  \int_0^{n_12T}\varphi(V(s))ds +  \int_{n_12T}^{n2T}\varphi(V(s))ds\right). 
\]
Since the solution $V(s)$ converges to the periodic solution $P_{m,T}(t)$ for $n_1$ large enough we can replace in the second integral $V(s)$ by  $P_{m,T}(s)$ and then
\[
\lim_{t \to \infty} \frac{U(t)}{t} =\lim_{n \to \infty} \frac{1}{n 2T}\left(U(0)+  \int_0^{n_12T}\varphi(V(s))ds\right) + \frac{1}{n 2T}\left( \int_{n_12T}^{n2T}\varphi(P_{m,T}(s))ds\right).
\]
The first term tends to $0$ as $n \to \infty$ and the second reads
\[
\frac{1}{n 2T}\left( \int_{n_12T}^{n2T}\varphi(P_{m,T}(s))ds\right)  = \frac{(n-n_1)2T}{n_12T+ (n-n_1)2T} \frac{1}{2T} \int_0^{2T} \varphi(P_{m,T}(s))ds,
\]
which limit for $n \to \infty$ is just
$$ \frac{1}{2T} \int_0^{2T} \varphi(P_{m,T}(s))ds  = \frac{1}{2T}  \int_0^{2T} m \cosh(2P_{m,T} (s) )-m - \varepsilon ds.$$
It remains to prove the claim. For $t \in [2nT, 2(n+1)T)$, we have
\[
\frac{U(t)}{t} = \frac{2nT}{t} \frac{U(2nT)}{2nT} + \frac{1}{t} \int_{2nT}^t \varphi(V(s)) ds
\]
When $t$ goes to infinity, $\frac{2nT}{t}$ goes to $1$, whereas, since $V(s)$ is bounded, the second term in the right hand side goes to $0$. This entails the claim and concludes the proof.
$\Box$

\begin{figure}
  \begin{center}
 \includegraphics[width=0.9\textwidth]{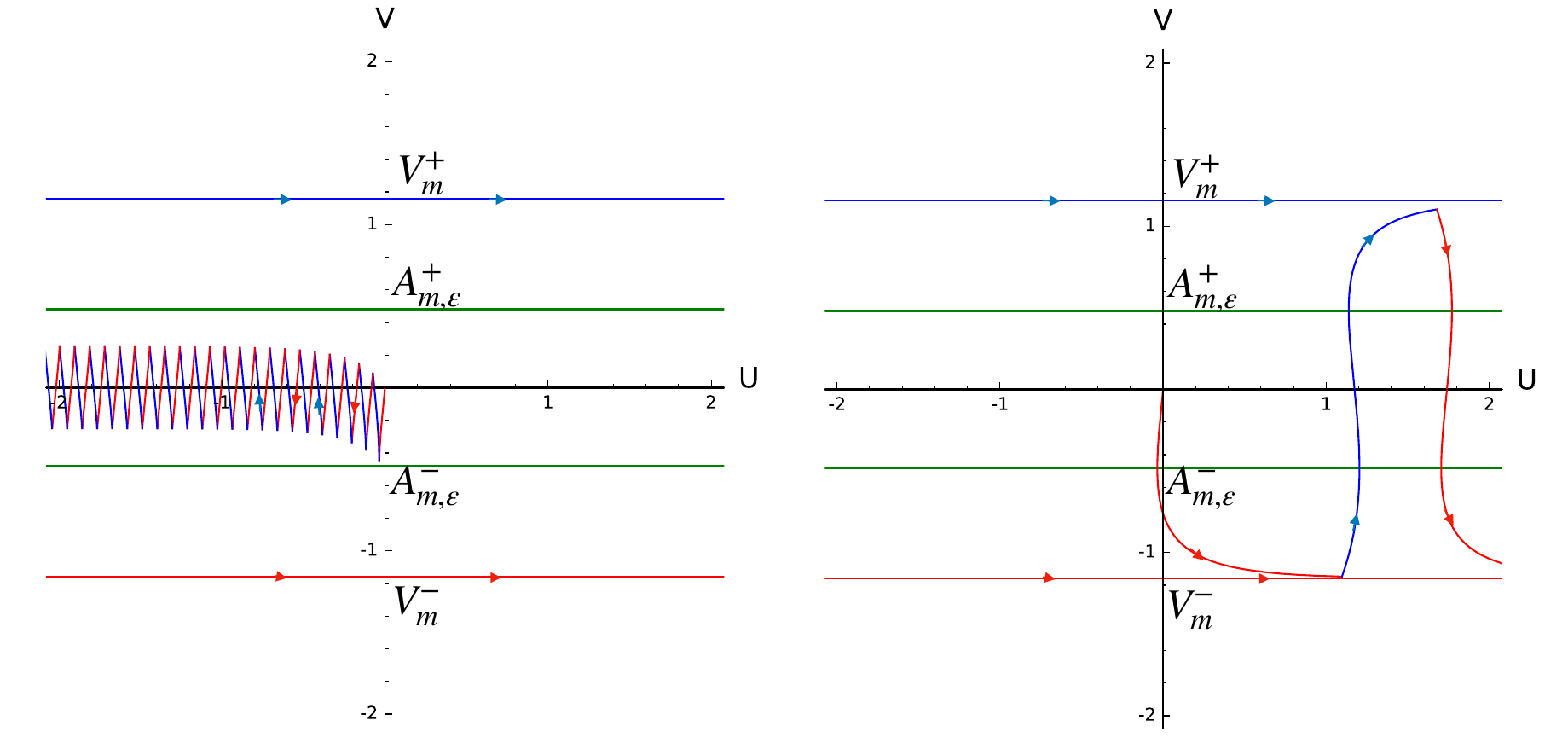}
 \caption{The switched system  $S(\varepsilon, m,T)$: in blue $u(t) = +1$, in red $u(t) = -1$. $\varepsilon = 0.1\;m = 0.2$,  $T = 0.5$ (left), $T =3$, (right).} \label{periodic} 
 \end{center}
 \end{figure}

Let us evaluate $\Delta(\eps,m,T)$.
Since the function $V \mapsto \left(m\,\ch(V)-m-\varepsilon \right) $ is even and, as  it is easily seen,  $S^+(\varepsilon,m)$ and $S^-(\varepsilon,m)$ are symmetric  with respect to the horizontal axe, we have
\beq \label{formuledelta}
\Delta(\eps,m,T)  =  \frac{1}{T}\int_0^{T} \ch(2P_{m,T}(s))-m-\varepsilon ds.
\feq
Recall (see \eqref{plusoumoins1enUV}) that $U$ is solution to 
\begin{equation}
\label{eq:Uintegral}
U(t) = U_0 + \int_0^t   \left(m\,\ch(V(s))-m-\varepsilon \right) ds.
\end{equation}
Elementary computations shows that
\begin{proposition}$\,$
 The right member of Equation \eqref{eq:Uintegral} is strictly negative on the interval $ ]A_{\eps,m}^ - , A_{\eps,m}^+[$ and positive outside with 
$$ A_{m,\eps}^- = - \frac{1}{2} \cosh^{-1} \left(1+\frac{\eps}{m}\right) \quad \quad A_{m,\eps}^+ = \frac{1}{2} \cosh^{-1} \left(1+\frac{\eps}{m}\right).  $$
Moreover, 
\ben
\item     $m  >  \frac{1 - \eps}{2 \eps} \Longrightarrow  [V^-_m,V^+_m] \subset [ A^-_{\epsilon,m}, A^+_{\epsilon,m}]$
\item $m < \frac{1 - \eps}{2 \eps} \Longrightarrow [V^-_m,V^+_m] \supset [A^-_{\epsilon,m}, A^+_{\epsilon,m}]$
\fen
\end{proposition}
From this proposition one can see what is going on. On figure \ref{periodic} one sees a simulation in the plane $(U,V)$ of trajectories in the case n° 2. On the left, the period is short which results in the amplitude of the oscillations of the periodic solution being small. In this case the trajectory of the periodic solution remains largely inside the stripe $V \in [V^-_m,V^+_m]$ where $\frac{dV}{dt} < 0$ which results in the decrease of $U$. On the contrary, when the period is large ($T = 3$ in the simulation on the right) $v(t)$ has time to approach and stay close to $V^-_m$ or $V^+_m$ where $\frac{dV}{dt} > 0$ which leads to the growth of $U$.

In other words, since $V = \ln(\sqrt{x_1/x_2})$ measures the ''asymmetry'' between the abundances on the two patches, we can see that when the period is important, the most loaded site is little diminished by the migration towards the less loaded patch, and, on the contrary, this last one sees its population strongly increase, which increases the product $x_1x_2$ and then $U$.

A more precise description of the behavior  of $\Delta(\eps,m,T)$ is given by the following proposition which details of the proof are given in Appendix \ref{appendicedelta}.

\begin{proposition}
 \label{propgenerale} 
 \textbf{ Properties of $\Delta(\eps,m,T)$}
\ben
\item For fixed $T >0$, for both small and large values of  $m$, $\Delta(\varepsilon,m,T) <0$ and thus if there is inflation it must be for some intermediate value of the migration $m$.

\item For fixed $\eps> 0$ and $m <\frac{1-\eps^2}{2\eps}$,  there exists a threshold $T^*(\eps,m)$ such that for $T < T^*(\eps,m)$,  $\Delta(\varepsilon,m,T) <0$ and there is no inflation  while for  $T > T^*(\eps,m)$ , $\Delta(\varepsilon,m,T) >0$ and there is inflation.

\item For every $\eps > 0$, the minimum of $T^*(\eps,m)$ over $m$ is strictly positive. In other words there exists a threshold $T^{**}> 0$ such that for $T < T^{**}$, for all values of $m$, $\Delta(\eps, m, T) < 0$ and there is no inflation.
\fen\end{proposition}

\subsubsection{An explicit formula for $\Delta(\eps,m,T)$}

\begin{figure}
  \begin{center}
 \includegraphics[width=0.8\textwidth]{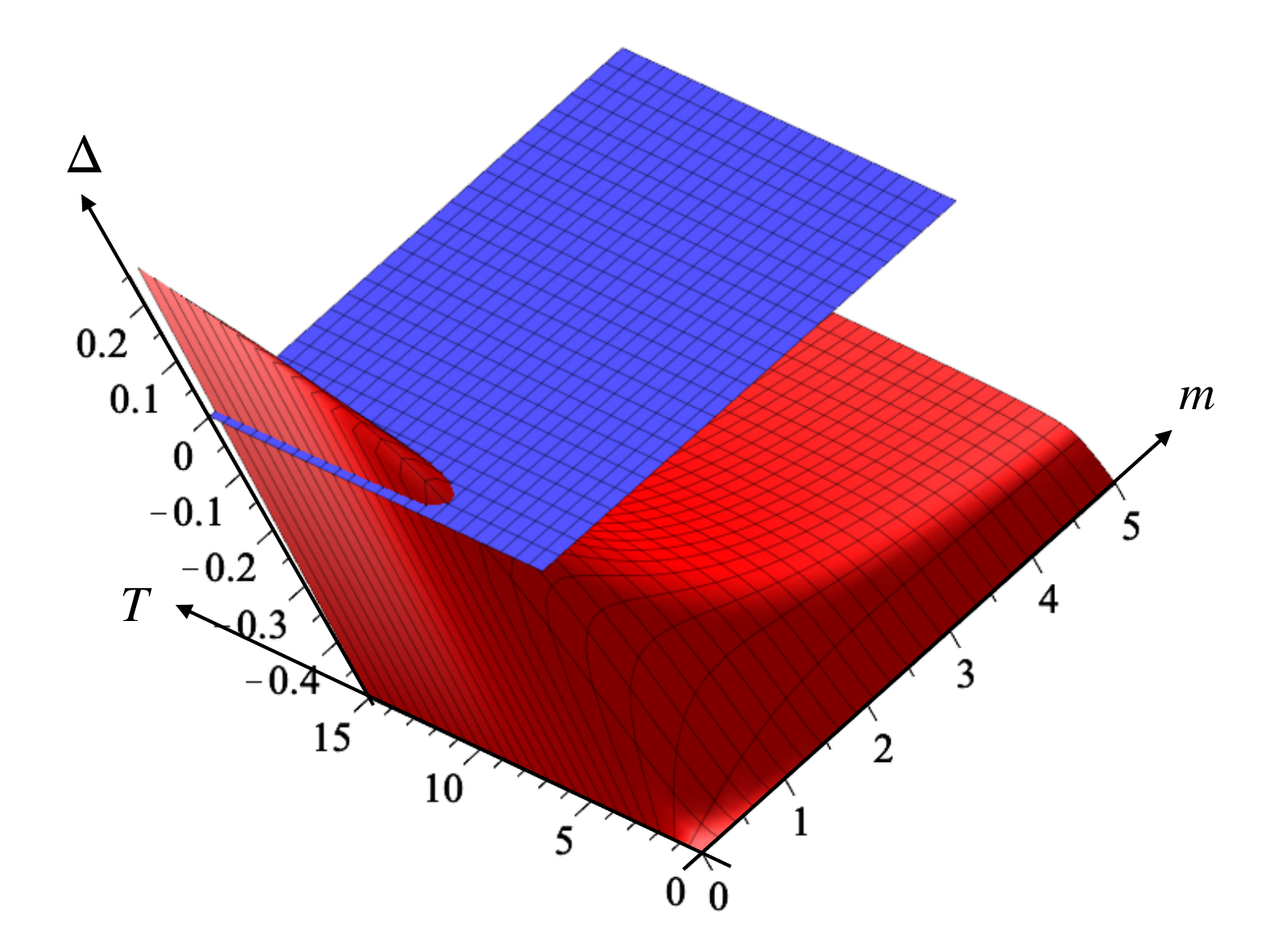}
 \caption{Graph of $\Delta(0.5,m,T)$} \label{graphdelta}
 \end{center}
 \end{figure}

By an elementary but not immediate computation (see appendix \ref{appendicedeltaexplicite}) one proves the following explicit formula for $\Delta(\eps,m,T)$ 
\begin{proposition}\label{explicitdelta}
Let us denote
$$  b=e^{T\sqrt{1+m^2}} \quad \quad C=m^2b^4+2m^2b^2+4b^2+m^2.$$
Then one has
\begin{equation}\label{Delta1}
\Delta(\eps,m,T)=\frac{1}{2T}\ln\frac{m^2b^4+2b^2+m^2+
m(b^2-1)\sqrt{C}}{2(1+m^2)b^2}
-(m+\varepsilon)
\end{equation}
\end{proposition}
On Figure \ref{graphdelta} this formula is used to draw the picture (using the software Maple) of the graph of $\Delta(0.5,m,T)$ with respect to the variables $(m,T)$.

\paragraph{Small and large values of $m$ or $T$.}
On Figure. \ref{graphdelta} one sees that 
$\Delta(0.5,0,T) = - 0.5 =  -  \varepsilon$,  
which is easily understandable: for $m = 0$ one has $\frac{dx_1}{dt} = (u(t)-\eps)x_1,\, \frac{dx_2}{dt} = (- u(t)-\eps)x_1)$ taking the mean of the two logarithm $\frac{U}{dt} = -  \eps$. For $T = 0$, one sees that $\Delta(0.5,m,0)$ is equal to $-1 = - \eps$ which is explained by the general fact (see e.g. \cite{FW,JUR97}) that if we consider a  switched system at a rate witch tends to infinity (i.e. $T\rightarrow 0$ ) then the solutions tend to solutions of the system which is the mean of the two systems ; in our case the mean of the two systems is :
\beq
\begin{array}{lcl}
\displaystyle \frac{dU}{dt}& =& \displaystyle m\cosh(2V)-m-\eps \\[8pt] 
\displaystyle \frac{dV}{dt}& =&\displaystyle  \frac{(u(t)-m\sinh(2V)) + (- u(t) -m\sinh(2V)) }{2}= -m\sinh(2V)   
\end{array}
\feq 
which, after a transient, are just $\displaystyle \frac{dU}{dt} = -  \eps$

The asymptotic behaviors for small and large values of $m$  or $T$ can also be derived by basic development on the explicit formula \eqref{Delta1}, as shown by the next proposition.
\begin{proposition}
\label{prop:asymptoticDelta}
For fixed value of $m$, one has
$$
\lim_{T \to \infty} \Delta(\varepsilon, m, T) =   \left( \sqrt{1+m^2} - (m + \varepsilon) \right)
\quad \quad 
\lim_{T \to 0} \Delta(\varepsilon, m, T) = -  \varepsilon.
$$
For fixed value of $T > 0$, 
$$
\lim_{m \to 0}\Delta(\varepsilon, m, T) = \lim_{ m \to +\infty}  \Delta(\varepsilon, m, T) = -  \eps.
$$

\end{proposition}
This proposition tells us that, if $m < \frac{1 - \varepsilon^2}{2 \varepsilon}$, then for $T$ sufficiently large, inflation occurs. In addition, if $m$ is fixed and $T$ is small or if $T$ is fixed and $m$ too small or too large, there is no inflation.

\paragraph{Threshold value of $m$ for large $T$.}
\begin{figure}
  \begin{center}
 \includegraphics[width=1\textwidth]{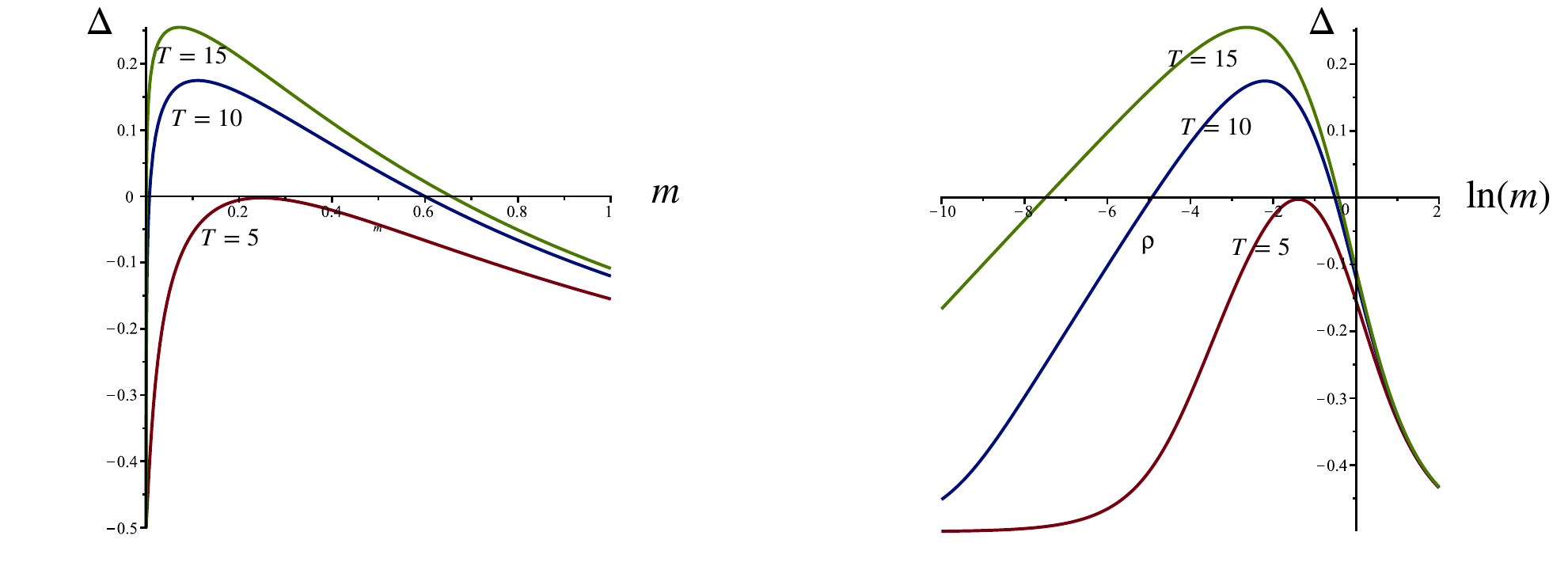}
 \caption{Graphs of $m \mapsto \Delta(\eps,m,T)$ : $T = 5,10,15$ ;  $\eps = 0.5$ }
  \label{seuil}
 \end{center}
 \end{figure}
  On figure \ref{graphdelta} one sees that, for large values of $T$ the dependence with respect to $m$ is very sharp close to $0$ ; in order to have a better understanding of what is going on around $0$ we ask to Maple to draw the graph of $m \mapsto \Delta(\eps,m,T)$ for three values of $T$ : $5, 10, 15$. The result is shown on figure \ref{seuil}. On the left one sees that the threshold for the appearance  of positive values of $\Delta$ is very  small  and on the right we see the same graphs but with a logarithmic scale for $m$.  From the picture we guess the following property,  which is confirmed by the mathematical derivation of appendix \ref{asymptotics}:

  \begin{proposition} \label{approxseuil}
  When $T$ is large ($T\rightarrow +\infty$) the threshold value at which $m \mapsto \Delta(\eps,m,T)$ becomes positive is the exponentially small value :
  \begin{equation}
  \label{eq:approxseuil}
  m^*(\eps,T) \sim \emat^{-(1-\eps)T}
\end{equation} 
  \end{proposition} 
  Recall that if $m = 0$, there is no inflation. The above proposition states that the threshold value of $m$ at which inflation occurs is exponentially small for $T$ large enough. As is appears from Figure \ref{graphdelta}, the approximation \eqref{eq:approxseuil} works extremely well for $T$ larger than $5$. This result is quite striking. For example, for $T = 10$, the system goes from no inflation for $m=0$ to inflation for $m = 10^{-5}$ !
 \begin{remarque}{\rm 
 This proposition  gives an affirmative answer to \textbf{Conjecture 3} of Katriel's paper \cite{KAT22} in the particular case of  piecewise constant model which takes advantage of explicit formulas for the solutions. 
 In a forthcoming paper we shall prove this conjecture in the general case \cite{LOB22}.
 }
 \end{remarque}
\subsubsection{Back to the variables $x_1,x_2$}\label{variables initiales}

One has $\exp(U(t)) =\sqrt{ x_1(t)x_2(t)}$ and we know that for large values of  $t$  one has  $-V^+_m<V(t) < +V^+_m $ which means $\frac{1}r < \frac{x_2}{x_1} < r$, with $r = \exp(V^+_m)$. From this we deduce that :
\beq \label{Sdelta}
\Sigma(\varepsilon,m,T)\;\mathrm{is \; stable} \Longleftrightarrow \Delta(\varepsilon,m,T)<0
\feq
On the other hand let us consider  the  ''period mapping'' of $\Sigma(\varepsilon,m,T)$, that is to say the linear mapping which, to an initial condition $(x_{1}(0),x_{2}(0))$ at time $0$,  assigns the solution of $\Sigma(\varepsilon,m,T)$ at time $2T$ and let us denote it by 
\beq 
\left(
	\begin{array}{l}
	x_1(2T) \\
	x_2(2T) 
	\end{array}
	\right)
= M(\varepsilon,m,T)
\left(
	\begin{array}{l}
	x_1(0) \\
	x_2(0) 
	\end{array}
	\right)
\feq 
The stability of our system $\Sigma(\varepsilon,m,T)$ is equivalent to the stability of the linear discrete system of $\mathbb{R}^2$ 
\beq
X_{n+1} = M(\varepsilon,m,T)X_n
\feq 
For $u = - 1, 1$, let us denote by $M _{\varepsilon,m}^u$ the matrix 
\beq
\label{Muepsm}
M _{\varepsilon,m}^u = 
\left(
\begin{array}{cc}
u-m-\varepsilon &+m\\[14pt]
+m & -u-m -\varepsilon
\end{array}
\right)
\feq
With this  notation the matrix  $M(\varepsilon,m,T)$ is given by 
\beq
\label{MepsmT}
M(\varepsilon,m,T) = \emat^{TM_{\varepsilon,m}^{-1}} \emat^{TM_{\varepsilon,m}^{+1}} 
\feq
The stability of the discrete system $X_{n+1} = M(\varepsilon,m,T)X_n$    is decided by the spectral radius 
\beq \label{rayonspectral} 
\sigma(\varepsilon,m, T) = \max |\lambda_i(\varepsilon,m, T)|\quad i=1,2
\feq
where $\lambda_i(\varepsilon,m, T)$ are the two real eigenvalues of $M(\varepsilon,m,T)$ (note that since $M_{\varepsilon,m}^{+1} $ and $M_{\varepsilon,m}^{-1}$ are symmetric,  so is $ M(\varepsilon,m,T)$ and its eigenvalues are real).  Thus 
\beq \label{Ssigma}
\Sigma(\varepsilon,m,T)\;\mathrm{is \; stable} \Longleftrightarrow \sigma(\varepsilon,m,T)<1
\feq
In view of (\ref{Sdelta}) and (\ref{Ssigma}) there must be a connection between $\Delta(\varepsilon,m,T)$ and $\sigma(\varepsilon,m,T)$. The connection is given by the following proposition, which is proved in appendix \ref{sigmadelta}.
\begin{proposition}\label{formule}
\beq
\Delta(\varepsilon,m,T) = \frac{1}{2T}\ln( \sigma(\varepsilon,m,T))
\feq
\end{proposition} 

\subsection{The $(\pm 1)$ model in stochastic environment}
\label{sec:sto}
\subsubsection{Random choice of switching times}
\label{sec:random_switch}
In the previous section, the switching from system $\Sigma^{+}( \varepsilon, m)$ to $\Sigma^{-}( \varepsilon, m)$
and vice versa, occurs after a fixed deterministic time $T$. Considering that switching from one system to the other models is a change in the environment, it makes sense to deal with the case where switching occur after a random time. More precisely, we consider a sequence of iid\footnote{'' iid '' means Independent, Identically Distributed.} random variables $(S_n)_{n \geq 0}$, with common law $\mu$ on $[0, \infty)$, and a random function $t \mapsto \bm{u}(t)$ such that,
\[
t \in [ T_{2n}, T_{2n+1} [ \Rightarrow \bm{u}(t) = 1\quad \quad  t \in [T_{2n+1}, T_{2n+2}[ \Rightarrow \bm{u}(t) = -1
\]
where $T_0 = 0$ and for $n \geq 1$, $T_n = \sum_{k=1}^n S_k$. From this function $\bm{u}$, we build a stochastic process solution to
\beq \label{systeme1random}
\bm{\Sigma}(\varepsilon,m,\mu) \quad
\left\{
\begin{array}{lcr}
\displaystyle  \frac{dx_1}{dt} &=&(+\bm{u}(t)-\varepsilon)x_1+m(x_2-x_1)\\[8pt]
\displaystyle  \frac{dx_2}{dt} &=& (-\bm{u}(t)-\varepsilon)x_2+m(x_1-x_2) 
 \end{array} 
 \right.
\feq
In other words, after the $n$-th switching, we draw a random variable $S_{n+1}$ with law $\mu$ independent from anything else, 
and we integrate system $\Sigma^{+}( \varepsilon, m)$ when $n$ is even ($\Sigma^{-}( \varepsilon, m)$ when $n$ is odd) for a time $S_n$. The periodic system studied in the previous section is the particular case when the law  $\mu$ is the Dirac mass at $T$, i.e., $S_n = T$ almost surely for all $n \geq 1$. 

We let $\mathbb{E}(Y)$ denote the expectation of a random variable $Y$ and $\mathbb{P}(A)$ the probability of an event $A$. Since $\mathbb{E}(S_1)$ represents the mean time spent in each regime, we assume that $ \mathbb{E}(S_1) < + \infty$. In addition, we assume that $\mathbb{P}(S=0) = 0$, in order to avoid instantaneous change of regime.

As in the case of periodic environment, we perform the change of variable $V=\frac{1}{2}(\ln(x_1) - \ln(x_2))$ and $U = \frac{1}{2}(\ln(x_1) + \ln(x_2))$ to get 
\beq 
\label{Srandom} 
  \mb{S}(\varepsilon,m,\mu)\quad \quad\quad \quad
\left\{
\begin{array}{lcl}
\displaystyle  \frac{dU}{dt} &=&\displaystyle   m\,\ch(2V)-m-\varepsilon   \\[8pt]
\displaystyle  \frac{dV}{dt} &=&\displaystyle   \bm{u}(t) - m \sh(2V)  
 \end{array} 
 \right. 
\feq
The system $\mb{S}(\varepsilon,m,\mu)$ is composed of the one-dimensional system
 \beq 
  \bm{F}(m,\mu)\quad \quad
\begin{array}{lcl}
\label{eq:VR}
\displaystyle  \frac{dV}{dt} &=&\displaystyle   \bm{u}(t) - m \sh(2V)   
 \end{array} 
\feq
which then gives the solution of $U$,
\beq \label{UdeR}
\begin{array}{lcl}
\displaystyle  U(t)  &=&\displaystyle  U_0 + \int_0^t \left(m\,\ch(2V(s)) -m-\varepsilon \right)ds  
 \end{array} 
\feq
\begin{remarque}
\label{rem:semiMarkov}
{\rm The process $(\bm{u}(t))_{t \geq 0}$ is a so-called \emph{semi-Markov process}. It is in general not a Markov process, but the sequence of post-jump locations make a Markov chain, which explains the name semi-Markov. In the particular case where $\mu$ is an exponential law (and only in this case),  $(\bm{u}(t))_{t \geq 0}$  is a Markov process, as well as the processes $(V, U, \bm{u})$ and $(V, \bm{u})$. These two latter processes are called \emph{Piecewise Deterministic Markov Processes} (see Section \ref{sec:PDMP} below).
}
\end{remarque}
$\,$

Recall that we say that a sequence of random variables $Y_n$ converges in distribution (or in law) to a variable $Y_{\infty}$ if for all bounded continuous function $f$, $\mathbb{E}(f(Y_n)) \to \mathbb{E}(f(Y_{\infty}))$ when $n$ goes to infinity.  For all $n \geq 0$, set $\hat V_n = V(T_n)$.

\begin{lemme}
\label{lem:V2k}
There exist random variables $V_{\infty}^{-} = V_{\infty}^{-}(m, \mu) $ and  $V_{\infty}^{+} = V_{\infty}^{+}(m, \mu)$ such that
\begin{enumerate}
\item $V_{\infty}^{+}$ and $V_{\infty}^{-}$ lie almost surely in $[ V_m^-, V_m^+]$; where $V_m^+$ and $V_m^-$ are given by \eqref{Vm+Vm-};
\item $\hat V_{2n}$ and $\hat V_{2n+1}$ converge in distribution to $V_{\infty}^{-}$ and $V_{\infty}^{+}$ respectively;
\item Let $S$ be a random variable with law $\mu$, independent from $V_{\infty}^-$ and  $V_{\infty}^+$. Then, with probability $1$, for all bounded measurable function $f : \mathbb{R} \times \mathbb{R}_+ \to \mathbb{R}$, 
\[
\lim_{n \to \infty} \frac{1}{n} \sum_{k=0}^{n-1} f ( \hat V_{2k}, S_{2k+1} ) = \mathbb{E}[f(V_{\infty}^{-}, S)]; \quad  \lim_{n \to \infty} \frac{1}{n} \sum_{k=0}^{n-1} f ( \hat V_{2k+1}, S_{2k+2} ) = \mathbb{E}[f(V_{\infty}^{+}, S)]
\]

\end{enumerate}
\end{lemme}

\begin{remarque}{\rm
If $\mu = \delta_T$ is the Dirac mass at $T$, $V_{\infty}^{-}(m, \delta_T) = P_{m,T}(0)$ and $V_{\infty}^{+}(m, \delta_T) = P_{m,T}(T)$, where $P_{m,T}$ is the unique periodic solution of the system $F(m,T)$ (see Equation \eqref{SV}) granted by Proposition  \ref{prop1}.
}
\end{remarque}
$\,$\\
This lemma proved in Appendix \ref{app:lemV2k}, tells us that the location of $V$ after an even number of jumps is asymptotically close, in distribution, to a variable $V^{-}_{\infty}$, and to a variable $V^{+}_{\infty}$ after an odd number of jumps. This enable us to give the asymptotic growth rate of $U$, as in Proposition \ref{prop2}. Recall that  $\varphi^+_t(v)$ and $\varphi^-_t(v)$ are the solutions to $F_m^+$ and $F_m^-$ at time $t \geq 0$, starting from $v$ at time $0$, respectively.
\begin{proposition}
\label{prop:Delta_random_time}
Let $S$ be a random variable with law $\mu$, independent from $V_{\infty}^{+}$ and $V_{\infty}^{-}$. Set 
\begin{equation}
\label{eq:Delta_random_time}
\bm{\Delta}( m, \mu) = \frac{\mathbb{E} \left( \int_0^S m ( \ch(2 \varphi_s^+ ( V_{\infty}^-) ) - 1 ) ds \right) + \mathbb{E} \left( \int_0^S m ( \ch( 2\varphi_s^- ( V_{\infty}^+) ) - 1 ) ds \right) }{2 \mathbb{E}(S)}
\end{equation}
Then, for all initial condition $V(0), U(0)$, one has, with probability one, 
\[
\lim_{ t \to \infty} \frac{U(t)}{t} =  \bm{\Delta}( \varepsilon, m, \mu):= \bm{\Delta}( m, \mu) - \varepsilon.
\]
\end{proposition}
The proof of Proposition \ref{prop:Delta_random_time}, given in Appendix \ref{App:Delta_random_switch}, is similar to the proof of Proposition \ref{prop2}, using probabilistic tools and the law of large numbers given by the third point of Lemma \ref{lem:V2k}.

\begin{remarque}{\rm 
If $\mu = \delta_T$ is the Dirac mass at $T$, we have by definition of $P_{m,T}$ that for all $s \in [0,T]$, $\varphi_s^+ ( V_{\infty}^-)  = P_{m,T}(s)$ and $\varphi_s^- ( V_{\infty}^+) = P_{m,T}(T+s)$. Thus, 
\[
\mathbb{E} \left( \int_0^S m ( \ch(2 \varphi_s^+ ( V_{\infty}^-) ) - 1 ) ds \right) = \int_0^T m ( \ch(2 P_{m,T}(s) ) - 1 ) ds
\]
and
\begin{align*}
\mathbb{E} \left( \int_0^S m ( \ch( 2\varphi_s^- ( V_{\infty}^+) ) - 1 ) ds \right) & = \int_0^T m ( \ch(2 P_{m,T}(s+T) ) - 1 ) ds\\
& = \int_T^{2T} m ( \ch(2 P_{m,T}(s) ) - 1 ) ds,
\end{align*}
so that
\[
\bm{\Delta}( \varepsilon, m, \delta_T) = \frac{1}{2T} \int_0^{2T} m ( \ch(2 P_{m,T}(s) ) - 1 - \varepsilon) ds = \Delta(\varepsilon,m,T).
\]
Therefore, we retrieve the growth rate computed in the previous section.
}
\end{remarque}
$\,$
\begin{remarque}{\rm 
\label{rem:additionalmu}
Under some additional assumptions on the law $\mu$, it is possible to prove that $V(t)$ converges in distribution to a variable $V_{\infty}$ as $t$ goes to infinity (see the forthcoming paper \cite{HS22} for general conditions and Section \ref{sec:PDMP} below for a particular case). In that case, we can express $\bm{\Delta}( \varepsilon, m, \mu)$ as 
\[
\bm{\Delta}( \varepsilon, m, \mu) = \mathbb{E} \left[m \left(\ch(2 V_{\infty}) - 1\right) \right]- \varepsilon.
\]
}
\end{remarque}
\subsubsection{The particular case of PDMP}
\label{sec:PDMP}
In this section, we detail the particular case where $\mu$ is an exponential law with parameter $\sigma$, i.e., $\mu$ is absolutely continuous with respect to the Lebesgue measure with density defined on $\mathbb{R}_+$ by $g(x) = \sigma e^{ - \sigma x}$.  Since $\mathbb{E}(S_1) = \frac{1}{\sigma}$, we rather use the parametrization $T = 1/\sigma$ in the sequel.
In this situation of an exponential law, as noticed in Remark \ref{rem:semiMarkov}, the process $(V_t, \bm{u}(t))_{t \geq 0}$ is a Piecewise Deterministic Markov process. In addition, some explicit computations are made possible in that case.

First, it can be proven easily, using e.g.  \cite[Theorem 4.6]{BMZIHP} that $V(t)$ converges in distribution to a random variable $V_{\infty}$, whose law admits a density with respect to the Lebesgue measure. In addition, this density is explicitly computable, and given by  (see e.g \cite[Proposition 3.12]{FGR09} for the general formula)	
\begin{equation}
\label{eq:density}
\bm{\rho}_{m,T}(v) =  C(m) \left( \frac{1}{|F_m^+(v)|} + \frac{1}{|F_m^-(v)|} \right)  \left( \frac{e^{V_m^+} - e^v}{e^v + e^{V_m^-} }  \frac{e^v - e^{V_m^-} }{e^v  + e^{V_m^+} } \right)^{\frac{1}{2T \sqrt{m^2 + 1}}};
\end{equation}
for all $v \in [V_m^-, V_m^+]$, where $C(m)$ is a normalisation constant. Moreover, the following strong law of large numbers is satisfied. For all bounded measurable function $f : [V_m^-, V_m^+] \mapsto \mathbb{R}$
\[
\lim_{t \to \infty} \frac{1}{t}\int_0^t f(V(s)) ds = \mathbb{E}( f(V_{\infty})).
\]
This entails, as noticed in Remark \ref{rem:additionalmu}, that
\[
\lim_{t \to \infty} \frac{U(t)}{t} = m \left(\ch(2 V_{\infty}) - 1\right) - \varepsilon=\bm{\Delta}( \varepsilon, m, \mu),
\]
which can be rewritten as
\[
\bm{\Delta}( \varepsilon, m, \mu) = \int_{[V_m^-, V_m^+]} m \left(\ch(2 v) - 1\right) \bm{\rho}_{m,T}(v) dv.
\]
\begin{remarque}{\rm
\label{rem:explicite}
From the explicit expression 
\eqref{eq:density} of  $\rho_{m,T}$, it is possible to prove that, for fixed $m$ and $T$, there exist constants $C_-(m,T), C_+(m,T)$ such that,  as $v \to V_m^+$,
\[
\rho_{m,T}(v) \sim C_+(m,T) (e^{V_m^+} - e^v )^{\frac{1}{2T \sqrt{1+m^2} } - 1},
\] 
while as $v \to V_m^-$,
\[
\rho_{m,T}(v) \sim C_-(m,T) (e^{v} - e^{V_m^-} )^{\frac{1}{2T \sqrt{1+m^2} } - 1}.
\]
In particular, $\rho_{m,T}$ is bounded in  neighbourhoods of $V_m^+$ and $V_m^-$ if and only if $ 1 \geq  2 T \sqrt{m^2 + 1}$. This condition is consistent with the following heuristic: if $T$ is large, the environment does not switch often, and the process $V$ follows the vector fields $F_m^{+}$ and $F_m^-$ for a long time, and thus spend a large amount of time close to the equilibria $V_m^+$ and $V_m^{-}$. Hence, for large $T$, one expects that the distribution $\Pi_{m,T}$ give a lot of mass near $V_m^+$ and $V_m^{-}$. On the contrary, if $T$ is small, the environment switches frequently, and the process $V$ spend most of time in the middle of the interval $[V_m^-, V_m^+]$, and therefore one expects the distribution $\Pi_{m,T}$ to vanish at the extremity of the interval.
}
\end{remarque}
\subsubsection{When the mean switching time goes to infinity}
Let us assume that we have a family of law $(\mu^{(T)})_{ T > 0}$ such that, for all $T > 0$; $\int_{[0, \infty)} t \mu^{(T)} (dt) = T$. In  other words, the mean time spent in each environment is $T$.  We now prove that, when $m$ is fixed and $T$ goes to infinity, the asymptotic of $\bm{\Delta}( \varepsilon, m,  \mu^{(T)})$  is the same as in the periodic case. In particular, for $m$ small enough, one can choose $T$ large enough so that $\bm{\Delta}( \varepsilon, m,  \mu^{(T)}) > 0$ and there is inflation. This comes from the fact that, as $T$ goes to infinity, the time spent in each environment is large enough so that $\varphi_s^0(v)$ and $\varphi_s^1(v)$ become, uniformly in $v$ in a compact interval, arbitrarily close from $V_m^{+}$ and $V_m^{-}$, respectively.

\begin{proposition}
\label{prop:lim_Delta_switching_infty}
For fixed $m > 0$, 
\[
\lim_{ T \to \infty} \bm{\Delta}( \varepsilon, m,  \mu^{(T)}) = \sqrt{1 + m^2} - m - \varepsilon.
\]
In particular, whenever $m < \frac{1 - \varepsilon^2}{2\varepsilon}$, for $T$ large enough, $\bm{\Delta}( \varepsilon, m,  \mu^{(T)}) > 0$ and there is inflation.
\end{proposition}
\begin{proof}
The proof is similar to the proof of Cesaro's Lemma. Let $S^{(T)}$ be a random variable with law $\mu^{(T)}$. We claim that for all continuous function $g : [V_m^-, V_m^+] \to \mathbb{R}$, 
\[
\lim_{T \to \infty} \sup_{ v \in [V_m^-, V_m^+]} \big| \frac{\E\left(\int_0^{S^{(T)}} g( \varphi_r^+(v)) dr\right)}{\E(S^{(T)})} - g(V_m^+) \big| = 0,
\]
and similarly with $V_m^-$ instead of $V_m^+$ when $\varphi^+$ is replaced by $\varphi^-$ in the integral. Applying this result to the function $g(v) = m (\ch(2v) - 1)$ and using formula \eqref{eq:Delta_random_time} proves the proposition since $g(V_m^+) = g(V_m^-) = \sqrt{1 + m^2} - m$. We now prove the claim. Since $V_m^+$ is globally attractive for the flow $\varphi^+$, for all $\varepsilon > 0$ there exists $M > 0$ such that, for all $r \geq M$, $\sup_{v \in [V_m^-, V_m^+]} | \varphi_r(v) - V_m^+ | \leq \varepsilon$. Since $g$ is uniformly continuous on $[V_m^-, V_m^+]$, this entails that for $M$ large enough and $r \geq M$, $\sup_{v \in [V_m^-, V_m^+]} | g(\varphi_r(v)) - g(V_m^+) | \leq \varepsilon$. Hence, for all $v \in [V_m^-, V_m^+]$, 
\begin{align*}
\big| \frac{\E(\int_0^{S^{(T)}} g( \varphi_r^+(v) )dr)}{\E(S^{(T)})} - g(V_m^+) \big| & \leq \frac{\mathbb{E}( \int_0^{S^{(T)} \wedge M} | g(\varphi_r^+(v)) - g(V_m^+)| dr)}{\E(S^{(T)})} \\
& \quad   +  \frac{\mathbb{E}( \int_M^{S^{(T)}} | g(\varphi_r^+(v)) - g(V_m^+)| dr\1_{S^{(T)} > M} )}{\E(S^{(T)})} \\
& \leq \frac{2 M \|g\|_{\infty}}{T} + \varepsilon,
\end{align*}
where we have used that $\E(S^{(T)}) = T$. This proves the claim.
\end{proof}

\subsubsection{Random choice of $(\pm 1)$.}
\label{sec:random_pm}
In this section, we study another type of random process linked to the $(\pm 1)$ model. We assume that, in each patch, after each $T$ units of time, we select at random, independently from anything else, whether the growth rate within the patch will be $1 - \varepsilon$ or $-(1 + \varepsilon)$ for the next $T$ units of time. More formally, we consider the system such that; for all $n \geq 0$, for all $t \in [nT,(n+1)T)$, 
\begin{equation}
\label{system1RandomChoices}
\bm{\Sigma}(\varepsilon,p_1,p_2,m,T) \quad \left\{
\begin{array}{lcrcl}
\displaystyle \frac{dx_1}{dt}&=&\displaystyle (z_1^n - \varepsilon) x_1&+&m(x_2-x_1)\\[2mm]
\displaystyle\frac{dx_2}{dt}&=&\displaystyle (z_2^n -  \varepsilon) x_2&+&m(x_1-x_2)
\end{array}
\right.
\end{equation}
where $(z_1^n)_{n \geq 0}$ and  $(z_2^n)_{n \geq 0}$  are independent  sequences of i.i.d. random variables with values in $\{-1,1\}$ such that  for $i = 1, 2$, $\mathbb{P} ( z_i^n = 1) = p_i \in (0,1)$. Note in particular that the growth rates of the two patches are totally uncorrelated and that there is no temporal autocorrelation for the value of the growth rate within a given patch. In the $U - V$ variables, the system becomes, for  $t \in [nT,(n+1)T)$, 
\begin{equation}\label{aleatenUV}
\left\{
\begin{array}{lcl}
\displaystyle  \frac{dU}{dt} &=&\displaystyle \frac{z_1^n + z_2^n}{2}  + m\, \ch(2V)-m - \varepsilon   \\[8pt]
\displaystyle  \frac{dV}{dt} &=& \displaystyle \frac{z_1^n - z_2^n}{2}-m\,\sh(2V)
 \end{array} 
 \right. 
\end{equation}
Note that now, $V$ is switching between three autonomous system : $F_m^+$, $F_m^-$ and $F_m^0$, where
\beq \label{S0}
 F_m^0 \quad \quad\quad\quad \quad\quad \quad \left\{\frac{dV}{dt} =  -m\,\sh(2V) \right. \quad 
\feq
We denote by $\varphi^0$ the flow associated to $F_m^0$.
For $n\geq 0$, we let $\hat V_n = V(nT)$. As in the previous section, we can precise the asymptotic behavior of $\hat V_n$.
\begin{lemme}
The followings hold true:
\begin{enumerate}
\item The sequence $(\hat V_n)_{n \geq 0}$ is a Markov chain;
\item $\hat V_n$  converges in distribution to a random variable $\hat V_{\infty}$ which lies almost surely in $[V_m^-, V_m^+]$;
\item For $k \geq 0$, let $h_k = \frac{z_1^k-z_2^k}{2}$ and let $\bm{h}$ a random variable independent from $\hat V_{\infty}$, with the law of $h_1$. Then, with probability $1$, for all bounded measurable function $f : [V_m^-, V_m^+] \times \{ -1, 0, 1\} \to \mathbb{R}$,
\[
\lim_{n \to \infty} \frac{1}{n} \sum_{k=0}^{n-1} f(\hat V_k, h_k) = \mathbb{E}( f( \hat V_{\infty}, \bm{h}).
\]
\end{enumerate}

\end{lemme}
Now, similarly to the proof of Proposition \ref{prop:Delta_random_time}, we can use the previous lemma to show that there exists an asymptotic growth rate for $U$.
\begin{proposition}
\label{prop:Delta_random_choice}
Let $\bm{h}$ be a random variable, independent from $\hat V_{\infty}$, with the law of $\frac{z_1^1-z_2^1}{2}$. Set 
\[
\bm{\Delta}(p_1, p_2, m, T) = \frac{\mathbb{E}(\int_0^T m \left( \ch( 2 \varphi^{\bm{h}}_r(\hat V_{\infty})) - 1\right) dr) }{T}+ p_1 + p_2 - 1
\]
Then, 
\[
\lim_{t \to \infty} \frac{U(t)}{t} = \bm{\Delta}(p_1, p_2, m, T, \varepsilon) = \bm{\Delta}(p_1, p_2, m, T) - \varepsilon.
\]
\end{proposition}
 The proof of Proposition \ref{prop:Delta_random_choice} is very similar to the proof of Proposition \ref{prop:Delta_random_time} and left to the reader.
  \begin{figure}
  \begin{center}
 \includegraphics[width=1\textwidth]{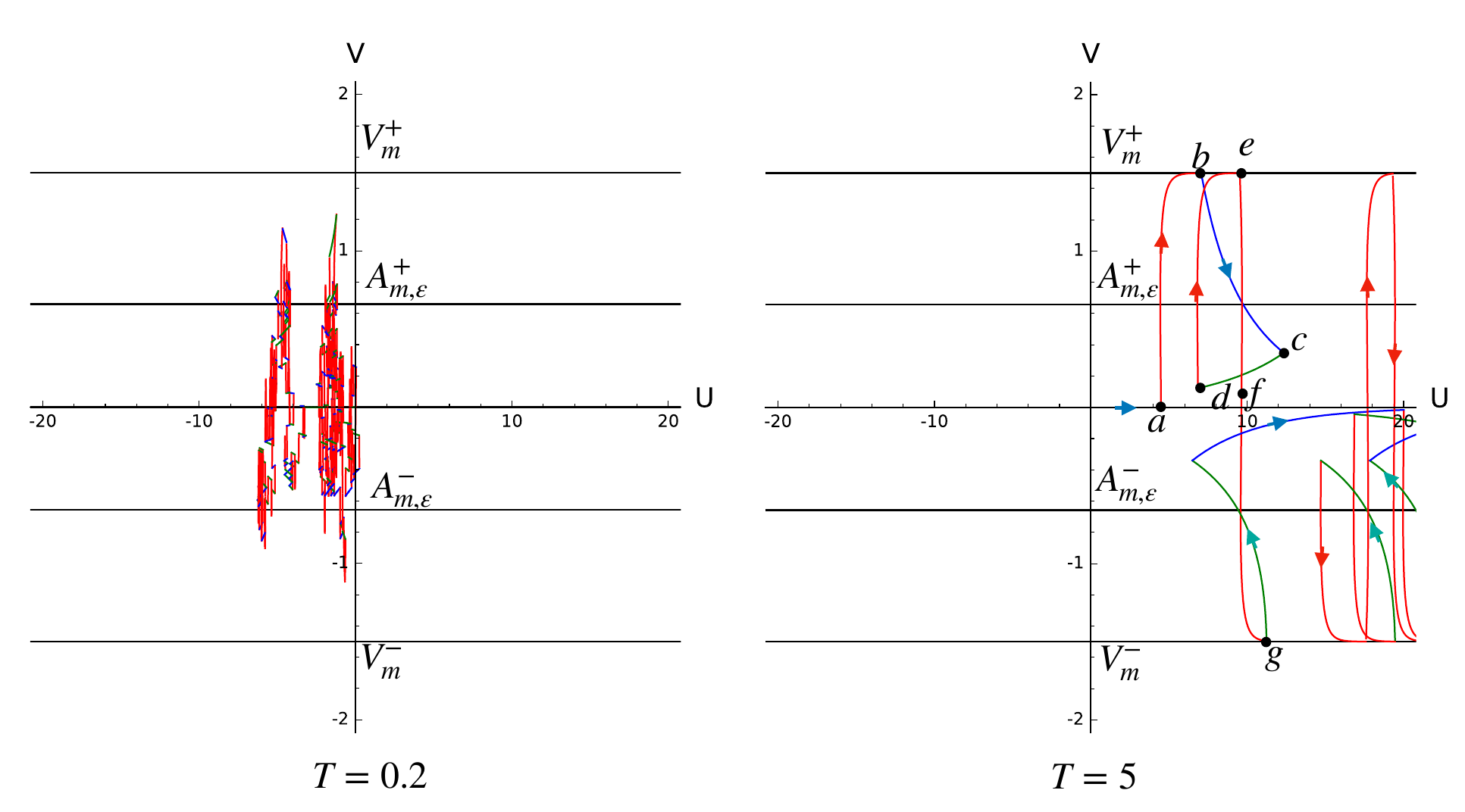}
 \caption{System  \eqref{aleatenUV}. On the left one realization with duration between switches = 0.2 ; on the right duration = 5. More comments in the text. }\label{PMPPMM}
 \end{center}
 \end{figure}
We illustrate this proposition by the simulations of figure \ref{PMPPMM} which is are counterpart of those of figure \ref{periodic} in the deterministic periodic case. We have considered the system \eqref{aleatenUV} with $p_1 = p_2 = \frac{1}{2}$ and $m = 0.1$ and  $\eps = 0.1$. On the right one sees a realization of the process with $T=5$ ; one sees:
\bito
\item From $0$ to $a$ : $z_1 = z_2 = +1$. The abondance of population is increasing on both patches; trajectory in blue.
\item From $a$ to $b$ : $z_1 = 1,\;z_2 = -1$. Site  1 is favorable, site  2 unfavorable.
\item From $b$ to $c$ : $z_1 = z_2 = +1$. The abondance of population is increasing on both patches.
\item From $c$ to $d$ : $z_1 = z_2 = - 1$. The abondance of population is decreasing on both patches; trajectory in green.
 \item From $d$ to $e$ : $z_1 = 1,\;z_2 = -1$. Site  1 is favorabble, site  2 unfavorable.
 \item From $e$ to $f$ : $z_1 = -1,\;z_2 = +1$. Site  1 is unfavorable, site  2 favorable.
 \item From $f$ to $g$ : $z_1 = -1,\;z_2 = +1$. Site  1 is favorable, site  2 unfavorable.
 \item $\dots$
 \fit
We see that, with respect to $U$, green and blue trajectories almost compensate while, since $T$ is large enough, the red trajectories spent enough time in the strips $\mathbb{R} \times [A^+_{m,\varepsilon},V^+_m] $ and $\mathbb{R} \times [V^-_m,A^-_{m,\varepsilon},] $ where $U$ is increasing. on the left simulation we chose $T = 0.2$  which gives little chance to the trajectory to reach the strips where $U$ is increasing.  
 
 In the situation of random switching time studied previously, we have shown that, provided $m$ is small enough and the mean switching time is large enough, there is inflation. The present case of random choice of $(\pm 1)$ is a bit different. Depending on the parameters $p_1$, $p_2$ and $\varepsilon$, it might happen that inflation never occurs, whatever the values of $m$ and $T$ are. We use the notations of  Katriel to state this result.
\begin{proposition}
\label{prop:chi_random_choice}
Let 
\[
\chi = \chi(p_1, p_2,\varepsilon) = p_1 ( 1 - p_2) + p_2 ( 1 - p_1) +p_1 + p_2 - 1  - \varepsilon
\]
Then
\begin{itemize}
\item If $\chi < 0$, then for all $(m,T)$, $\Delta(p_1, p_2, m, T, \varepsilon) < 0$, and there is no inflation
\item If $\chi > 0$,  there exists $m^*( \varepsilon)$ such that, for all $m \in (0, m^*(\varepsilon))$, there exists $T^*(m)$ such that, for all $T \geq T^*(m)$,  $\Delta(p_1, p_2, m, T, \varepsilon) > 0$ and there is inflation.
\end{itemize}
The second assertion is a consequence of the fact that, for all $m > 0$, 
\begin{align*}
\lim_{T \to \infty} \Delta(p_1, p_2, m, T, \varepsilon) & = \left[p_1 ( 1 - p_2) + p_2 ( 1 - p_1)\right]( \sqrt{1 + m^2} - m) +p_1 + p_2 - 1  - \varepsilon\\
& = \chi - \left[p_1 ( 1 - p_2) + p_2 ( 1 - p_1)\right]m + \po (m).
\end{align*}
\end{proposition}
\begin{remarque}{\rm
Note that
\[
\chi =  \left[ p_1 ( 1 - p_2) + p_2 ( 1 - p_1) + p_1 p_2\right] ( 1 - \varepsilon) - \left[ (1-p_1)(1-p_2) \right](1+\varepsilon).
\]
The term $p_1 ( 1 - p_2) + p_2 ( 1 - p_1) + p_1 p_2$ is the proportion of time where a least one patch is favourable, while the term $ (1-p_1)(1-p_2) $ is the proportion of time where the two patches are unfavourable. In particular, $\chi > 0$ if and only if
\[
\frac{p_1 ( 1 - p_2) + p_2 ( 1 - p_1) + p_1 p_2}{(1-p_1)(1-p_2)} > \frac{1+\varepsilon}{ 1 - \varepsilon},
\]
that is, the ratio of the time in favorable states and the time in unfavorable state is higher than the ratio of the rates of decrease and of increase.
}
\end{remarque}
The proof of the first item of Proposition \ref{prop:chi_random_choice} is similar to the proof of the result of  Katriel, with the use of the law of large numbers. It is remarkably simple, as we detail now
\begin{proof}
For all $t \in [nT, n(T+1))$, we have
\begin{align*}
\frac{d( x_1 + x_2)}{dt} & = (z_1^n - \varepsilon) x_1 + (z_2^n - \varepsilon) x_2\\
& \leq \left( \max(z_1^n, z_2^n) - \varepsilon\right)( x_1 + x_2).
\end{align*}
This implies that 
\[
\ln( (x_1 + x_2)((n+1)T)) \leq \left( \max(z_1^n, z_2^n) - \varepsilon\right) T + \ln((x_1 + x_2)(nT)),
\]
and thus for all $n \geq 1$,
\[
 \frac{\ln((x_1 + x_2)(nT))}{nT} \leq- \varepsilon +  \frac{1}{n} \sum_{k=0}^{n-1} \max(z_1^k, z_2^k).
\]
Since the sequence  $(\max(z_1^k, z_2^k))_{k \geq 0}$ is i.i.d., the strong law of large numbers implies that with probability $1$,
\begin{align*}
\lim_{ n \to \infty} \frac{1}{n} \sum_{k=0}^{n-1} \max(z_1^k, z_2^k)  & = \mathbb{E}( \max(z_1^k, z_2^k))\\
& = p_1 ( 1 - p_2) + p_2 ( 1 - p_1) +p_1p_2 - (1-p_1)(1 - p_2)
\end{align*}
Hence, 
\[
\lim_{n \to \infty} \frac{\ln((x_1 + x_2)(nT))}{nT} \leq \chi,
\]
and this entails the first point of the proposition.

The proof of the second point is very similar to the proof of Proposition \ref{prop:lim_Delta_switching_infty}, where we also use that, as $t$ goes to infinity, $\varphi^0(v) \to 0$, uniformly in $v \in [V_m^-, V_m^+]$. Thus, the proof is omitted.
\end{proof}
\begin{remarque}{\rm 
We could also have considered the case where $z_1^n$ is not necessarily independent from $z_2^n$. In that case, we give the law  $\mathbf{p} = (p_{1,1}, p_{1,-1}, p_{-1,1}, p_{-1,-1})$ of the couple $Z^n = (z_1^n, z_2^n)$ : for $h, h' \in \{-1, +1\}$,
\[
\mathbb{P}\left(Z^n = (h,h') \right) = p_{h,h'},
\]
for some $p_{1,1}, p_{1,-1}, p_{-1,1}, p_{-1,-1}$ that sum to $1$. The formula become
\[
\lim_{T \to \infty} \Delta(\mathbf{p}, m, T, \varepsilon) = \left[p_{+,-} + p_{-,+}\right]( \sqrt{1 + m^2} - m) + p_{++} - p_{--} - \varepsilon
\]
\[
\chi(\mathbf{p},\varepsilon) = [ p_{+,-} + p_{-,+} + p_{+,+} ] ( 1 - \varepsilon)  - p_{-,-} ( 1 + \varepsilon)
\] 
and $\chi > 0 $ if and only if, either $p_{-,-} = 0$ or
\begin{equation}\label{condition}
\frac{p_{+,-} + p_{-,+} + p_{+,+}}{p_{-,-}}  > \frac{1+\varepsilon}{ 1 - \varepsilon}.
\end{equation}
We notice that when both $p_{1,1}$ and $p_{-1,-1}$ are null, i.e.  the case where the patches are always in opposite growth,
\[
\lim_{T \to \infty} \Delta(\mathbf{p}, m, T, \varepsilon) = ( \sqrt{1 + m^2} - m)  - \varepsilon,
\]
which is the same limit as in the periodic case with alternating $-+$ and $+-$. 
}
\end{remarque}

 \subsubsection{Link with the top Lyapunov exponent}\label{sec:TLE}
Let $X_t = (x_1(t), x_2(t))$ the solution to $\bm{\Sigma}(m, \varepsilon, \mu)$. With the notation of Section \ref{variables initiales}, one can rewrite $\bm{\Sigma}(m, \varepsilon, T)$ as
\beq
\frac{d X_t}{dt} = M_{\varepsilon, m}^{\mb{u}(t)} X_t.
\feq
Since $\mb{u}(t)$ is constant on interval $[T_n, T_{n+1}[$ of length $S_{n+1}$, one has
\[
X(T_{n+1}) = e^{ S_{n+1} M_{\varepsilon,m}^{\mb{u}(T_n)} } X(T_n)
\]
Hence, setting $\hat X_n = X(T_n)$ and $\hat{\mb{u}}_n = \mb{u}(T_n)$, one can write $\hat X_n$ as the product
\[
\hat X_n = \left( \prod_{i=0}^{n-1} B_i \right) \hat X_0,
\]
where $B_i$ is the random matrix $e^{ S_{i+1} M_{\varepsilon,m}^{\hat{\mb{u}}_i} }$. Since the sequence $(\hat{\mb{u}}_n)_{n \geq 0}$ forms a Markov chain and since the $(S_n)_{n \geq 1}$ are i.i.d., one can prove (see \cite[Propsition 3.8]{CM19}) that the classical Oseledet's Multiplicative ergodic theorem can be applied. According to this theorem, the limit\footnote{Here $\| \cdot \|$ stands for the euclidian norm on $\mathbb{R}^2$, but the limit is independent of the choice of the norm.}
\[
\lim_{n \to \infty} \frac{1}{n} \ln \| \hat X_n \|
\]
exists, and can take at most two different values $\lambda_1 \geq \lambda_2$, called \emph{Lyapunov exponent} (see e.g Chapter 1.4 in \cite{these}). 
Since the matrices $M_{\varepsilon, m}^{h}$ are irreducible and \textit{ Metzler }, ie have non-negative off-diagonal coefficients, a random version of Perron - Frobenius Theorem (see \cite{ADG94}), and Proposition 2.13 in \cite{BS19}) implies that the top  Lyapunov exponent $\lambda_1$,  is such that, for all $X_0 \in \mathbb{R}_+^2 \setminus \{0 \}$, almost surely, 
\[
\lim_{n \to \infty} \frac{1}{n} \ln \| \hat X_n \| = \lambda_1.
\]
Moreover, by Proposition 3.4 in \cite{CM19}, we can define the growth rate $\bm{\Lambda}(\varepsilon, m, \mu)$ of the continuous-time model and related it to the Lypaunov exponent of the discrete-time model:
\[
\bm{\Lambda}(\varepsilon, m, \mu) := \lim_{t \to \infty} \frac{1}{t} \ln \|X_t \| = \frac{1}{\mathbb{E}(S_1)} \lim_{n \to \infty} \frac{1}{n} \ln \| \hat X_n \| =\frac{\lambda_1}{\mathbb{E}(S_1)}.
\]
Obviously, $\bm{\Lambda}(\varepsilon, m, \mu)$ and $\bm{\Delta}(\varepsilon, m, T)$ are linked. Indeed, note that the compact set $[V_m^-, V_m^+]$ is positively invariant for $V = \frac{1}{2}(\ln(x_1) - \ln(x_2))$ and attracts all trajectories. Hence, for all initial condition $(x_1(0), x_2(0))$, there exists a time $t_0$ such that, for all $t \geq t_0$, $V(t) \in [V_m^-, V_m^+]$. In particular, for $t \geq t_0$; 
\[
e^{2V_m^-} \leq \frac{x_1(t)}{x_2(t)} \leq e^{2V_m^+}
\]
This yields 
\[
\left( e^{2V_m^-}  + e^{ - 2V_m^+} \right) x_1(t)x_2(t) \leq x_1(t)^2 + x_2(t)^2 \leq  \left( e^{-2V_m^-}  + e^{ 2V_m^+} \right) x_1(t)x_2(t)
\]
Taking the logarithm and sending $t$ to infinity proves the following:
\begin{proposition}
\label{prop:lyapDelta}
One has
\[
\bm{\Delta}(\varepsilon, m, \mu) =  \bm{\Lambda}(\varepsilon, m, \mu).
\]
\end{proposition}
\begin{remarque}{\rm 
For the system $\bm{\Sigma}(\varepsilon,p_1,p_2,m,T)$ considered in Section \ref{sec:random_pm}, one can prove similarly that $\hat X_n = X(nT)$ is described by a random product of matrices, and that there exists a top Lyapunov exponent $\lambda_1$ such that $\bm{\Delta}(p_1, p_2, m, T, \varepsilon) = \frac{\lambda_1}{T}$.
}
 \end{remarque}
\section{Some extensions to more complex situations} \label{morecomplex}
\subsection{The case of partial phase shift}\label{phaseshift}
 In the preceding section we considered the case where the two patches where always in opposite conditions during the whole period $2T$. A more realistic situation is when the two patches are ruled by the same periodic environment $r(t)$ shifted of $\varphi T$ with $\varphi \in (0,1)$.

Hence, we consider the system
$\Sigma(r_1,r_2,m,T)$,
given by \eqref{BLSS1}, where $r_1(t)$ and $r_2(t)$ are the $2T$-periodic functions defined by
$$
r_1(t)=\left\{
\begin{array}{rll}
r&\mbox{if}&t\in[0,T)\\
-d&\mbox{if}&t\in[T,2T)
\end{array}
\right.
\quad
r_2(t)=r_1(t-\varphi T).
$$
{As in \eqref{+r-d}, we assume that $d>r>0$ which means that the mean of the growth rate on each patch is negative (each patch is a sink). We have :
\[
\chi = \frac{1}{2T} \int_0^{2T} \max( r_1(s), r_2(s) ) ds = r \frac{1+ \varphi}{2} - d \frac{1 - \varphi}{2}
\]
where $\frac{1+ \varphi}{2}$ is the proportion of time where at least one of the patches is increasing, while $\frac{1 - \varphi}{2}$ is the proportion of time where both patches are decreasing. Hence, $\chi > 0$ if and only if  
\[
\frac{1 + \varphi}{ 1 - \varphi} > \frac{d}{r},
\]
as in the stochastic $(\pm 1)$ model, see \eqref{condition}.

Let us consider now the special case where $r=1-\varepsilon$ and $d=1+\varepsilon$ corresponding to $(\pm 1)$ model. We have
$
\chi = \varphi-\varepsilon
$.
Therefore $\chi>0$ if and only if $\varphi>\varepsilon$. For illustration let us  plot the function $(m,T)\mapsto \Delta(m,T)$.
 As previously we have the two systems~:
  \beq \label{pm}
 \begin{array}{lcl}
	\Sigma^{+-}(\varepsilon,m,T,\varphi)&\quad\quad &
	\left\{
	\begin{array}{lcl}
		\displaystyle  \frac{dx_1}{dt}& = &(+1-\varepsilon)x_1+m(x_2-x_1) \\[10pt]
		\displaystyle  \frac{dx_2}{dt}&=&  (-1-\varepsilon)x_2+m(x_1-x_2) 
 	\end{array} 
	 \right.
 \end{array}
 \feq
and :
 \beq  \label{mp}
 \begin{array}{lcl}
\Sigma^{-+}(\varepsilon,m,T,\varphi)&\quad\quad &
\left\{
	\begin{array}{lcl}
	\displaystyle  \frac{dx_1}{dt}& = &(-1-\varepsilon)x_1+m(x_2-x_1) \\[10pt]
	\displaystyle  \frac{dx_2}{dt}&=&  (+1-\varepsilon)x_2+m(x_1-x_2) 
	 \end{array} 
 \right.
 \end{array}
 \feq
to which we add:
  \beq \label{pp}
 \begin{array}{lcl}
	\Sigma^{++}(\varepsilon,m,T,\varphi)&\quad\quad &
	\left\{
	\begin{array}{lcl}
		\displaystyle  \frac{dx_1}{dt}& = &(+1-\varepsilon)x_1+m(x_2-x_1) \\[10pt]
		\displaystyle  \frac{dx_2}{dt}&=&  (+1-\varepsilon)x_2+m(x_1-x_2) 
 	\end{array} 
	 \right.
 \end{array}
 \feq
and :
 \beq  \label{mm}
 \begin{array}{lcl}
\Sigma^{--}(\varepsilon,m,T,\varphi)&\quad\quad &
\left\{
	\begin{array}{lcl}
	\displaystyle  \frac{dx_1}{dt}& = &(-1-\varepsilon)x_1+m(x_2-x_1) \\[10pt]
	\displaystyle  \frac{dx_2}{dt}&=&  (-1-\varepsilon)x_2+m(x_1-x_2) 
	 \end{array} 
 \right.
 \end{array}
 \feq
 We switch from one system to the other according to the following scheme :
 \beq
 \begin{array}{|c||c|c|c|c|}
 \hline
 t \in&[0,\quad \varphi T [&[\varphi T,\quad T[&[T,T(1+\varphi),[&[T(1+\varphi),2T[\\
 \hline
 \hline
 \Sigma&+-&++&-+&--\\
 \hline
 \end{array}
 \feq 

Using notations similar to those we used in Subsection \ref{variables initiales}, let us define~:
 $$M^{+-}_{\varepsilon,m} = 
\left[
\begin{array}{cc}
1-m-\varepsilon &+m\\
+m & -1-m -\varepsilon
\end{array}
\right]
\quad
M^{-+}_{\varepsilon,m} = 
\left[
\begin{array}{cc}
-1-m-\varepsilon &+m\\
+m &1-m -\varepsilon 
\end{array}
\right]
$$
$$M^{++}_{\varepsilon,m} = 
\left[
\begin{array}{cc}
1-m-\varepsilon &+m\\
+m & +1-m -\varepsilon
\end{array}
\right]
\quad
M^{--}_{\varepsilon,m} = 
\left[
\begin{array}{cc}
-1-m-\varepsilon &+m\\
+m &-1-m -\varepsilon 
\end{array}
\right]
$$
The spectral radius of the matrix :
$$
M(\varepsilon,m,T,\varphi) =\emat^{T(1-\varphi)M_{\varepsilon,m}^{--}}  \emat^{\varphi T M_{\varepsilon,m}^{-+}} \emat^{T(1-\varphi)M_{\varepsilon,m}^{++}}  \emat^{\varphi T M_{\varepsilon,m}^{+-}} 
$$
decides of the stability of the switched system associated to these four systems, $T$ and $\varphi$.  Once again, we ask to Maple to compute the eigenvalues of $M(\varepsilon,m,T,\varphi)$, we select the largest one $\lambda_1(\varepsilon,m,T,\varphi)$  and look for le mapping $(m,T) \mapsto 1/T \ln( \lambda_1(\varepsilon,m,T,\varphi))$ for $\eps = 0.1$ and various values. When the shift $\varphi T$ is not equal to $T$ our intuition is that the inflation effect will be proportional  to the shift and will be maximum when $\varphi = 1$. This is confirmed by Figure \ref{VPshift}.
  \begin{figure}
  \begin{center}
 \includegraphics[width=1\textwidth]{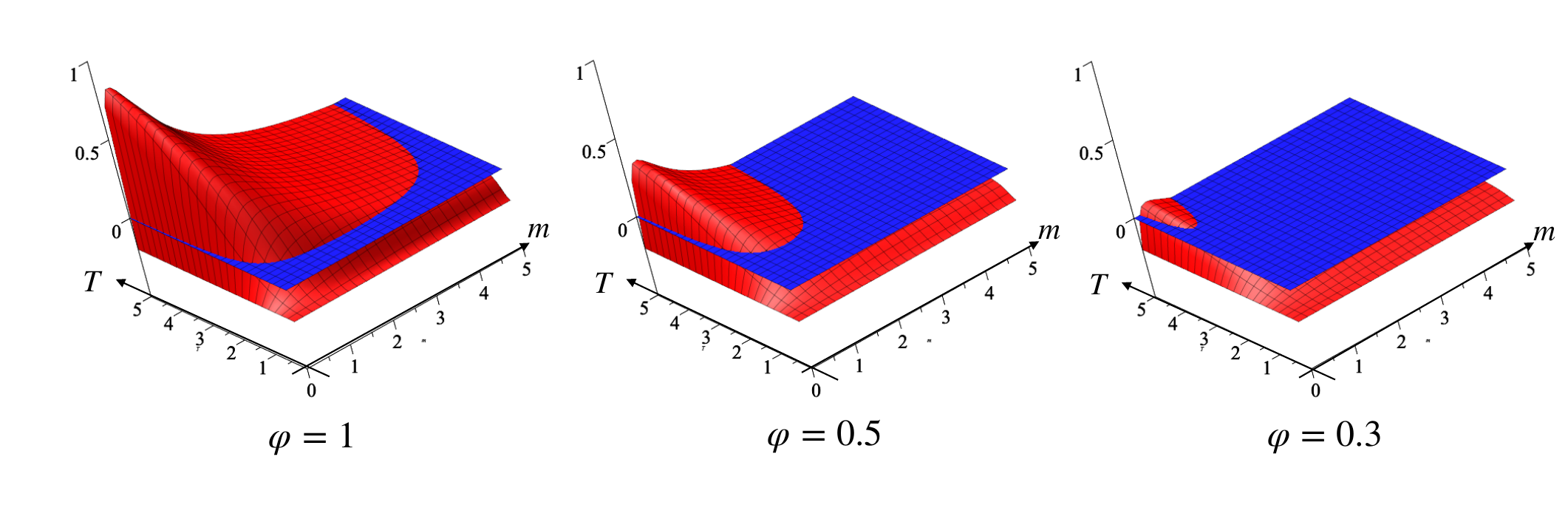}
 \caption{Graphs of $(m,T), \mapsto \frac{1}{T}\ln(\lambda_1(0.1,m,T,\varphi))$ for three values of $\varphi$ } \label{VPshift}
 \end{center}
 \end{figure}

\subsection{Migration between different patches}
{As discussed at the end of Section \ref{modelpm}, the $(\pm 1)$-model given by \eqref{systeme1} encompasses the more general case of two identical patches that are in phase opposition. Let us show now that the patches do not need to be identical
and that our approach applies in the more general case of model \eqref{BLSS1}, where the functions $r_1(t)$ and $r_2(t)$ are given by 
\begin{equation}\label{r1r2}
r_1(t)=\left\{
\begin{array}{rcl}
r_1&\mbox{if}&t\in[0,T]\\
-d_1&\mbox{if}&t\in[T,2T]
\end{array}
\right.
\quad
r_2(t)=\left\{
\begin{array}{rcl}
-d_2&\mbox{if}&t\in[0,T]\\
r_2&\mbox{if}&t\in[T,2T]
\end{array}
\right.
\end{equation}
where $r_1$, $r_2$, $d_1$ and $d_2$ are real parameters. The system $\Sigma(\eps,m,T)$, defined by \eqref{systeme1} corresponds to the case where $r_1=r_2=1-\eps$ and $d_1=d_2=1+\eps$. On the other hand, the system \eqref{+r-d} corresponds to the case where $r_1=r_2=r$ and $d_1=d_2=d$.

Using notations similar to those we used in Subsection \ref{variables initiales}, let us define~:
 $$M^{1}_{r_1,d_2,m} = 
\left[
\begin{array}{cc}
r_1-m &+m\\
+m & -d_2-m 
\end{array}
\right],
\qquad
M^{2}_{r_2,d_1,m} = 
\left[
\begin{array}{cc}
-d_1-m &+m\\
+m &r_2-m  
\end{array}
\right]
$$
The spectral radius of the matrix :
$$
M(r_1,d_1,r_2,d_2,m,T) =\emat^{TM_{r_2,d_1,m}^2}  \emat^{T M_{r_1,d_2,m}^1} 
$$
decides of the stability of the switched system. 
Once again, we ask Maple to compute the eigenvalues of $M(r_1,d_1,r_2,d_2,m,T)$, we select the largest, denoted $\lambda_1(r_1,d_1,r_2,d_2,m,T)$ and draw the graph of the function 
$(m,T) \mapsto \frac{1}{2T} \ln( \lambda_1(r_1,d_1,r_2,d_2,m,T)$ for various values of the parameters, see Fig. \ref{lambdasitesdifferents}.

  \begin{figure}
  \begin{center}
 \includegraphics[width=1\textwidth]{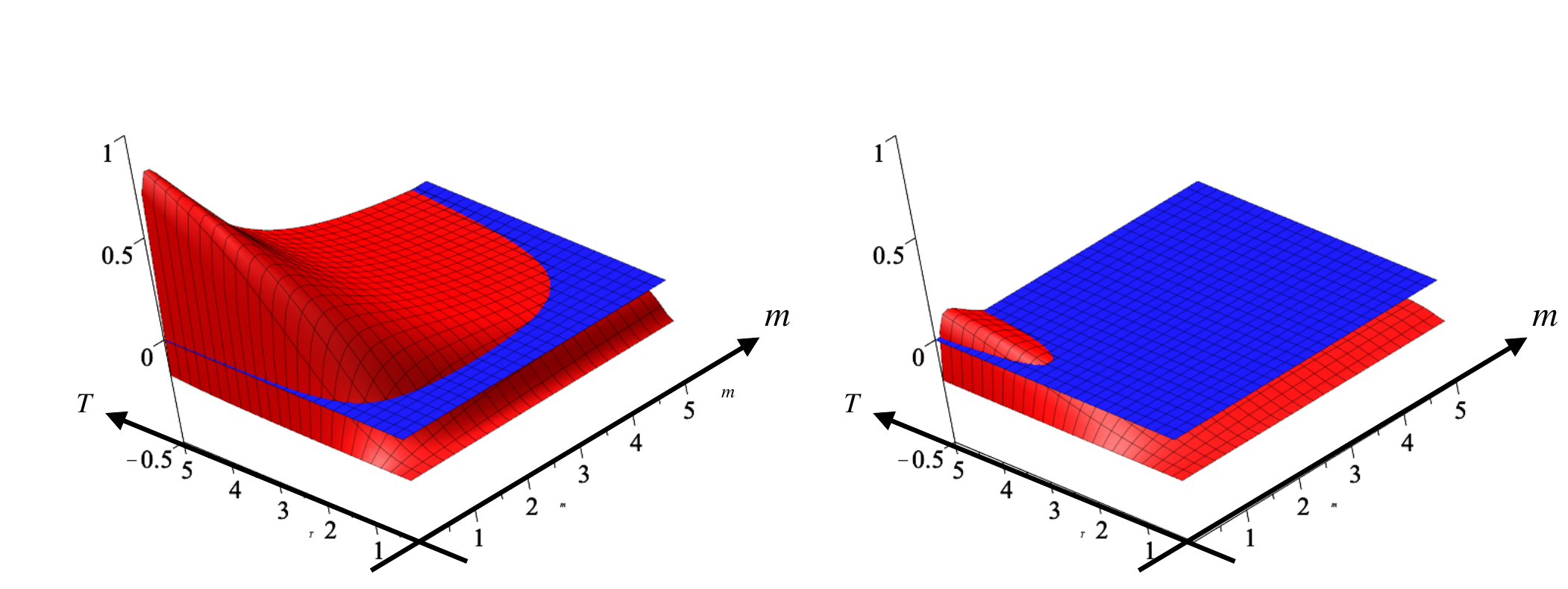}
 \caption{Graphs of $(m,T), \mapsto \frac{1}{T}\ln(\lambda_1(r_1,d_1,r_2,d_2,m,T))$ for two values of $(r_1,d_1,r_2,d_2)$ } \label{lambdasitesdifferents}
 \end{center}
 \end{figure}
\paragraph*{Comments on figure \ref{lambdasitesdifferents}.}

{We look to $\frac{1}{2T}\ln(\lambda_1(r_1,d_1,r_2,d_2,m,T))$ in two different cases.  On the left we consider the case :
\beq
\begin{array}{|l|l|}
\hline 
r_1=0.9&d_1 = 1.1 \\
\hline
r_2=0.9&d_2 =1.1\\
\hline
\end{array}
\feq
which is the case of the $(±1)$model for $\eps = 0.1$ which we already considered.
We compare this case to the case :
\beq
\begin{array}{|l|l|}
\hline 
r_1=0.9&d_1 = 1.1 \\
\hline
r_2=-0.1&d_2 =0.1\\
\hline
\end{array}
\feq
In this case the patch  1 is unchanged and the patch  2 represent some place without seasonality.
In this case inflation is smaller but   still observable.
\begin{figure}
  \begin{center}
 \includegraphics[width=0.6\textwidth]{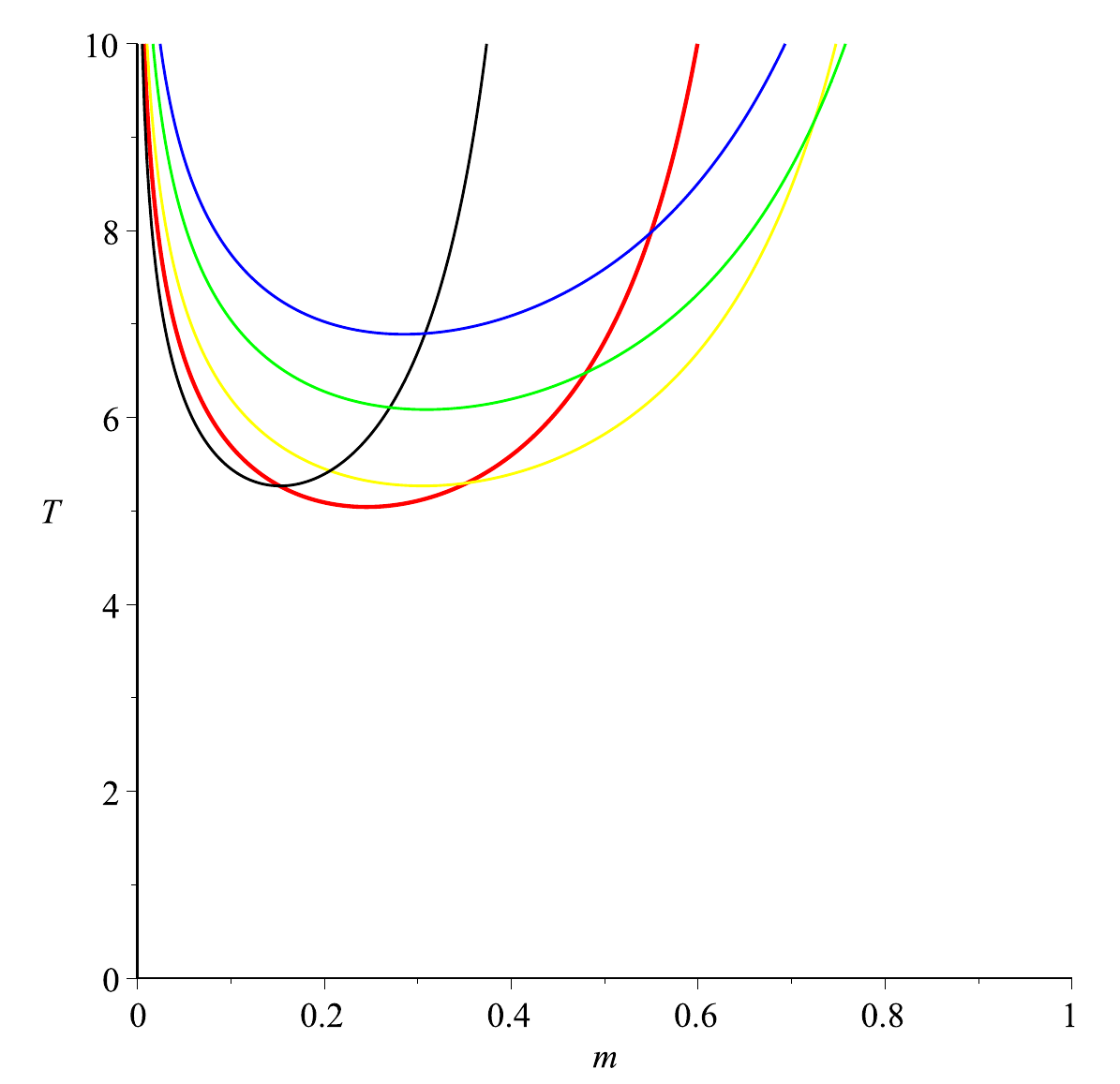}
 \caption{The set  $\{(m,T):\Delta(0.5,m,T,\gamma)=0\}$ for the values of $\gamma$: $\gamma=2$ (Black), 
 $\gamma=1$ (Red), $\gamma=0.5$ (Yellow), $\gamma=0.2$ (Green), $\gamma=0.1$ (Blue).} \label{gammadifferents}
 \end{center}
 \end{figure}

\subsection{The case of non symmetric dispersal}\label{nonsym}
The symmetric rate of dispersal between the two patches is a very special (and unlikely) situation.
A non symmetric dispersal like in the model :
\begin{equation}\label{r_igamma}
\left\{
\begin{array}{lcrcl}
\displaystyle \frac{dx_1}{dt}&=&\displaystyle r_1(t)x_1&+&m(\gamma x_2-x_1)\\[2mm]
\displaystyle\frac{dx_2}{dt}&=&\displaystyle r_2(t)x_2&+&m(x_1-\gamma x_2)
\end{array}
\right.
\end{equation}
with $\gamma>0$, is certainly more realistic.
Using the change of variables  
$$\begin{array}{l}
U = \ln\left(x_1^\gamma x_2\right)^\frac{1}{\gamma+1}=\frac{\gamma\ln x_1+\ln x_2}{\gamma +1}
,\qquad
V = \ln\sqrt{\frac{x_1}{ x_2}}
=\frac{\ln x_1-\ln x_2}{2},
\end{array}
$$
 one obtains :
\begin{equation}\label{r_iBetaUV}
\left\{
\begin{array}{lcl}
\frac{dU}{dt}&=&\frac{\gamma r_1(t)+r_2(t)}{\gamma+1}+\frac{2\gamma m}{\gamma+1} \left(\cosh\left(2V-\ln\gamma\right)-1\right)\\[2mm]
\frac{dV}{dt }&=&\frac{r_1(t)-r_2(t)}{2}-m\left(\sqrt{\gamma}\sinh\left(2V-\frac{\ln\gamma}{2}\right)+\frac{\gamma-1}{2}\right)
\end{array}
\right..
\end{equation}
This system reduces to  (\ref{systvarUV}) in the symmetric case $\gamma=1$,
and its study will follow the same lines than the study of \eqref{systvarUV}. In particular, since 
$\cosh(\alpha)\geq 1$, we have 
$$
U(t)\geq U(0)+\int_0^t
\frac{\gamma r_1(s)+r_2(s)}{\gamma+1}ds
$$
Therefore, we have
$$\liminf\frac{U(t)}{t}\geq \frac{\gamma \bar r_1+\bar r_2}{\gamma+1}$$
Note that, as in Appendix \ref{tycho}, we can use singular perturbation theory 
\cite{TYK52,LOB98} to show that 
\[
\lim_{m\to \infty}\Delta(r_1(\cdot),r_2(\cdot),m,T)=
\frac{\gamma \bar r_1+\bar r_2}{\gamma+1}
\]
Therefore, we have
$$\inf\Delta(r_1(\cdot),r_2(\cdot),m,T)=\frac{\gamma \bar r_1+\bar r_2}{\gamma+1}.$$

On the other hand, for the $(\pm 1)$ model associated to the asymmetric dispersal \eqref{r_igamma}, we simply consider the matrices~:
 $$M^{+-}_{\varepsilon,m,\gamma} = 
\left[
\begin{array}{cc}
1-\varepsilon-m &\gamma m\\
m & -1-\varepsilon-\gamma m 
\end{array}
\right],
\qquad
M^{-+}_{\varepsilon,m,\gamma} = 
\left[
\begin{array}{cc}
-1-\varepsilon-m &\gamma m\\
m &1-\varepsilon -\gamma m  
\end{array}
\right]
$$
The spectral radius of the matrix :
$$
M(\varepsilon,m,T,\gamma) =\emat^{TM^{-+}_{\varepsilon,m,\gamma}} \emat^{T M^{+-}_{\varepsilon,m,\gamma}} 
$$
decides of the stability of the switched system. 
Once again, we ask Maple to compute the eigenvalues of $M(\varepsilon,m,T,\gamma)$, we select the largest, denoted $\lambda_1(M(\varepsilon,m,T,\gamma))$. To have a better understanding of the role of $\gamma$, we depict in Fig. \ref{gammadifferents} the zero level-set 
$$
\{(m,T):\Delta(\varepsilon,m,T,\gamma)=0\}
$$
of the function 
$\Delta(\varepsilon,m,T,\gamma)=
\frac{1}{2T} \ln( \lambda_1(M(\varepsilon,m,T,\gamma))$, 
for various values of $\gamma$.

 \subsection{A  density dependent deterministic model}\label{ddmodel}
In \cite{ARD15,ARD18} a complete description of the asymptotic behavior of the model :
 \beq\label{ARD18}
\begin{array}{rcl} 
  \displaystyle \frac{dx_1}{dt} &=& \displaystyle  r_1x_1\left( 1- \frac{x_1}{K_1}\right)  + \beta\left( \frac{ x_2}{\gamma_2}- \frac{x_1}{\gamma_1}\right) \\[8pt]
    \displaystyle \frac{dx_2}{dt} &=& \displaystyle  r_2x_1\left( 1- \frac{x_2}{K_2}\right)  + \beta\left( \frac{ x_1}{\gamma_1}- \frac{x_2}{\gamma_2}\right) 
 \end{array}
\feq
is given in the space of the six independent  parameters $\{ r_i, K_i, (i = 1,2),  \beta/\gamma_2,  \gamma_1/\gamma_2 \}$, the focus being on the comparison between   the total equilibrium population with the sum $K_1+K_2$ of the two carrying capacities.  Here we complement this study by considering the question of persistence when $r_1$ and  $r_2$ vary in time for specific values of the parameters. Namely, we consider  the system
\beq \label{logistic1}
D(\varepsilon,\alpha,m,T) \quad \quad
\left\{
\begin{array}{lcl}
\displaystyle  \frac{dx_1}{dt} &=& (+u(t)-\varepsilon)x_1-\alpha x_1^2+m(x_2-x_1)\\[12pt]
\displaystyle  \frac{dx_2}{dt} &=&(-u(t)-\varepsilon)x_2-\alpha x_2^2+m(x_1-x_2) 
 \end{array} 
 \right.
 \feq
where $0 \leq \varepsilon \leq 1$, $\alpha \geq 0$, $  m \geq 0$, $ T\geq 0$ and the function   $t \mapsto u(t)$ is periodic of period $2T$, with
$$ t \in [0,\,T[ \Rightarrow u(t) = 1\quad \quad  t \in [T,\,2T[ \Rightarrow u(t) = -1$$

We are interested in the {\em persistence} of (\ref{logistic1}). 
Recall that the system $D(\varepsilon,\alpha,m,T)$ is {\em uniformly persistent} (see for instance \cite{BUT86}) if there exist  strictly positive constants $a < b$  such that  every solutions $(x_1(t),x_2(t))$ of  $D(\varepsilon,\alpha,m,T)$ is asymptotically bounded from below by $a$ and from above by $b$ (i.e.  $a \leq x_i(t) \leq b $ for $t$ sufficiently large).

When $\alpha = 0$ the system $D(\varepsilon,0,m,T)$ is just the $(±1)$model $\Sigma(\varepsilon,m,T)$. When $\alpha$ is not 0,  but $m = 0$, on each patch the dynamic is :
\beq \label{up1}
\frac{dx_i}{dt} = (u(t)-\varepsilon)x_i-\alpha x_i^2\quad \quad i = 1,2
\feq
with $u(t) = ±1$. 
In both cases one has a logistic equation with a globally stable equilibrium equal to $\frac{1- \varepsilon }{\alpha}$ or $0$. One sees easily that in the space $(\mathrm{R}^+)^2$ the square S = $[0,\frac{1-\eps}{\alpha}]\times [0,\frac{1-\eps}{\alpha}]$ is an invariant global attractor~; this implies that every trajectories of (\ref{logistic1}) are bounded from above.

Regarding boundedness from below we can say intuitively that  the system 
$D(\varepsilon,\alpha,m,T)$ behaves around the origin like its linear approximation, namely the system $\Sigma(\varepsilon,m,T)$ and thus is persistant if and only if $\Sigma(\varepsilon,m,T)$ is unstable.
Actually the following proposition can be proved (see appendix \ref{densitedependant}) 
  \begin{figure}
  \begin{center}
 \includegraphics[width=1\textwidth]{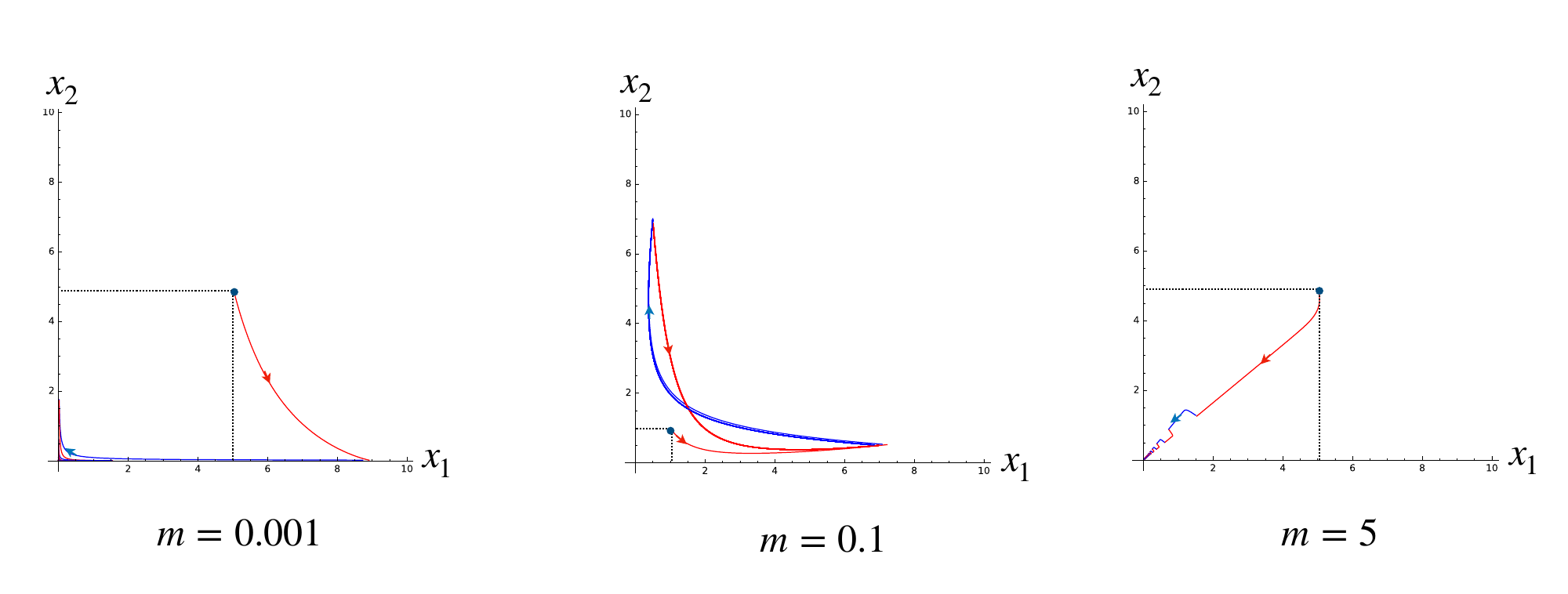}
 \caption{Simulations of $D(0.1,0.1,m,5)$  showing persistence for intermediate values of $m$.} \label{persis}
 \end{center}
 \end{figure}
{\proposition \label{persistance} If the parameters $(\eps,m,T)$ are such that the system $D(\varepsilon,0,m,T) =\Sigma(\varepsilon,m,T)$ is~:
\bit
\item stable, then the solutions of $D(\varepsilon,\alpha,m,T)$ tend to $0$ (extinction),
\item unstable, then $D(\varepsilon,\alpha,m,T)$ is persistent.
\fit}
Thus we have, for $T$ large enough,  the sequence : small $m$ : extinction $-$ intermediate $m$ : persistence $-$ large $m$ : extinction/ This is illustrated by the simulations of Figure  \ref{persis}.
\subsection{A density dependent stochastic model}
\label{stochasticdd}
In this short section, we show that the Proposition \ref{persistance} is still true under a random signal $\bm{u}$. More precisely, we consider the system 
\beq \label{logistic1r}
\bm{D}(\varepsilon,\alpha,m,T) \quad \quad
\left\{
\begin{array}{lcl}
\displaystyle  \frac{dx_1}{dt} &=& (+\bm{u}(t)-\varepsilon)x_1-\alpha x_1^2+m(x_2-x_1)\\[12pt]
\displaystyle  \frac{dx_2}{dt} &=&(-\bm{u}(t)-\varepsilon)x_2-\alpha x_2^2+m(x_1-x_2) 
 \end{array} 
 \right.
 \feq
where $\bm{u}$ switches from $1$ to $-1$ and conversely at random exponential time, as described in Section \ref{sec:PDMP}. Like in the periodic case described above, when $\alpha = 0$, $\bm{D}(\varepsilon,\alpha,m,T)$ is just the stochastic $(±1)$model $\bm{\Sigma}(\varepsilon,m,T)$. Using a terminology borrowed to Schreiber and Chesson, we say that the system  $D(\varepsilon,\alpha,m,T)$ is \emph{stochastically persistent} if for all $\eta > 0$, there exists a compact set $K_{\eta} \subset \mathbb{R}_{++}^2$ such that, almost surely, 
\[
\liminf_{ t \to \infty} \frac{1}{t} \int_0^t \1_{(x_1(s), x_2(s)) \in K_{\eta}} \, ds \geq 1 - \eta.
\] 
We now give the stochastic counterpart of Proposition \ref{persistance}: 
{\proposition \label{persistancesto} We have the following dichotomy:
\bit
\item If $\bm{\Delta}(\varepsilon,m,T) \leq 0$, then system $D(\varepsilon,\alpha,m,T)$  goes to extinction;
\item If $\bm{\Delta}(\varepsilon,m,T) > 0$, then system $D(\varepsilon,\alpha,m,T)$ is stochastically persistent, and the process $(x_1,x_2,u)$ admits a unique stationary distribution $\nu$ such that $\nu( \mathbb{R}_{++}^2 \times \{\pm 1\} ) = 1$.
\fit}

  \subsection{An S.I.R. type  epidemic model} \label{holt}
 In  \cite{HOLT20} {\em Nicholas Kortessisa, Margaret W. Simon, Michael Barfield, Gregory Glass, Burton H. Singer and Robert D. Holt } consider the classical S.I.R. model for a population living in two patches connected by migration.
 The model is the following system :
 \beq\label{Holt1}
\begin{array}{rcl} 
  \displaystyle \frac{dS_1}{dt} &=& -  \beta(t) S_1I_1 +m(S_2-S_1)  \\[8pt]
 \displaystyle \frac{dI_1}{dt} &=&+ \beta(t) S_1I_1  -( \gamma(t) +\mu)I_1   + m(I_2-I_I)        \\[8pt]
  \displaystyle \frac{dS_2}{dt} &=& -  \beta(t-\varphi) S_2I_2 +m(S_1-S_2)  \\[8pt]
 \displaystyle \frac{dI_2}{dt} &=&+ \beta(t-\varphi) S_2I_2  - (\gamma(t-\varphi)+\mu)I_2   +m(I_1-I_2)       
 \end{array}
\feq
where $S_i(t)$ represents the number of ''susceptible to be infected'' at time $t$ on each patch,   $I_i(t)$ represents the number of ''infected'' on each patch. The parameters $\beta(.) $ and $\gamma(.)$  are piecewise constant functions of period $2T$ varying according to the presence or absence of social distancing measures~; we examine the messages of this paper in light of our previous study of inflation phenomenon\footnote{ The authors of \cite{HOLT20} publish the same message in P.N.A.S. \cite{HOLTPNAS20}  but using, in our opinion, a less realistic $\beta(.) $ and  $\gamma(.)$ like continuous sinusoidal functions.  We prefer to refer to the initial paper but our discussion would be the same with the P.N.A.S. paper. }.

The first remark of the authors of \cite{HOLT20} is to consider that we are essentially interested by the beginning of the contamination when, as a first approximation, we can consider that $S(t)$ is almost equal to the initial total population $N$. Then the approximate model is :
 \beq\label{Holt2}
\begin{array}{rcl} 
   \displaystyle \frac{dI_1}{dt} &=&\big( \beta(t) N  -( \gamma(t) +\mu)\big)I_1   + m(I_2-I_I)        \\[8pt]
 \displaystyle \frac{dI_2}{dt} &=&\big( \beta(t-\varphi) N  - (\gamma(t-\varphi)+\mu)\big)I_2   +m(I_1-I_2)       
 \end{array}
\feq
They denote respectively by the subscripts $n$ and $s$ the values of parameters in ''normal'' periods and periods when the ''social distancing'' is in effect. They adopt, according to the current literature, the following realistic values .

\begin{center}
\begin{tabular}{|c||c||c|c|} 
\hline
$ \beta_nN$ = 0.1988&$\gamma_n$ = 0.098&$\mu_n$ = 0.002\\[8pt]
\hline
$ \beta_sN $ = 0.0288&$\gamma_s$ = 0.128&$\mu_s$ = 0.002\\[8pt]
\hline
 \end{tabular}
 \end{center}

  \begin{center}
\begin{tabular}{|c|} 
\hline
$ \beta_nN -(\gamma_n+\mu_n) $ = 0.0988\\[8pt]
\hline
$ \beta_sN -(\gamma_s+\mu_s) $ =  - 0.1012\\[8pt]
\hline
 \end{tabular}
 \end{center}
  \begin{figure}
  \begin{center}
 \includegraphics[width=1\textwidth]{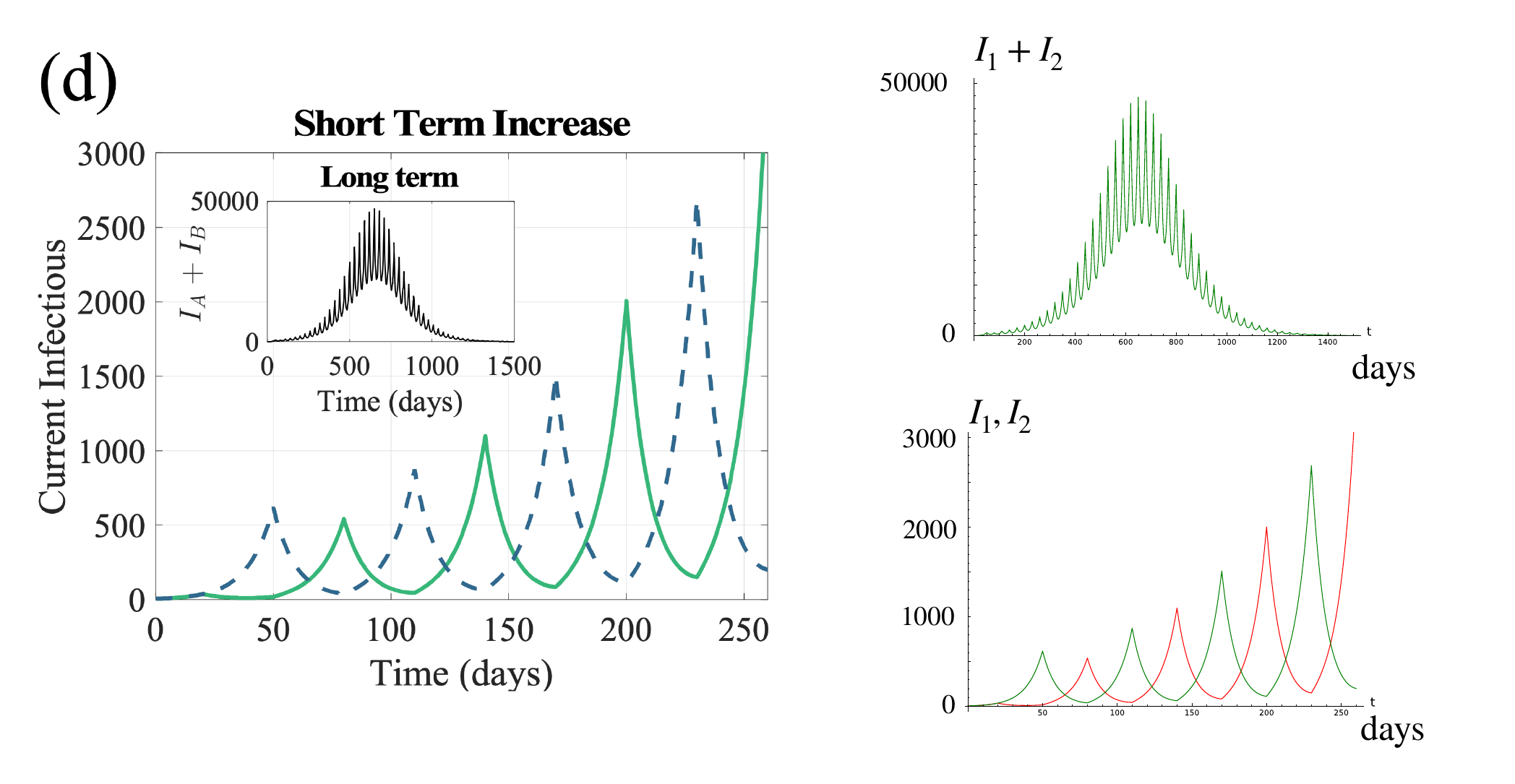}
 \caption{Simulation from \cite{HOLT20} (left) ; our simulation (right)} \label{Holtmoi}
 \end{center}
 \end{figure}
and they  discuss the case $T=30$.
 	We have done a simulation with these parameters and $m = 0.005$ as they did. We obtained the same picture than \cite{HOLT20} (see figure \ref{Holtmoi}) which confirms that we are actually running the same model but our objective is not to reproduce \cite{HOLT20} results but to complete them.  For this purpose we consider the effect of migration, in the case of a small phase shift in the application of social distancing. We assume that $\varphi = 4$ days.
 
 In the absence of migration the linear model is :
  \beq\label{Holt3}
\begin{array}{l} 
   \displaystyle \frac{dI_ i}{dt} =0.0988\, I_i \quad \mathrm{(normal)}\quad \quad \frac{dI_ i}{dt} =-0.1012 \, I_i \quad \mathrm{(social \;distancing)}
 \end{array}
\feq
If we multiply the dynamic of $I_i$ by the factor 10 (which means a change unit for the time) we have :
 \beq\label{Holt4}
\begin{array}{l} 
   \displaystyle \frac{dI_ i}{dt} =0.988\, I_i \quad \mathrm{(normal)}\quad \quad \frac{dI_ i}{dt} =-1.012 \, I_i \quad \mathrm{(social \;distancing)}
 \end{array}
\feq
which we read as the  ''$(±1)$model'':
  \beq\label{Holt3basic}
\begin{array}{l} 
   \displaystyle \frac{dx_ i}{dt} =(1-\eps)x_i \quad \mathrm{(normal)}\quad \quad \frac{dx_ i}{dt} = - (1+\eps)   x_i \quad \mathrm{(social \;distancing)}
 \end{array}
\feq
with $\eps = 0.012$. To $T = 30$ and a  phase shift of $4$ days in the model (\ref{Holt1}) correspond $T = 3$ a shift of $0.133$ in (\ref{Holt3basic}). 
  \begin{figure}
  \begin{center}
 \includegraphics[width=1\textwidth]{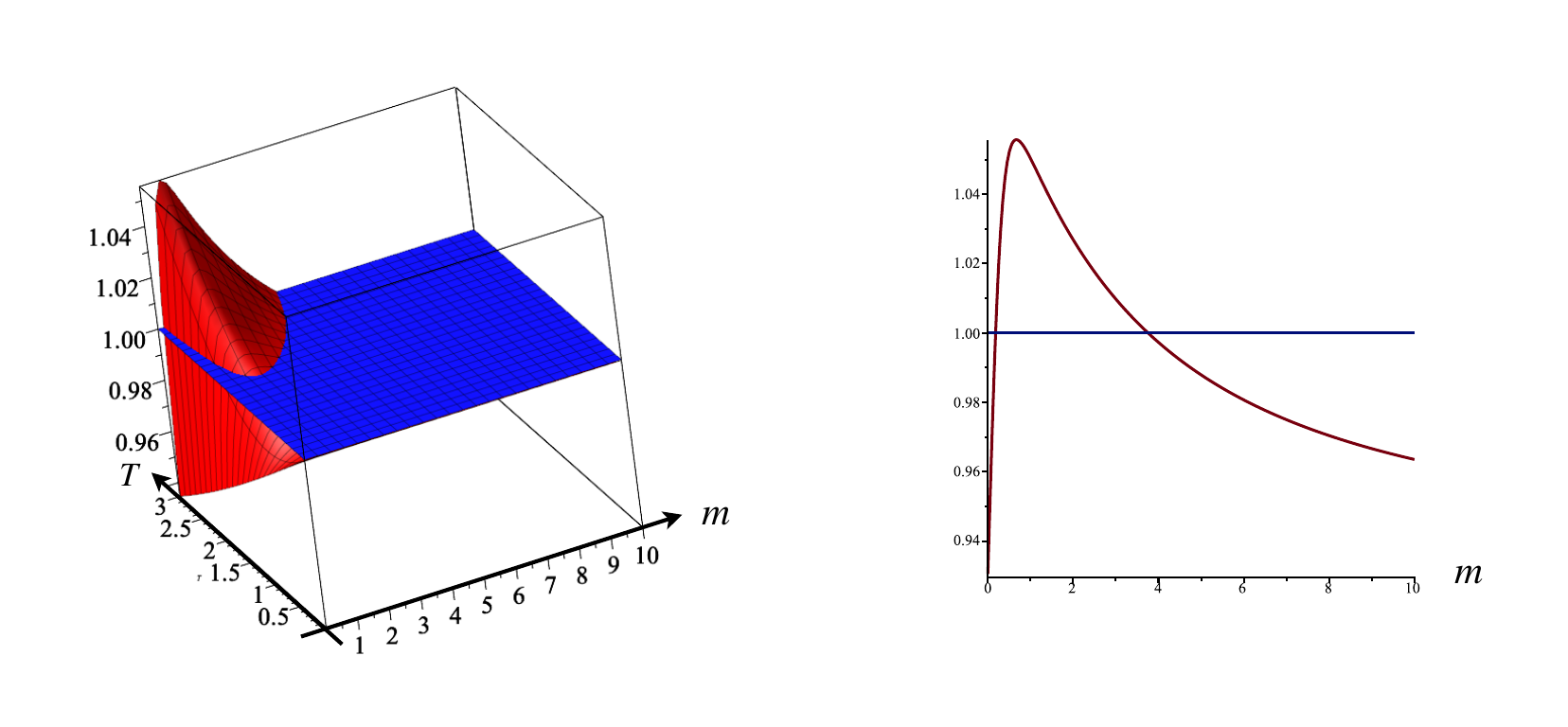}
 \caption{Graphs of $(m,T), \mapsto \lambda_1(0.012,m,T)$ (left) and  $ m \mapsto \lambda_1(0.012,m,3,0.133)$} \label{Hshift}
 \end{center}
 \end{figure}
For these values one sees on the graphs of $\lambda_1$ that there is no longer instability for $m> 2$ which means $m > 0.2$ in the original system (\ref{Holt1}). This  must be reflected on the epidemic. If one looks to the cumulative number of cases for a duration of 1500 days the simulation of the model gives the figure \ref{Hshiftbis}.
 \begin{figure}
  \begin{center}
 \includegraphics[width=0.8\textwidth]{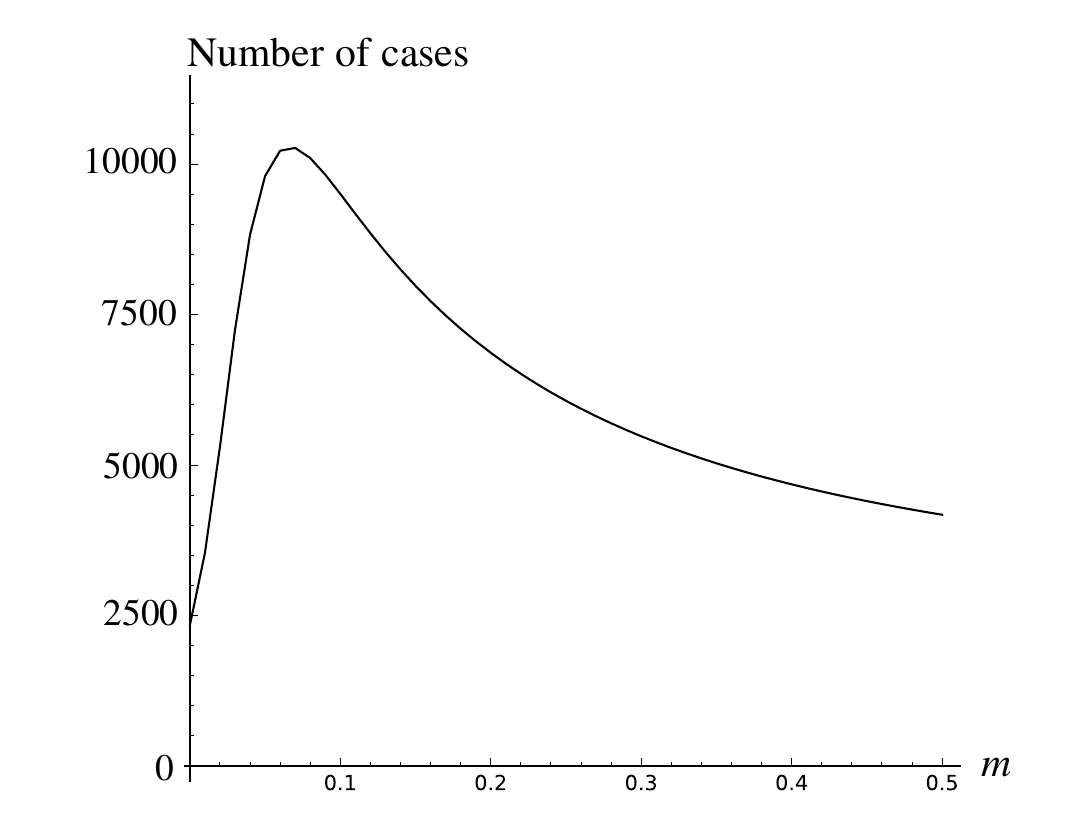}
 \caption{Cumulative number of cases up to 1500 days as a function of migration.} \label{Hshiftbis}
 \end{center}
 \end{figure}
 We can see that no migration at all is the best, when migration grows from 0 to approximately $0.1$ the number of cases is multiplied by 4 but and after that decreases. {\em Very low migration can increase dramatically the number of cases, while, if migration is unavoidable, comparatively large one has better effect}.
 
In their paper  published in the P.N.A.S. the authors \cite{HOLTPNAS20} say :
 "These findings highlight a need for integrated, holistic policy: Intensify mitigation locally, coordinate tactics among locations, and reduce movement."
 
 In the light of our work, we see that the latter recommendation, is not necessarily correct, 
depending on where you are located with respect to the maximum of $\Lambda(\varepsilon, m,T)$. This does not invalidate but reinforces the conclusion of their paper with which we fully agree. 
 "It is increasingly recognized that monitoring and controlling movement is essential for effective pandemic control. The impact of such actions is, however, contextual, because their dynamical effects are intertwined with the magnitude of asynchrony in local transmission across space. More-realistic, spatially structured epidemiological models including movement and asynchronous transmission $-$ at scales from local to international $-$ are essential to control this and future pandemics in the coupled metapopulations of humans and their pathogens."
\section{Discussion}\label{discussion}

 We first studied the simplest model likely to present the inflation phenomenon for continuous time models. 
   For that purpose we considered only two patches, a symmetric migration of rate $m$ between the two patches and  piecewise constant environments. On the patch 1, for a duration $T$ the environment is favorable, the growth rate is $1-\eps$ and it is followed by a period also of duration $T$ where the environment is unfavorable and the  rate of decay is $-1-\eps$ ; since  $\eps $ is strictly  positive , over the period $2T$ the environment is globally unfavorable, so patch 1 considered alone, is a sink. On patch 2 we consider an identical environment, so this site is also a sink, but we suppose that it is out of phase of half a period with patch 1; thus when the environment is favorable on site 1 it is unfavorable on patch 2 and conversely. On this minimal model that we noted $\Sigma(\eps,m,T)$ we showed that the Liapunov exponent $\Delta(\eps,m,T)$ which characterizes the growth of $\Sigma(\eps,m,T)$ has the following properties:
 \ben
 \item For all $T$, for $m$ small ($m \rightarrow 0$) and $m$ large ($m \rightarrow +\infty$) we have $ \Delta(\eps,m,T) < 0$ (no inflation).(see prop. \ref{propgenerale})
 \item For any $m <  (1-\eps^2)/2\eps$ there is a threshold $T^*(\eps,m)$  for the half period $T$ below which there is no inflation ($ \Delta(\eps,m,T) < 0$) and above which there is inflation ($ \Delta(\eps,m,T) > 0$).
 \item We have given an explicit formula (see prop. \ref{explicitdelta} ) for $ \Delta(\eps,m,T) $ from which we deduce (see prop. \ref{approxseuil}) that the 
 the threshold value at which $m \mapsto \Delta(\eps,m,T)$ becomes positive is an exponentially small value 
  $ \sim \emat^{-(1-\eps)T}$ for large values of $T$.
 \fen
 The graph of $(m,T) \mapsto \Delta(\eps,m,T)$ (see figure \ref{graphdelta}) summarizes the situation.

 Note that the two first points follow from the general result of Katriel  (\cite{KAT22}) but the latter is only proved for growth rates that depend continuously on time and thus cannot (formally) apply to our situation  we don't know if the methods of \cite{KAT22} apply to switched systems. By the way our results complement Katriel's results.
 
 Our results can be obtained (thanks to formal software like Maple) directly from the explicit (but rather obscure) formula for the dominant eigenvalue of the matrix which defines the Poincaré application of the periodic system. But we have preferred to establish them from the examination of the phase portrait of the transformed system in the variables 
   $(U,V)$, $U$ being the geometric mean of the abundances $(x_1,x_2)$ and $|V|$ their geometric standard deviation. 

These variables allows to easily extend our understanding of the inflation phenomenon on  the $(±1)$-model to more general situations. In this spirit 
we have also shown  that the properties of the $(±1)$ model are still valid for some kinds of  continuous time stochastic models which is important 
if we want to be more realistic.  First we have reconsidered the $(±1)$-model assuming that the successive sequences of constant environment are not fixed an equal to $T$ like in the $2T$
 periodic environment but are succession of  independent  random duration $S_n$  following the same law $\mu$.  When the law is an exponential law of parameter $\lambda$ this defines a so 
called Piecewise Deterministic Markov Process (see \cite{DAV84}), which is interesting since we are able to make explicit calculations in that case. We have shown that the
 same inflation phenomenon occurs with the expectation of $S_n$ playing the role of the period in the deterministic periodic case. More precisely, we have shown that there exists a
 unique (deterministic) Lyapunov exponent $\bm{\Delta}(\varepsilon, m, \mu)$ whose sign characterise the behaviour of the system. In addition, we prove that for all $m$, provided the expectation
 of $S_n$ is large enough, $\bm{\Delta}(\varepsilon, m, \mu)$ is positive and there is inflation.} We have also considered a different case of stochasticity.  Now the duration of the constant growth rate $(±1)$
 is of a fixed duration $T$ but the choice of $+1$ or $-1$ is random, these choices being correlated or not. For this model we also characterize the presence of inflation, by the sign of a
 deterministic Lyapunov exponent $\Delta(p_1, p_2, m, T, \varepsilon)$. Moreover, we prove that there exists a threshold $\chi(p_1, p_2,\varepsilon)$ such that
\ben
\item If $\chi(p_1, p_2,\varepsilon)<0$, then for all $(m,T)$, $\Delta(p_1, p_2, m, T, \varepsilon) < 0$ and there is no inflation,
\item if $\chi(p_1, p_2,\varepsilon)>0$, then for all $m <  (1-\eps^2)/2\eps$ there is a threshold $T^*(\eps,m)$, there is inflation whenever $T \geq T^*(\eps, m)$.
\fen

  All the above results apply to the $(±1)$ model which, of course, is totally unrealistic. This is why we have shown how our results extend to, or illuminate, the more general situations listed in the index of section \ref{morecomplex}, which ends with the application to the epidemiological model of \cite{HOLTPNAS20} that motivated this study.

This mathematical success of the $(U,V)$ variables also suggests that the geometric mean and deviation of abundances at each site are a better indicator of the metapopulation status than are the arithmetic mean and standard deviation.

    We have shown that the growth rate of the two patches is always higher than 
    the mean of the growth rates within each patch, see Remark \ref{rem:avantageU}.
    Therefore, in the $(\pm)$ model, where the growth rates within each patch 
    are equal (i.e. $\bar{r}_1=\bar{r}_2$), the growth rate 
    $\Delta(\varepsilon,m,T)$ is always
    higher than its limit $\Delta(\varepsilon,0,T)$, for $m=0$, see Fig. \ref{graphdelta} and 
    Proposition \ref{prop:asymptoticDelta}. This property is also satisfied in the examples depicted in 
    Figs. \ref{VPshift} and \ref{lambdasitesdifferents}. We can therefore deduce that migration always favors growth.  This positive effect may be sufficient to change a negative growth rate when no migration to 
     a positive growth rate (i.e., inflation), but may also not (i.e., 
    not inflation). But in either case, the growth rate is inflated.However, the property $\Delta(\varepsilon,m,T)> \Delta(\varepsilon,0,T)$, if $m>0$, is specific to the case where $\bar{r}_1=\bar{r}_2$. If these growth  rates are different then the limit of $\Delta(\varepsilon,m,T )$ when $T$
tends to 0 is a strictly decreasing function with respect to $m$, from 
$\max(\bar{r}_1,\bar{r}_2)$ to $(\bar{r}_1+\bar{r}_2)/2$, see \cite[Lemma 9]{KAT22}. Consequently, when $m$ and $T$ are sufficiently small, 
$\Delta(\varepsilon,m,T)<\Delta(\varepsilon,0,T)$ will result, and the 
migration will not favor growth. For a more in-depth analysis of this 
behavior, the reader is referred to \cite{BLSS23}. 
 
  If we try to understand intuitively the mechanisms that cause inflation {\em in the case of continuous time models } we see the following. 
  Let us say that the environment   is {\em positive} (respectively {\em negative}) at time $t$ on some patch if the abundance of the population is {\em increasing} (respectively {\em decreasing}).  The fluctuations of the environment on the two patches can be thought as a succession of regimes  $(++),(+-),(- +),(- -)$ and the key ingredients for inflation are :
  \bittiret 
 \item A sufficient (in mean) duration between two changes of regime (i.e. a period large enough in the deterministic periodic case), 
 \item a proportion of time spent in opposite regimes $(+-)$ or $(-+)$  as large as possible,
 \item a migration neither too weak nor too strong.
 \fit
 which we understand as follows. In the absence of migration, a  $(+-)$ or $(-+)$ regime, will create rapidly a dissymmetry between the two sites, since the abundance on the source is increasing and the abundance on the sink is decreasing as long as one remains in this ; if this duration is long enough the ratio of the two abundances, tends to $ \infty$ or$0$  ; if we now introduce a migration it induces a transfert from the positive environment to the negative one ; if it is too weak it does not allow the patch with negative environment to benefit from the growth on the other one and, conversely, a migration that is too strong slows down the growth the patch with positive environment.

 This description of the mechanisms at the origin of inflation contrasts with those put forward in the case of {\em discrete-time random models} 
 where the emphasis is generally placed on the need for a temporal correlation of the fitnesses at the two patches (\cite{HOLT97, ROY05, SCH10}, and, in some cases only, (\cite{SCH10} for instance ) the need for a migration that is neither too large nor too small. We believe that the origin of this difference lies in the following remark.

 What makes the intuition of the inflation phenomenon intricate, in both discrete and continuous time models, is to consider simultaneously growth and/or decay on each of the patches and a migration in both directions, but basically things are rather simple if we suppose that patch  2 is a patch where there is neither growth nor decay, that there is only migration from 1 to 2 and that we reverse the roles of patches 1 and 2 after a time $T$. It is not exactly the case we considered since now the migration depends on time, but it helps to understand.
We thus consider the following situation:
\ben
\item On patch 1 the population grows in an exponential way
$ \frac{dx}{dt} = \alpha x$
during a duration $ T$. 
\item A part of the population of patch 1 is transferred on  patch 2 according to two different modalities:
	\ben
	\item \textbf{Discrete mode} : At the end of the growth period $T$ a part $d$ ($0 \leq d \leq 1$) is transferred from patch 1 on the patch 2. 
	\item \textbf{Continuous mode} : Continuously ( $\frac{dx}{dt} = -m x$ ) a part of the population of the patch  1 is transferred on the patch 2.
	\fen
\item The patch  2 is a ''neutral'' patch where the population remains constant.
\fen
Let us suppose that one seeks to maximize the size of the population on patch  2 at the end of T time units because we know that in the following period this patch will be the source.

In the discrete mode, it is obvious that the most efficient thing to do is to transfer with $d = 1$ the totality of the population of patch 1, i.e. $\emat^{\alpha T}x(0)$ to the patch 2.
On the other hand, in the continuous mode, if it allows  to (almost) empty completely the patch 1 (by taking $m$ infinitely large), does not allow in this case to obtain a good result on the patch 2 because all the population subtracted at the beginning of the period will not undergo any more growth on patch 2 ; conversely, if $m$ is very small, there is a strong growth of the biomass on patch  1 of which only a very small part is transferred to patch  2; thus {\em for the continuous mode, small and large values of $m$ are inefficient to obtain a large population on patch 2} and intermediate values are to be considered. 

Let us complicate matters a little by imagining now, in the discrete case, two successive periods $T_1$ and $T_2$ where, to mimic a positive correlation in time, the population of patch n°1 grows exponentially at the rates $\alpha_1$ and then $\alpha_2$, but where {\em the same proportion} $d$ of the population is transferred at times $T_1$ and $T_1+T_2$ from patch 1 to patch  2. In total, the quantity transferred from $1$ to $2$ is 
$$d \emat^{\alpha_1 T_1} x(0)+ d \emat^{\alpha_2 T_2}(1-d) \emat^{\alpha_1 T_1} x(0)$$
and this time, it will  be a value of $d$ strictly between $0$ and $1$ that maximizes it since the formula is quadratic in $d$.

 It is this fundamental difference between continuous migration (for instance, planktonic microorganisms drifting between two reefs) and (near) instantaneous dispersal over two distinct sites (such as seed dispersal at flowering) that makes comparison of the inflation phenomenon in the discrete and continuous cases intricate. We believe that a complete understanding of similarities and differences of the continuous and the discrete time models  is beyond the scope of this discussion and requires further investigation.
 
    \section{Conclusion}
 
  The major limitation of our approach is obviously that it is quite specific to the two-patch situation. The general mathematical results of Katriel (continuous time) and Schreiber (discrete time) ensure that from a qualitative point of view the phenomenon of inflation (or DIG (according to Katriel's terminology) is present on any system of $N$ patches. In the case of two patches we have shown how the phenomenon can be can be accurately quantified. For more than two patches we do not know, for the moment, how to proceed but it seems likely that different assumptions about the topology of the sites (island-continent, stepping-stones, homogeneous dispersal etc...) will be necessary. In addition, as noted above, a better understanding of the differences and similarities of discrete and continuous time models needs to be worked on, both for reasons of mathematical aesthetics but more importantly for their biological isignificance.
  
  In classical mechanics, the harmonic oscillator, i.e. the linear differential equation of the second order with a constant coefficient,   plays a major organizing role. It is a simple mathematical object that can be taught and understood very early (at the end of high school) and that opens the door to various fields : the theory of nonlinear oscillators and endogenous oscillations, the theory of forced oscillations and  the phenomenon of resonance (that we could call inflation) etc.
  In a stimulating essay, {\em Ecological orbits, how planets move and populations grow}, Lev Ginzburg and M. Colyvan  \cite{GIN04} defend the idea that exponential growth is, for population dynamics, the equivalent of the principle of inertia in classical mechanics: a population whose growth rate is not limited grows exponentially, just as a material body subjected to no force keeps the same speed. 

If their vision is correct  the $(±1)$ model, in spite of its total unrealism but thanks to its mathematical simplicity, is a basic brick in the understanding of the mechanisms that govern the growth of a meta-population on various connected  patches, with temporally varying growth rates.

 \section*{Acknowledgements}

We warmly thank the editor of T.P.B.  S. Schreiber who provided us interesting references and the two anonymous reviewers of T.P.B. whose excellent reviews have produced a significantly improved version of our article.

\newpage
\section*{Appendix }
\appendix

\section{From $(x_1,x_2)$ to $(U,V)$} \label{chtvariables}

Consider :
\beq \label{cht1}
\begin{array}{lcl} 
\displaystyle \frac{dx_1}{dt} &=& r_1(t) x_1 +m( x_2-x_1) \\[6pt]
\displaystyle \frac{dx_2}{dt}& =&r_2(t) x_2 +m(x_1-x_2)
\end{array} 
\feq
which is the same a :
\beq \label{cht2}
\begin{array}{lcl} 
\displaystyle \frac{dx_1}{dt} &=&\left( r_1(t)  +m( x_2/x_1 -1)\right) x_1 \\[6pt]
\displaystyle \frac{dx_2}{dt} &=&\left( r_2(t)  +m(   x_1/x_2 -1)\right) x_2
\end{array} 
\feq
Let :
$$\xi_1 = \ln(x_1) \Leftrightarrow x_1 = \emat^{\xi_1}\quad \quad  \xi_2 = \ln(x_2) \Leftrightarrow  x_2 = \emat^{\xi_2}$$
Then one has :
\beq \label{cht3}
\begin{array}{lcl} 
\displaystyle \frac{d\xi_1}{dt} =\frac{1}{x_1}\frac{dx_1}{dt} &=& r_1(t)  +m( \emat^{\xi_2-\xi_1}- 1) \\[6pt]
\displaystyle \frac{d\xi_2}{dt}  =\frac{1}{x_2}\frac{dx_2}{dt} & =&r_2(t) +m(\emat^{\xi_1-\xi_2}-1) 
\end{array} 
\feq
Let :
$$ U = \frac{\xi_1 + \xi_2 }{2} \quad \quad V = \frac{\xi_1 - \xi_2 }{2}$$
Then, adding and substracting the above equations one gets :
 \beq \label{cht4}
\begin{array}{lcl}
\displaystyle  \frac{dU}{dt} &=&\displaystyle \frac{r_1(t)+r_2(t)}{2} +m \,( \ch(2V) -1)  \\[8pt]
\displaystyle  \frac{dV}{dt} &=&\displaystyle  \frac{r_1(t)- r_2(t)}{2}-m\,\sh(2V)
 \end{array} 
\feq

\section{The  system $\Sigma(\eps,m,T) $ for large $m$ \label{tycho}}
Consider the system :
$$
\begin{array}{lcl}
\displaystyle  \frac{dx_1}{dt} &=& (+u(t)-\varepsilon)x_1+m(x_2-x_1)\\[8pt]
\displaystyle  \frac{dx_2}{dt} &=& (-u(t)-\varepsilon)x_2+m(x_1-x_2) 
 \end{array} 
$$
and put :
$$S = x_1+x_2\quad \quad D = x_1 - x_2$$
One has :
$$
\begin{array}{lcl}
\displaystyle  \frac{dS}{dt} &=& (+u(t)-\varepsilon)x_1+m(x_2-x_1)+ (-u(t)-\varepsilon)x_2+m(x_1-x_2)= u(t)D-\varepsilon S\\[8pt]
\displaystyle  \frac{dD}{dt} &=& \displaystyle u(t) -2(m+\varepsilon)D = -2(m+\varepsilon)\left(D - \frac{u(t)}{2(m+\varepsilon)}\right)
 \end{array} 
$$
From  Tychonov theorem (\cite{TYK52,LOB98}), when $2(m+\varepsilon) \rightarrow \infty$ the solution to this system (called a slow-fast system), after a  short transient, tends to : 
$$ S(t) = S(0) \emat^{-\varepsilon t}\quad \quad D(t) = 0$$
Thus, in the variables $(x_1,x_2)$ one has :
$$x_1(t) \approx x_2(t)  \approx \frac{x_1(0) + x_2(0)}{2} \emat^{-\varepsilon t}$$
This is also called by physicists, the method of the {\em quasi stationary state}.

\section{The switched system $F(m,T)$ has a single globally stable periodic orbit}\label{apendiceperiodic}
We consider the one dimensional switched system :
 \beq \label{commutdim1}
  F(m,T)\quad \quad
\begin{array}{lcl}
\displaystyle  \frac{dV}{dt} &=&\displaystyle  \frac{1+u(t)}{2}F^+_m(V)+ \frac{1-u(t)}{2} F^-_m(V)  
 \end{array} 
\feq
with : 
$$F^+_m (V) = +1-m\,\sh(2V) \quad \quad F^-_m (V) = - 1-m\,\sh(2V)$$
and :
$$ V_m^+ =\frac{1}{2}  \sh^{-1}(+1/m)  \quad \quad \quad  V_m^- = \frac{1}{2}   \sh^{-1}(-1/m)$$

Since $F^+_m(V)$ is continuous, differentiable and  such that $(V-V^+_m)F^+_m(V) < 0$ except for $V = V^+_m$,
from the elementary theory of differential equations we know that,  if we denote by $\varphi^+_t(v)$ the unique  solution of :
$$\frac{dV}{dt} = F^+_m(V) \quad \quad \quad  V(0) = v$$
then $\varphi^+_t(v)$ is defined for every positive $t$ and the mapping $v \mapsto \varphi^+_t(v)$ is continuous and differentiable.
{\lemme Let $T >0$. The mapping $V \mapsto\varphi^+_T(v)$ is a continuous mapping, strictly increasing, from $[V^-_m\,,\,V^+_m]$ into $[V^-_m\,,\,V^+_m]$ ; moreover its derivative is strictly smaller than 1. }\\
\textbf{Proof } Assume that $v \mapsto \varphi^+_T(v)$ is not strictly increasing. Then it exists $v_1 < v_2$ such that $\varphi^+_t(v_2) \leq \varphi^+_t(v_1)$ ; and, by the way, some $t \leq T$ for which $\varphi^+_t(v_1) = \varphi^+_t(v_2)$ and, thus, two solutions starting from different initial conditions meet at some point. This contradicts the uniqueness of solutions. The derivative of $V \mapsto \varphi^+_T(v)$ at the point $v_o$ is obtained by integrating the linearized equation along the trajectory $t \mapsto\varphi^+_t(v_o)$ up to time $T$, that is to say :
$$ \frac{d\delta \varphi^+_t}{dt} = D F^+_m(\varphi^+_t(v_0))\delta \varphi^+_t \quad \quad \quad \delta \varphi^+_0 = 1$$
where $DF^+_m(V)$ is the derivative of $F^+_m(V)$ at point $V$.
$$ \delta \varphi^+_T = \mathrm{exp} \left(\int_0^T DF^+_m(\ \varphi^+_t(v_o))dt \right)dt$$
One has :
$$\int_0^T DF^+_m(\varphi^+_t(v_0))dt = \int_{v_o}^{\varphi^+_t(v_0)}\frac{DF^+_m(V)}{F^+_m(V)} dV = \ln(F^+_m(\varphi^+_t(v_0))-\ln(F^+_m(v_0))$$
Since the function $F^+_m(V)$ is decreasing and $V\varphi^+_t(v_0) > v_o$ the integral is negative and thus $\delta \varphi^+_T < 1$.\\$\Box$

For the same reasons we have the following  lemma \ref{b2} : recall that We denote by $\varphi^-_t(v)$ the unique  solution of :
$$\frac{dV}{dt} = F^-_m(V) \quad \quad \quad  V(0) = v$$
then $V\varphi^-_t(v)$ is defined for every positive $t$ and the mapping $v \mapsto \varphi^-_t(v)$ is continuous and differentiable.
{\lemme \label{b2}Let $T >0$. The mapping $V \mapsto \varphi^-_T(v)$ is a continuous mapping, strictly increasing, from $[V^-_m\,,\,V^+_m]$ into $[V^-_m\,,\,V^+_m]$ ; moreover its derivative is strictly smaller than 1. }\\
\begin{figure}
  \begin{center}
 \includegraphics[width=1\textwidth]{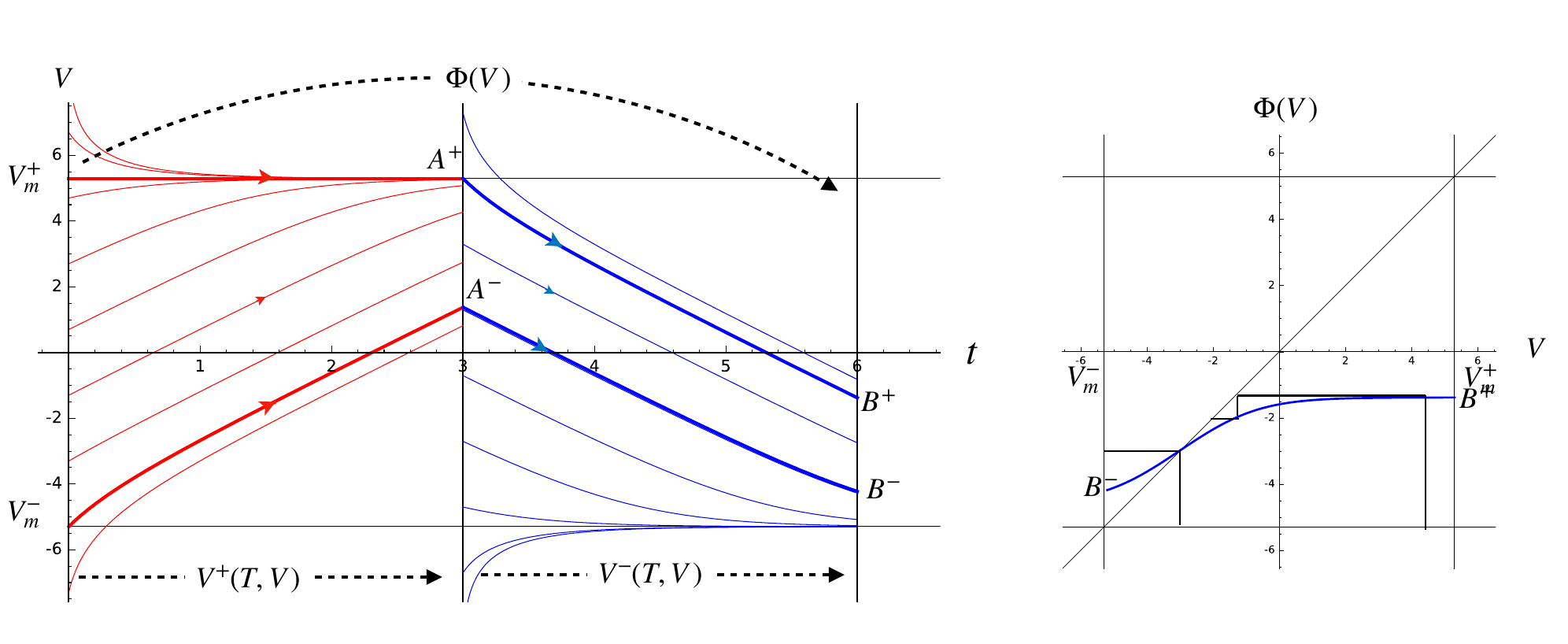}
 \caption{ The switched system $F(m,T)$: $m=0.01\;;\; T = 3$. The segment $[A^-,A^+]$ is the image of $[V^-_m,V^+_m]$ by the mapping $V \mapsto V^+(T,V)$, and $[B^-,B^+]$ is the image of $[A^-_m,A^+_m]$ by the mapping $V \mapsto V^-(T,V)$ }\label{pointfixe}
 \end{center}
 \end{figure}

Now consider ''period-map'', that is the composite mapping $ v \mapsto \Phi(V)= \varphi^-_T \circ \varphi^+_T(v)$ from $[V^-_m\,,\,V^+_m]$ into $[V^-_m\,,\,V^+_m]$  (see figure \ref{pointfixe}, left); from the preceding  lemmas it turns out that  it is strictly increasing, with $\Phi'(V) < 1$ such that $V^-_m <\Phi(V^-_m)$ and  $\Phi(V^+_m) < V^+_m $. From elementary calculus  the discrete dynamical system defined by :
$$ V(n+1, V_o ) = \Phi(V(n,V_o)), \quad \quad V(0,V_o) = V_o$$
has a unique equilibrium $V_e$ , i.e. the unique solution of $\Phi(V) = V$, this equilibrium is globally asymptotically stable  (see figure \ref{pointfixe}, right). Since $ \Phi(V(n,V_o) = V(n2T, V_o) $,  where $V(n2T, V_o)$ is the solution of the switched system (\ref{commutdim1}) we have proved :
{\proposition \label{apandiceprop1} The switched system $F(m,T)$ has a unique periodic solution, denoted $P_{m,T}(t)$,  globally asymptotically stable which oscillates between two values $P^-_{m,T}$, and $P^+_{m,T}$ contained in the interval $[V_m^-,\; V_m^+]$ ;  $P^-_{m,T}= - P^+_{m,T}$ and 
the function $T \mapsto P^+_{m,T}$ is an increasing function of $T$ which tends to $V_m^+$ when $T$ tends to infinity.}\\\\
\begin{figure}
  \begin{center}
 \includegraphics[width=1\textwidth]{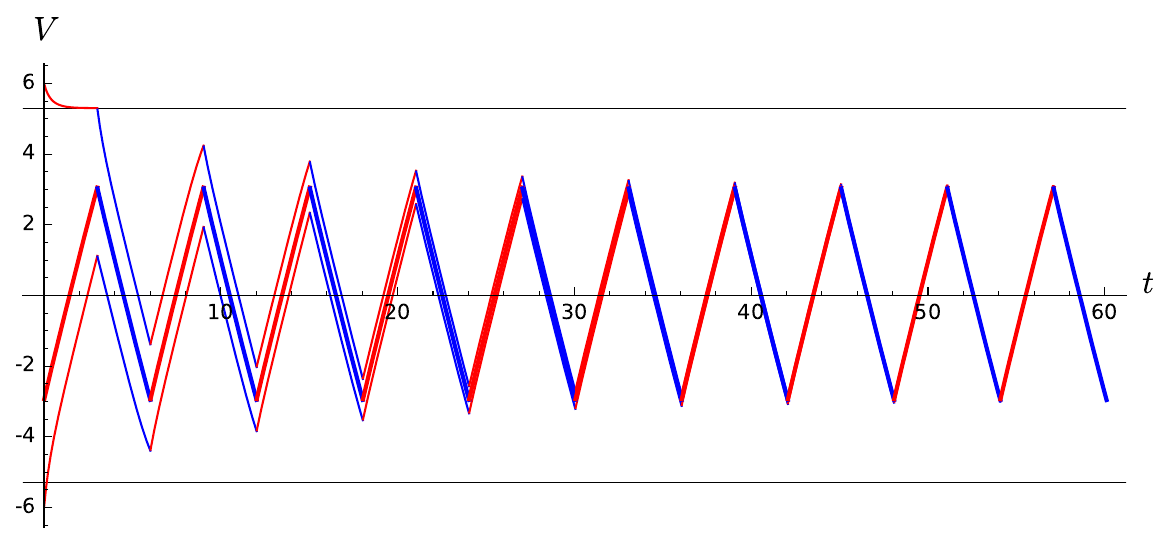}
 \caption{ The switched system $F(m,T)$: $m=0.01\;;\; T = 3$. Solutions converge to a unique periodic orbit. }\label{solper}
 \end{center}
 \end{figure}
 The solutions of (\ref{commutdim1}) are explicitly computable as we show now.\\
 On $[0,T]$ one has :
  \beq \label{}
\begin{array}{lcl}
\displaystyle  \frac{dV}{dt} &=&1-m\,\sh(2V))\quad \quad V(0) = V_0
 \end{array} 
\feq
thus $dt = \frac{dV}{1-m\,\sh(2V))}$ and by the way :
   \beq \label{}
\begin{array}{lcl}
t &=& \displaystyle \int_{V_o}^{V(t)}  \frac{dV}{1-m\,\sh(2V)} 
 \end{array} 
\feq
Since the function that we have to integrate is a rational fraction with respect to $\emat^V$ we can integrate it explicitly (by hand or with the help of some formal software)  and the result is :
$$
\int\frac{dV}{1-m\sinh(2V))}=\frac{1}{A} \tanh^{-1}\left(\frac{\tanh(V)+m}{A}\right)
$$
where $$A = \sqrt{1+m^2}$$
from which we have $V$ as a function of $t$.

The periodic solution oscillate between $-P^+_{m,T}$ and $+P^+_{m,T}$ solutions of the equation :
\beq \label{TV}
T=\int_{-P^+_{m,T}}^{P^+_{m,T}}\frac{dV}{1-m\sinh(2V))}
\feq 
 and thus $P^+_{m,T}$ is a solution of the equation :
$$
\tanh^{-1}\left(\frac{\tanh(P^+_{m,T})+m}{A}\right)-
\tanh^{-1}\left(\frac{\tanh(-P^+_{m,T})+m}{A}\right)=
TA
$$
Thus, if we put $x = \tanh(P^+_{m,T})$, we are searching for the solutions of the equation
$$
\tanh^{-1}\left(\frac{x+m}{A}\right)+
\tanh^{-1}\left(\frac{x-m}{A}\right)=TA
$$
From the formula $\tanh^{-1}(a)+\tanh^{-1}(b)=\tanh^{-1}\left(\frac{a+b}{1+ab}\right)$,  one obtains the equation

$$TA = \tanh^{-1}\left( \frac{\frac{x+m}{A} +\frac{x-m}{A} }{1+ \frac{x+m}{A} \frac{x-m}{A} }   \right)$$
$$
TA = \tanh^{-1}\left(\frac{2Ax}{A^2-m^2+x^2}\right)
$$
and :
$$
x^2\tanh(TA)-2Ax+\tanh(TA)(A^2-m^2)=0
$$
but since  $A = \sqrt{1+m^2}$ one has 
$$
x^2\tanh(TA)-2Ax+\tanh(TA)=0
$$
Put :
$$ B = \tanh(TA)$$
This equation admits two solutions :
$$
x=\frac{A-\sqrt{A^2-B^2}}{B},
\quad
x=\frac{A+\sqrt{A^2-B^2}}{B},
$$
The second solution is not acceptable since it is grater than $1$.We have :
\begin{equation}\label{thV*}
\tanh(P^+_{m,T})=\frac{A-\sqrt{A^2-B^2}}{B},
\quad\mbox{with }B=\tanh(TA)
\end{equation}
Thus we have proved the :
\begin{proposition} \label{valmax}
The maximum $P^+_{m,T}$ (respectively the minimum $P^-_{m,T}= -P^+_{m,T}$) of the periodic solution of (\ref{commutdim1})  is given by :
\beq
\boxed{
\begin{array}{lcl}
P^+_{m,T}=\displaystyle \tanh^{-1} \left(\frac{A-\sqrt{A^2-B^2}}{B}\right)\\[12pt]
\mbox{with} \;A = \sqrt{1+m^2} \,\mbox{and} \; B = \tanh(TA)
\end{array} }
\feq
\end{proposition} 

\section{Qualitative properties of $\Delta(\eps,m,T)$}\label{appendicedelta}
Recall that : 
$$\Delta(\eps,m,T) = \frac{1}{2T} \int_0^{2T} \varphi(P_{m,T}(s))ds \quad \mathrm{ with}\quad \varphi(V) = m\ch(2V)-m-\eps $$
We prove :
{\proposition \label{propgenerale2} \textbf{Qualitative properties of $\Delta(\eps,m,T)$}
\ben
\item For small  $m$, $\Delta(\eps,m,T) < 0$
\item For large $m$, $\Delta(\eps,m,T) > 0$
\item For fixed $\eps>0$ and $m <\frac{1-\eps^2}{2\eps}$  there exists a threshold $T^*(\eps,m)$ such that for $T < T^*(\eps,m)$,  $\Delta(\eps,m,T) < 0$  and $\Delta(\eps,m,T) > 0$  for $T>T^*(\eps,m)$ 
\item For every $\eps > 0$, the minimum of $T^*(\eps,m)$ over $m$ is strictly positive. In other words there exists a threshold $T^{**}> 0$ such that for $T < T^{**}$, for all values of $m$, $\Delta(\eps, m, T) < 0$ : there is no inflation.
\fen }

\noindent \textbf{Proof of 1)} Since in the interval $[V^-_m,\,V^+_m]$, where the periodic solution lives, one has $|\frac{dV}{dt}| < 1 $ we know that $P^+_{m,T} < T$. Since $\displaystyle \lim_{m\rightarrow 0} A^+_{\eps,m} = \lim_{m\rightarrow 0} \ch^{-1}\left(1+\frac{\eps}{m}\right) = +\infty$, for small enough $m$, the periodic solution $ P_{m,T}(t)$ lives in the interval $]A^-_{\eps,m} ,\;A^+_{\eps,m} [$ where the function $\varphi$ is strictly negative an hence so is $\Delta(m,T)$.\\

\noindent \textbf{Proof of 2)} Given $\eps > 0$ the relative positions of $V^+_m = \sh^{-1}(\frac{1}{m})$ and $A^+_{\eps,m}=\ch^{-1}\left(1+\frac{\eps}{m}\right) $ depends on $m$. One easily compute  that 
\beq \label{ineg}  
A^+_{\eps,m} < V^+_m \Longleftrightarrow m < \frac{1-\eps^2}{2\eps}
\feq
hence, if $m >  \frac{1-\eps^2}{2\eps} $ one has $ [ V^-_m,\, V^+_m]\subset [A^-_{\eps,m} ,\, A^+_{\eps,m} ]$ and $\varphi(P_m(t))$ is always negative.\\

\noindent \textbf{Proof of 3)} One has :
$$ \int_0^{2T} \varphi(P_{m,T}(s))ds = 2\int_{-P^+_{m,T} }^{P^+_{m,T} }\frac{m\,\ch(2V)-m-1}{1-m\,\sh(2V)}dV$$
The conclusion follows from the fact that $T \mapsto P^+_{m,T}$ is an increasing function of $T$ such that :
$$\lim_{T \to 0} P^+_{m,T} = 0\quad \mathrm{and} \quad  \lim_{T \to \infty } P^+_{m,T} = V^+_m = \frac{1}{2} \sh^{-1}\left(\frac{1}{m}\right) $$
and :
$$ \lim_{V \to V^+_m }  \int_{_V }^{+V} \frac{m\,\ch(2V)-m-1}{1-m\,\sh(2V)}dV = +\infty$$

\noindent \textbf{Proof of 4)} Let $\eps > 0$ given, consider $m_o = \frac{1-\eps^2}{2\eps}$ and set  $ T_o = A^+_{\eps,m_o} = V^+_{m_o}$. Let $T < T_o$. For $m>m_o$ we already know that $\Delta(\eps,m,T) < 0$. If $m < m_o$ one has  $A^+_{\eps,m} > T_o > T$ and, by the way, $\varphi(P_{m,T}(t))$ is always negative. Thus for $ T < \sh^{-1}\left(\frac{2\eps}{1-\eps^2}\right)$, whatever the value of $m$, one has 
$\Delta(\eps,m,T) < 0$ which proves 4).

\section{Explicit formula for  $\Delta(\eps,m,T)$}\label{appendicedeltaexplicite}

We consider the periodic solution  $P_{m,T}(t)$ to $F(m,T)$.  We are interested by the sign of :
\beq \label{solperannexe} 
\Delta(\eps,m,T) =   \frac{1}{T} \int_{P^-_{m,T}}^{P^+_{m,T}} \frac{m\ch(2V)-m-\varepsilon}{1-m\sh(2V)}dV
\feq

From the formula (\ref{solperannexe}) and proposition \ref{valmax} we can  deduce an explicit formula for $\Delta(\epsilon,m,T)$. First, 
if we use formula (\ref{TV}) in the definition of $\Delta(\eps,m,T)$ we get :
\beq
\Delta(\eps,m,T)=\frac{1}{T}\int_{-P^+_{m,T}}^{P^+_{m,T}}\frac{m\cosh(2V)dV}{1-m\sinh(2V)}-(m+\varepsilon)
\feq 
Since  $\frac{d}{dV}\sinh(V)=\cosh(V)$, one can explicitly compute the integral to get :
\beq
\int\frac{m\cosh(2V)dV}{1-m\sinh(2V)}=- \frac{\ln(1-m\sinh(2V))}{2}
\feq
and, by the way :
\begin{equation}\label{Delta}
\Delta(\eps,m,T)=\frac{1}{2T}\ln\frac{1+m\sinh(2P^+_{m,T})}{1-m\sinh(2P^+_{m,T})}-(m+\varepsilon)
\end{equation}
Using the formula  $\sinh(a)=\frac{2\tanh(a/2)}{1-\tanh^2(a/2)}$, from (\ref{thV*}) one gets :

$$\sinh(2P^+_{m,T}) = \frac{2 \tanh(P^+_{m,T})}{1 - \tanh^2(P^+_{m,T})}$$

Now, replacing $P^+_{m,T}$ by its value given by prop \ref{valmax} :

$$\tanh(P^+_{m,T })=\displaystyle  \left(\frac{A-\sqrt{A^2-B^2}}{B}\right)$$

$$
\sinh(2P^+_{m,T})=\frac{B(A-\sqrt{A^2-B^2})}{B^2-A^2+A\sqrt{A^2-B^2}}
$$
If we replace in (\ref{Delta}) we have :
$$\Delta(\eps,m,T)=\frac{1}{2T}\ln\frac{ B^2-A^2+A\sqrt{A^2-B^2}+mB(A-\sqrt{A^2-B^2}) }
{  B^2-A^2+A\sqrt{A^2-B^2}-mB(A-\sqrt{A^2-B^2})}-(m+\varepsilon)$$
and, after a multiplication by the conjugate quantity of the denominator one have the more simple expression :
$$
\Delta(\eps,m,T)=\frac{1}{2T}\ln\frac{A^2-B^2+m^2B^2+
2mB\sqrt{A^2-B^2}}{A^2-B^2-m^2B^2}
-(m+\varepsilon)
$$
Using $A^2=1+m^2$ and  $B=\tanh(TA)=\frac{e^{2TA}-1}{e^{2TA}+1}$, one gets :
\begin{equation}
\label{Delta2}
\Delta(\eps,m,T)=\frac{1}{2T}\ln\frac{m^2b^4+2b^2+m^2+
m(b^2-1)\sqrt{C}}{2(1+m^2)b^2}
-(m+\varepsilon)
\end{equation}
with   $b=e^{T\sqrt{1+m^2}}$ and
$C=m^2b^4+2m^2b^2+4b^2+m^2$.

\section{Asymptotics of $\Delta(\eps,m,T)$ for large $T$ } \label{asymptotics}
We are looking for solutions of $m \mapsto \Delta(\eps, m, T) = 0$ which are exponentially small with respect to $T$, that is to say for $x$ solutions of :
\beq 
\Delta(\eps, \emat^{xT}, T) = 0\quad \quad x < 0
\feq
We use Landau notation $\mathrm{o}$ for any quantity that tends to $\mathrm{o}$ when $t$ tends to $\infty$. From (\ref{Delta2}), 

$\Delta(\varepsilon,m,T)=0$ is equivalent to :
\beq\label{asy3}
\disp \frac{m^2b^4+2b^2+m^2+ m(b^2-1)\sqrt{C}}{2(1+m^2)b^2 } = \emat^{2(m+\varepsilon)T }
\feq 
with   
\beq\label{asy2}
\disp b=e^{T\sqrt{1+m^2}}, \quad C=m^2b^4+2m^2b^2+4b^2+m^2, \quad m = \emat^{xT}
\feq 
From (\ref{asy2}) one have :
\beq
\disp mb = \emat^{Tx}\emat^{T\sqrt{1+\emat^{2Tx}}}= \emat^{T(1+x+\frac{1}{2} \emat^{2Tx} (1+\po))}
\feq 
since for $x < 0$ we have $T \emat^{2Tx} = \po$ we deduce $mb = \emat^{T(1+x +\po})$ which tends to $\infty$ as long as $x > -1$.
from which we deduce that as long as $x > -1$   :
\beq
m^2b^4+2b^2+m^2 = m^2b^4(1+\po(1)) \quad \quad m(b^2-1)\sqrt{C} = m^2b^4(1+\po(1))
\feq 
which introduced in \eqref{asy3} gives :
\beq
m^2b^2(1+\po(1)) = \emat^{2T(1+x)(1+\po(1))} =  \emat^{2 T(\eps+\po(1))}
\feq 
from which we deduce :
\beq
2T(1+x)(1+\po(1)) =  2 T(\eps+\po(1)) \Longrightarrow  x = -(1-\eps)+ \po(1)
\feq
which is the evaluation of proposition \ref{approxseuil}.

\section{Connection between $\Delta(\eps,m,T)$ and $\sigma(\eps,m,T)$.}\label{sigmadelta}
Let $(x_1(t), x_2(t))$ be any solution of $\Sigma(\eps,m,T)$ ; let  $U(t) =\ln(\sqrt{x_1(t)x_2(t)})$ and $V(t) =  \ln(\sqrt{x_1(t)/x_2(t)})$.Then $V(t)$ is a solution of $F(m,T)$ and since the periodic solution of $F(m,T)$ is globally asymptotically stable $V(t)$ converges to $P_{m,T}(t)$,  thus :
\beq \label{lien1}
\begin{array}{l}
\displaystyle \lim_{n\rightarrow +\infty} U((n+1)2T)-U(n2T) =\cdots \\ [8pt]
\displaystyle \quad \quad \quad \lim_{n \rightarrow +\infty} \int_{n2T} ^{(n+1) 2T} 2(m\ch(V(s))-m-\eps)ds = \cdots \\ [8pt]
\displaystyle\quad \quad \quad \quad\quad  \quad\quad \int_0^{2T} m\ch(P_{m,T}(s))-m-\eps ds =2T \Delta(\eps,m,T)
\end{array} 
\feq \label{lien1}
and hence :
\beq
\Delta(\eps,m,T) = \frac{1}{2T}\displaystyle \lim_{n\rightarrow +\infty} \ln\left(\sqrt{ \frac{x_1((n+1)2T)x_2((n+1)2T)}{x_1(n2T)x_2(n2T)}  } \right)
\feq
Now, choose $(x_1(0), x_2(0)) = Z_1$ where $Z_1$ is the positive eigenvector of $M(\eps,m,T)$ associated with $\lambda_1$ (note that $M(\eps,m,T)$ has positive entries). Then, for all $n \geq 0$, $x_1((n+1)2T) \approx \lambda_1 x_1(n2T)$ and $x_2((n+1)2T) =\approx \lambda_1 x_2(n2T)$, thus
\[
\sqrt{ \frac{x_1((n+1)2T)x_2((n+1)2T)}{x_1(n2T)x_2(n2T)} }\approx  \lambda_1^2
\]
and thus $\Delta(\eps,m,T) = \frac{1}{2T} \ln( \lambda_1),$ which we wanted to prove.

We can also prove this equality directly from the explicit formulas of $\Delta$ and $\lambda_1$. The value of $\lambda_1$ given by Maple is :
$$\lambda_1=
\frac{e^{-2(m+\varepsilon)T}}{2A^2b^2}\left(
m^2b^4+2b^2+m^2+\sqrt{C_1}\right)$$
with 
$$C_1=
b^8m^4+4b^6m^2-2b^4m^4-8b^4m^2+
4b^2m^2+m^4
$$
$$\frac{1}{2T} \ln(\lambda_1)=\ln\frac{m^2b^4+2b^2+m^2+\sqrt{C_1}}{2A^2b^2}-(m+\varepsilon)$$

which is the value of $\Delta$ given by the proposition  \eqref{explicitdelta} since one has :
$$C_1=m^2(b^2-1)^2
(m^2b^4+2m^2b^2+4b^2+m^2)$$.

\section{Limit in distribution of $\hat V_n$ in the stochastic $(\pm 1)$ model}
 \label{app:lemV2k}

We prove Lemma \ref{lem:V2k}. The first point follows from the same observation as in the deterministic case, that the solution of $\bm{F}(m, \mu)$ are trapped in the interval $[V_m^-, V_m^+]$. Now, 
note that, for all $n \geq 0$;
\[
\hat V_{2n+2} = \varphi^-_{S_{2n+2}} \circ \varphi^+_{S_{2n}} ( \hat V_{2n} )
\]
Therefore, since $(S_n)_{n \geq 0}$ is a sequence of i.i.d. variables, the sequence $(\hat V_{2n})_{n \geq 0}$ is a Markov chain. Let us denote $\hat V_{2n}^v$ 
for the position of the chain after $n$ steps, whenever $\hat V_0 = v$. Since (see \ref{apendiceperiodic}}) $\delta \varphi^+_t \leq e^{ - 2 m t}$ and  
$\delta \varphi^-_t \leq e^{ - 2 m t}$, one has, for all $v, v' \in [V_m^-, V_m^+]$, 
\[
| \hat V_2^v - \hat V_2^{v'} | \leq e^{ - 2 m( S_2 + S_1) } | v - v'|.
\] 
Thus, using that $\mathbb{P}(S_1 = 0) = 0$, and therefore, $\mathbb{E}( e^{ - 2 m( S_2 + S_1) } ) < 1$, 
the Markov chain $\hat V_{2n}$ is contracting for the Wasserstein distance on the complete space $[V_m^-, V_m^+]$. As such, 
it admits a unique stationary distribution $\nu_{\infty}$, which is the law of a variable $V_{\infty}^{-}$, and $\hat V_{2n}$ converges 
geometrically fast in distribution to $V_{\infty}^-$. This proves the second point of Lemma \ref{lem:V2k}. As for the third point, this is a consequence of Birkhoff's ergodic theorem.

 
 \section{Existence of the growth rate in the random switching time case}
 \label{App:Delta_random_switch}
By \eqref{UdeR}, the asymptotic behaviour of $\frac{U(t)}{t}$ is given by those of
\[
 \frac{1}{t} \int_0^t m\, \left(\ch(2V(s)) -1 \right)ds  -\varepsilon := H(t) - \varepsilon.
\]
Let $t \geq 0$, there exists $n = N_t \in \mathbb{N}$ (random) such that $T_{2n} \leq t < T_{2n+2}$. This means that a time $t$, at least $2 n$ and at most $2n+1$ jumps have occurred. Therefore, we can write 
\begin{align*}
H(t) & = \frac{1}{t} \sum_{k=0}^{2N_t - 1} \int_{T_k}^{T_{k+1}}  m\, \left(\ch(2V(s)) -1 \right) ds +R(t) \\
& = \frac{1}{t} \sum_{k=0}^{N_t - 1} \left( \int_{T_{2k}}^{T_{2k+1}}  m\, \left(\ch(2V(s)) -1 \right) ds + \int_{T_{2k+1}}^{T_{2k+2}}  m\, \left(\ch(2V(s)) -1 \right) ds\right) +R(t),
\end{align*}
where $R(t)$ is a rest term given by
\[
R(t) = \frac{1}{t} \int_{T_{2N_t}}^t m\, \left(\ch(2V(s)) -1 \right)ds.
\]
Note that $R(t) \leq C \frac{S_{2n+1}}{t}$, for some constant $C$, and thus $R(t) \to 0$ as $t \to \infty$. We now use Lemma \ref{lem:V2k} to give the asymptotic behaviour of the term in the sum.
First, note that on a interval $[T_{2k}, T_{2k+1})$, in \eqref{eq:VR}, we are integrating system  $\Sigma^{+}( \varepsilon, m)$, with initial condition $V_{T_{2k}}$ and thus, for $s \in [T_{2k}, T_{2k+1})$, one has $V(s) = \varphi^+_{s - T_{2k}}(\hat V_{2k})$. Similarly, if $s \in [T_{2k+1}, T_{2k+2})$, one has $V(s) = \varphi^-_{s - T_{2k+1}}(\hat V_{2k+1})$, so that the first term in $H(t)$ can be rewritten as
\begin{equation}
\frac{1}{t} \sum_{k=0}^{N_t - 1} \left( \int_{0}^{S_{2k+1}}  m\, \left(\ch(2\varphi^+_{s}(\hat V_{2k+1})) -1 \right) ds + \int_{0}^{S_{2k+2}}  m\, \left(\ch(2\varphi^-_{s}(\hat V_{2k+1})) -1 \right) ds\right)
\label{eq:H}
\end{equation}
Next, for $(v,s) \in \mathbb{R} \times \mathbb{R}_+$, set 
\[
f^+(v,s) = \int_0^s  m\, \left(\ch(2\varphi^+_{r}(v)) -1 \right) dr
\]
and
\[
f^-(v,s) = \int_0^s  m\, \left(\ch(2\varphi^-_{r}(v)) -1 \right) dr
\]
Thus, 
\[
H(t) = \frac{1}{t}  \sum_{k=0}^{N_t - 1}  f^+( \hat{V}_{2k}, S_{2k+1}) + \frac{1}{t}  \sum_{k=0}^{N_t - 1}  f^-( \hat{V}_{2k+1}, S_{2k+2}) + R(t)
\]
Now, we have
\[
 \frac{1}{t}  \sum_{k=0}^{N_t - 1}  f^+( \hat{V}_{2k}, S_{2k+1})  =  \frac{N_t - 1}{t}  \frac{1}{N_t-1}   \sum_{k=0}^{N_t - 1}  f^+( \hat{V}_{2k}, S_{2k+1}) 
\]
Now, classical renewal theorem\footnote{see e.g. \url{https://en.wikipedia.org/wiki/Renewal_theory}} implies that, with probability one,
\[
\lim_{t \to \infty}  \frac{N_t - 1}{t} = \frac{1}{2 \mathbb{E}(S)},
\]
and in particular, $N_t \to \infty$ as $t \to \infty$. Hence, the third assertion of Lemma \ref{lem:V2k} implies that 
\[
\lim_{t \to \infty}  \frac{1}{N_t-1}   \sum_{k=0}^{N_t - 1}  f^+( \hat{V}_{2k}, S_{2k+1})  = \mathbb{E}[f^+( V_{\infty}^-,S)],
\]
and as a consequence, 
\[
\lim_{t \to \infty} \frac{1}{t}  \sum_{k=0}^{N_t - 1}  f^+( \hat{V}_{2k}, S_{2k+1})  = \mathbb{E}[f^+( V_{\infty}^-,S)].
\]
The same reasoning can be done for the term with $f^-$ and yields the expected result.

\section{Density dependent model} \label{densitedependant}
\subsection*{The deterministic case}
{\proposition \label{persistance1} If the parameters $(\eps,m,T)$ are such that the system $D(\varepsilon,0,m,T) =\Sigma(\varepsilon,m,T)$ is stable, then the solutions of $D(\varepsilon,\alpha,m,T)$ tend to $0$.}\\
\textbf{Proof.} We denote by $(x_1(t,x_{1_0},x_{2_0}),x_2(t, x_{1_0},x_{2_0}))$ the solutions of $D(\varepsilon,\alpha,m,T)$ and by $(\xi_1(t,\xi_{1_0},\xi_{2_0}),\xi_2(t, \xi_{1_0},\xi_{2_0}))$ the solutions of $\Sigma(\varepsilon,m,T)$. Let $(x_{1_0}, x_{2_0})$ be any initial condition for $D$ and choose $(\xi_{1_0},\xi_{2_0})$ such that:  $$x_{i_0}<\xi_{i_0}\quad i = 1,2$$
then, for every $t$ one has~:
$$x_i(t,x_{1_0},x_{2_0}) < \xi_i(t,\xi_{1_0},\xi_{2_0}) \quad \quad i = 1,2$$
Assume this is not the case ; let $t^*$ be the first time for which one has $x_i^*=x_i(t,x_{1_0},x_{2_0}) = \xi_i(t,\xi_{1_0},\xi_{2_0}) $  for at least one of the two indices ; assume for the shake of definitiveness that this index is 1~; one has:
$$\displaystyle \frac{dx_1(t^*)}{dt} =
  (±1-\varepsilon - m) x_1^*- \alpha x_1^{*^2} + mx_2(t^*) <  (±1-\varepsilon - m) x_1^* + m\xi_2(t^*)=
 \frac{d\xi_1(t^*)}{dt}$$
 The inequality $ \frac{dx_1(t^*)}{dt} <  \frac{d\xi_1(t^*)}{dt}$ contradicts the fact that $t^*$ is the first time for which $x_1(t,x_{1_0},x_{2_0}) = \xi_1(t,\xi_{1_0},\xi_{2_0}) $. Since $\Sigma(\varepsilon,m,T)$ is stable $ \xi_i(t,\xi_{1_0},\xi_{2_0})\;( i = 1,2)$ tends to 0 and also $x_i(t,x_{1_0},x_{2_0})\;(i=1,2) $.\\
 $\Box$
{\proposition \label{persistance2} If the parameters $(\eps,m,T)$ are such that the system $D(\varepsilon,0,m,T) =\Sigma(\varepsilon,m,T)$ is unstable, then the system $D(\varepsilon,\alpha,m,T)$ is uniformly persistant.}\\
In order to prove proposition \ref{persistance2} we need two lemmas. Let :
$$ U = \frac{1}{2} \ln(x_1x_2)\quad \quad V =\frac{1}{2} \ln(x_1/x_2)$$
In the $(U,V) $ variables the system $D$ is :
 \beq 
 D (\varepsilon,\alpha,m,T)\quad \quad\quad \quad
\left\{
\begin{array}{lcl}
\displaystyle  \frac{dU}{dt} &=&\displaystyle   \ch(V)-m-\varepsilon  -\alpha \emat^{U}\ch(2V) \\[8pt]
\displaystyle  \frac{dV}{dt} &=&\displaystyle  u(t)-m\,\sh(V)- \alpha \emat^{U}\sh(2V)
 \end{array} 
 \right. 
\feq
which is the system $S(\varepsilon,m,T)$  perturbed by the term :
\beq \label{perturbation}
\alpha \emat^{U} \left(
\begin{array}{c}
\ch(2V)\\
\sh(2V)
\end{array}
\right)
\feq
 It is easily seen that  the solutions of the  system $D$ enters in finite time the strip $\mathbb{R} \times [V^-_m,V^+_m]$ and thus persistence of $D$ is equivalent to the fact that for any solution  $\liminf_{t \rightarrow +\infty}U(t) > -\infty$. 
 
  Denote by :
  $$(U(t,U_0,V_0,t_0,S), F(t,U_0,V_0,t_0),S) \quad \quad(\mathrm{resp.} (U(t,U_0,V_0,t_0,D), F(t,U_0,V_0,t_0),D)) $$
 the solution of $S$ (resp. $D$) with initial condition $(U_0,V_0)$ at time $t_0$ and simply by $(U(t,D),V(t,D))$ (resp. $(U(t,S),V(t,S))$ the solution of $D$ (resp. $S$) when the reference to the initial condition is not needed.
 
 {\lemme \label{mino} Assume that  $S(\varepsilon,m,T)$ is unstable. Let $a>0$. Then there is $\theta > 0$   such that :
 $$  \forall\, V_0\in   [V^-_m,V^+_m],\forall\,U_0, \,\forall\,t_0 \;:\quad U(t_0+\theta,U_0,V_0,t_0,S) \geq U_0+a$$
  }
  \textbf{Proof}: Fix some $a > 0$. Since $S$ is unstable, for each $U_0,V_0,t_0$ such a $\theta$ exists~; it follows from the compactness of $[V^-_m,V^+_m]$, the periodicity of $S$ and the property $U(t,U_0,V_0,t_0) = U_0+U(t,0,V_0,t_0)$ that a universal $\theta$ does exist. $\Box$\\
 {\lemme \label{approximation}For any $\delta > 0$ there exists $\overline{U} $ such that :
  \beq
  \begin{array}{l}
   \left\{ \max_{t \leq t^*} U(t+t_0,U_0,V_0,t_0,D) \leq \overline{U} \right\}\Longrightarrow\cdots\\\cdots  | U(t^*+t_0,U_0,V_0,t_0,D)-U(t^*+t_0,U_0,V_0,t_0,S)|\leq t^*\delta
   \end{array}
   \feq}
 \textbf{Proof.} Since the perturbation  (\ref{perturbation}) tends to $0$ when $U$ tends to $-\infty$ uniformly with respect to $V \in [V^-_m,V^+_m]$ this is easily deduced from Gronwall inequality.\\\\
 \textbf{Proof of the proposition \ref{persistance2}} Fix some $a > 0$ and let $\theta$ be given by lemma \ref{mino} and $\overline{U}$ given by lemma \ref{approximation} such that  $\delta = \frac{a}{2\theta}$. The proof goes by contradiction.  Assume that: $$\liminf_{t \rightarrow +\infty}U(t,D) =  -\infty$$
then there exist (see Figure. \ref{ddep} ) $t_1$ and $t_2$ such that :
 \beq
 \overline{U} > U(t_1,D)+a > U(t_1,D) >  U(t_1,D)-\frac{\theta}{\pi} = U(t_2,D)
 \feq
 where $$-\pi = \min_{U\leq \overline{U}, V \in[V^-_m,V^+_m]}  m\,\ch(2V)-m-\varepsilon  -\alpha \emat^{U}\ch(V) < -\varepsilon$$
 Since $U(t,D)$ is continuous, from the intermediate value theorem there is some $\tau > t_1$ such that :
 \beq \label{eq1}
 t \in [\tau, t_2] \Rightarrow U(t,D) \leq U(t_1)
 \feq 
 and since $\pi$ is the minimum of the velocity of $U(t,D)$ it takes a duration $t_2-\tau$ greater than $\theta = \pi \frac{\theta}{\pi}$ to cover the distance from $U(\tau,D)$ to $U(t_2)$.
 \bitbul
 \item From lemma \ref{mino} :
 $$U(\tau+\theta, U(t_1,D),V(\tau,D),\tau,S) > U(t_1,D)+a$$
 (red curve of Figure \ref{ddep}).
 \item From lemma \ref{approximation}:
  $$|U(\tau+\theta, U(t_1,D),V(\tau,D)),\tau,S)-U(\tau+\theta, U(t_1,D),V(\tau,D),\tau,D)| <\frac{a}{2}$$
 \fit
 These points imply $\displaystyle U(\tau+\theta, U(t_1,D),V(\tau,D),\tau,D) \geq U(t_1)+\frac{a}{2}$ which is a contradiction with (\ref{eq1}).$\Box$
\begin{figure}
 \begin{center}
 \includegraphics[width=1\textwidth]{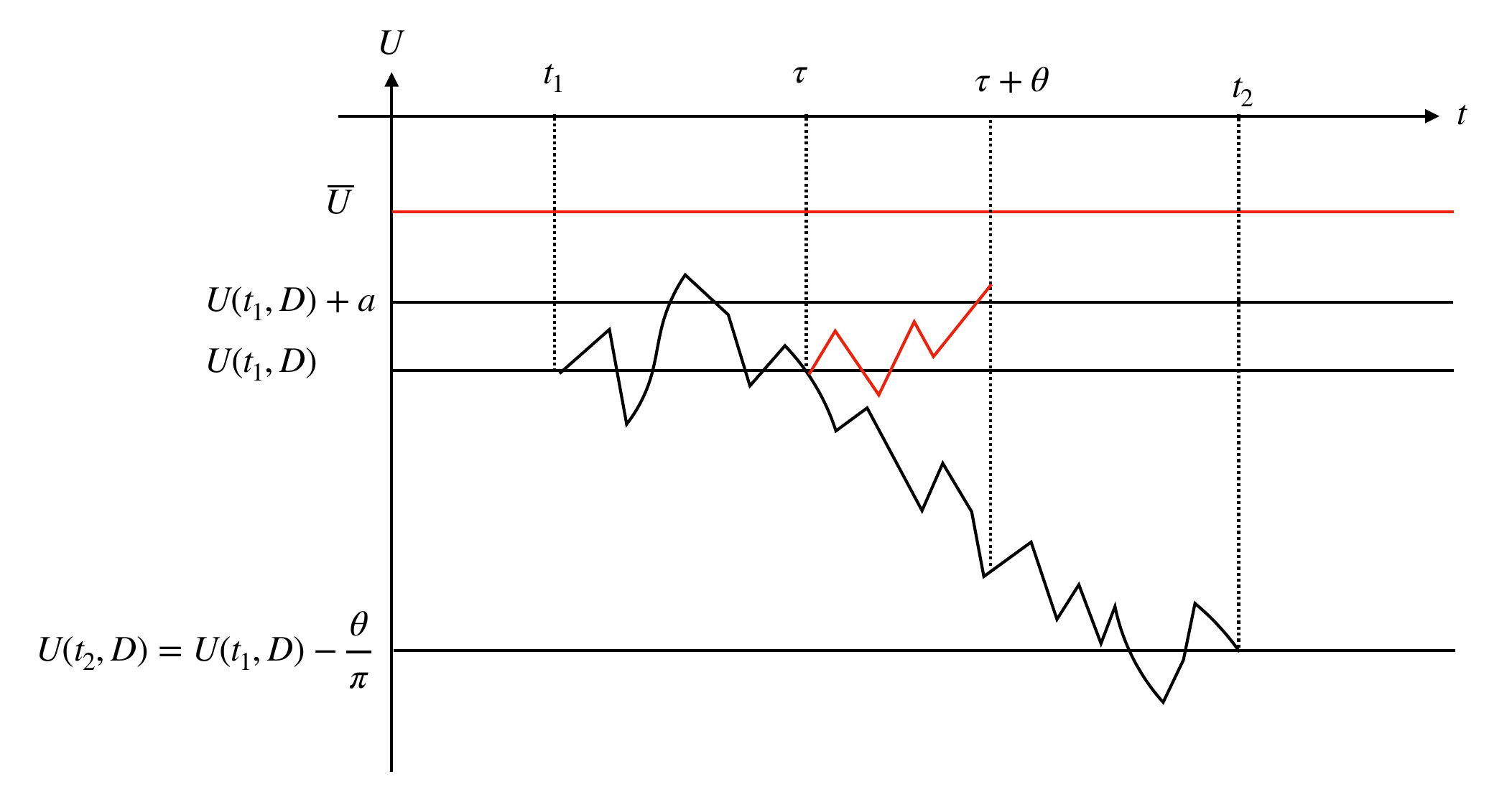}
 \caption{ Illustration of  the proof of  proposition \ref{persistance2} }\label{ddep}
 \end{center}
 \end{figure}

\subsection*{The stochastic case}
We now prove Proposition \ref{persistancesto} thanks to results in \cite{BS19}. Note that the vector fields in the right hand side of $\bm{D}(\varepsilon,\alpha,m,\sigma_-, \sigma_+)$ satisfy conditions E1, E2, E3, E4 and E5 in \cite[Section 4]{BS19} and admit a positively invariant compact set $K$ containing $0$. Thus, Proposition \ref{prop:lyapDelta},  \cite[Theorem 4.3]{BS19} (for the case $\Lambda(\varepsilon,m, T) <0$) and  \cite[Theorem 4.12]{BS19} (for the case $\Lambda(\varepsilon,m, T) > 0$) and \cite[Theorem 3.8]{NS20} (for the case $\Lambda(\varepsilon,m, T) = 0$) conclude the proof of proposition \ref{persistancesto}.

\section{Notations and glossary}


We give in this appendix the list of notations used in the paper
\bigskip

\begin{tabular}{ll}
$x_i$ ($i=1,2$)&
Abundance of population in patch $i$ ($i=1,2$)
\\
$U$, $V$&
$U=\ln\sqrt{x_1,x_2}$, $V=\ln\sqrt{x_1/x_2}$, see \eqref{varUV}
\\
$r_i(t)$ ($i=1,2$)&
Local growth of population in patch $i$ ($i=1,2$)
\\
$\bar{r}_i$ ($i=1,2$)&
Local average growth of population in patch $i$ ($i=1,2$)
\\
$\chi$&
$\chi=\frac{1}{2T}\int_0^{2T}\max(r_1(t),r_2(t))dt$
\\
$\Sigma(r_1,r_2,m,T)$&
Model of growth in two patches, see \eqref{BLSS1}
\\
$S(r_1,r_2,m,T)$&
System $\Sigma(r_1,r_2,m,T)$ in variables $(U,V)$, see \eqref{systvarUV}
\\
$\Sigma(\varepsilon,m,T)$&
Deterministic $(\pm 1)$-model, see \eqref{systeme1}\\
$\Sigma^+(m,T)$&
$(+ 1)$ system , see \eqref{Sigma+}
\\
$\Sigma^-(m,T)$&
$(- 1)$ system, see \eqref{Sigma-}
\\
$S(\varepsilon,m,T)$&
System $\Sigma(\varepsilon,m,T)$ in variables $(U,V)$, see \eqref{plusoumoins1enUV}
\\
$F(m,T)$&
Equation of $V$ in $S(\varepsilon,m,T)$, see \eqref{SV}
\\
$F^+_m$ and $F^-_m$&
Equation $F(\varepsilon,m,T)$ in environment ($+1$) and ($-1$), see \eqref{S1} and \eqref{S2}
\\
$V^+_m$ and $V^-_m$&
Equilibria of $F^+_m$ and $F^-_m$, respectively, see \eqref{Vm+Vm-}
\\
$P_{m,T}$ &
Periodic solution of $F(\varepsilon,m,T)$, see Proposition \ref{prop1}
\\
$\Delta(\eps,m,T)$&
Asymptotic growth of $U$, see \eqref{deltaT2} and Proposition \ref{prop2} 
 \\
$M^u_{\varepsilon,m}$, $u=\pm1$ &
Matrices of the linear systems   $\Sigma^+(m,T)$ and $\Sigma^-(m,T)$, see \eqref{Muepsm}
\\
$M(\varepsilon,m,T)$&
Period mapping of
$\Sigma(\varepsilon,m,T)$, see 
\eqref{MepsmT}
 \\
 $\sigma(\varepsilon,m,T)$&
 Spectral radius of $M(\varepsilon,m,T)$, see \eqref{rayonspectral}
 \\
$\bm{\Sigma}(\varepsilon,m,\mu)$&
Stochastic $(\pm 1)$-model (random switching times),
see \eqref{systeme1random}
 \\
$\bm{S}(\varepsilon,m,\mu)$&
System $\bm{\Sigma}(\varepsilon,m,\mu)$ in variables $(U,V)$,
see \eqref{Srandom}
 \\
$\bm{F}(m,\mu)$&
Equation of $V$ in $\bm{S}(\varepsilon,m,\mu)$,
see \eqref{eq:VR}
\\
$\bm{\Delta}( \varepsilon, m, \mu)$&
Asymptotic growth rate of $U$, see  \eqref{eq:Delta_random_time} and Proposition \ref{prop:Delta_random_time}
\\
PDMP&
Piecewise Deterministic Markov Processes, see Section \ref{sec:PDMP}
\\
$\bm{\Sigma}(\varepsilon,p_1,p_2,m,T)$
&
Stochastic $(\pm 1)$-model (random choices of $\pm 1$),
see \eqref{system1RandomChoices}
\\
$
\bm{\Delta}(p_1, p_2, m, T)$&
Asymptotic growth rate of $U$, see Proposition \ref{prop:Delta_random_choice}
\\
$\lambda_1$&
Top Lyapunov exponent, see Section \ref{sec:TLE}
\\
$\Sigma^{\pm\pm}(\varepsilon,m,T,\varphi)$&
Periodic $(\pm 1)$ model with phase shift, see \eqref{pm}, 
\eqref{mp}, \eqref{pp} and \eqref{mm}
\\
$D(\varepsilon,\alpha,m,T)$&
Density dependent deterministic $(\pm 1)$-model, see \eqref{logistic1}.
\\
$\bm{D}(\varepsilon,\alpha,m,T)$&
Density dependent stochastic $(\pm 1)$-model, see \eqref{logistic1r}.
\end{tabular}


\fin